\RequirePackage{fix-cm}
\documentclass[smallextended]{svjour3}       
\smartqed  

\usepackage[utf8]{inputenc}
\usepackage[dvipsnames]{xcolor}
\usepackage{graphicx}
\usepackage{times}
\usepackage[utf8]{inputenc}
\usepackage{mathtools}
\usepackage{amsfonts}
\usepackage[T1]{fontenc}
\usepackage{bm}
\usepackage{verbatim}
\usepackage{algorithm}
\usepackage{framed}
\usepackage{booktabs,array}
\usepackage{xcolor}
\usepackage{comment}
\usepackage{marginnote}
\usepackage{mathtools}
\usepackage{amssymb}
\usepackage{subfigure}
\usepackage{yhmath}

\usepackage{ulem}
\usepackage{cancel}
\usepackage{soul,xcolor}

\usepackage{tikz}
\usetikzlibrary{decorations.pathreplacing, decorations.pathmorphing, decorations.shapes}

\newcommand{\statevecT}{\underline{\rm \bf T}}
\newcommand{\horizon}{\delta}
\newcommand{\R}{\mathbb{R}}

\newcommand{\mcS}{\mathcal{S}}

\newcommand{\mcC}{\mathcal{C}}
\newcommand{\mcD}{\mathcal{D}}

\newcommand{\mcW}{\mathcal{W}}

\newcommand{\mcL}{\mathcal{L}}

\newcommand{\mcX}{\mathcal{X}}
\newcommand{\mcO}{\mathcal{O}}

\newcommand{\mbR}{\mathbb{R}}

\newcommand{\mbRn}{{\mathbb{R}^n}}

\newcommand{\Rn}{{\mathbb{R}^n}}

\newcommand{\abs}[1]{\left| #1 \right|}
\newcommand{\omg}{{\Omega}}
\newcommand{\omgb}{\mathcal{B}\Omega}
\newcommand{\ooomg}{\wideparen\Omega}
\newcommand{\oomg}{\overline{\omg}}
\newcommand{\omgl}{{\Omega_l}}
\newcommand{\omgn}{{\Omega_{nl}}}
\newcommand{\omgo}{{\Omega_o}}
\newcommand{\omgv}{{\Omega_v}}
\newcommand{\omgt}{{\Omega_t}}
\newcommand{\omgd}{{\Omega_p}}
\newcommand{\gammav}{{\Gamma_v}}
\newcommand{\gammad}{{\Gamma_p}}
\newcommand{\bb}{\mathbf{b}}
\newcommand{\alphab}{{\boldsymbol\alpha}}

\newcommand{\nub}{{\boldsymbol \nu}}

\newcommand{\phib}{{\boldsymbol \phi}}
\newcommand{\psib}{{\boldsymbol \psi}}

\newcommand{\epsb}{{\boldsymbol \varepsilon}}
\newcommand{\fb}{\mathbf{f}}
\newcommand{\sbp}{\mathbf{s}}
\newcommand{\gb}{\mathbf{g}}
\newcommand{\ub}{\mathbf{u}}

\newcommand{\vb}{\mathbf{v}}
\newcommand{\xb}{\mathbf{x}}
\newcommand{\xbp}{\mathbf{x}'}

\newcommand{\yb}{\mathbf{y}}

\newcommand{\zerob}{\mathbf{0}}
\newcommand{\xib}{{\boldsymbol\xi}}
\newcommand{\etab}{{\boldsymbol\eta}}

\newcommand{\Psib}{{\boldsymbol \Psi}}

\begin{document}

\title{A review of Local-to-Nonlocal coupling methods in nonlocal diffusion and nonlocal mechanics\thanks{M. D’Elia was supported by Sandia National Laboratories (SNL), SNL is a multimission laboratory managed and operated by National Technology and Engineering Solutions of Sandia, LLC., a wholly owned subsidiary of Honeywell International, Inc., for the U.S. Department of Energy's National Nuclear Security Administration contract number DE-NA0003525. This material is based upon work supported by the U.S. Department of Energy, Office of Science, Office of Advanced Scientific Computing Research under Award Number DE-SC-0000230927. This paper describes objective technical results and analysis. Any subjective views or opinions that might be expressed in the paper do not necessarily represent the views of the U.S. Department of Energy or the United States Government.\\
X. Li was supported by NSF - DMS 1720245 and a UNC Charlotte faculty research grant.\\
P. Seleson was supported by the Laboratory Directed Research and Development Program of Oak Ridge National Laboratory, managed by UT-Battelle, LLC, for the U. S. Department of Energy.
This manuscript has been co-authored by UT-Battelle, LLC under Contract No. DE-AC05-00OR22725 with the U.S. Department of Energy. The United States Government retains
and the publisher, by accepting the article for publication, acknowledges that the United
States Government retains a non-exclusive, paid-up, irrevocable, world-wide license to publish or reproduce the published form of this manuscript, or allow others to do so, for United States Government purposes. The Department of Energy will provide public access to these results of federally sponsored research in accordance with the DOE Public Access Plan (http://energy.gov/downloads/doe-public-access-plan).\\ 
X. Tian was supported in part by NSF grant DMS-1819233.\\
Y. Yu was supported by NSF - DMS 1620434 and the Lehigh faculty research grant.
}
}
%
\titlerunning{
Review of LtN coupling in nonlocal diffusion and  mechanics
}        

\author{Marta D'Elia \and Xingjie Li \and Pablo Seleson \and Xiaochuan Tian \and Yue Yu}

\authorrunning{D'Elia, Li, Seleson, Tian, and Yu}

\institute{
M. D'Elia \at Center for Computing Research, Sandia National Laboratories, Albuquerque, NM, 87123\\
\email{mdelia@sandia.gov}
\and X. Li \at
Department of Mathematics and Statistics, University of North Carolina-Charlotte, Charlotte, NC, 28223\\
\email{xli47@uncc.edu}
\and P. Seleson \at
Computer Science and Mathematics Division,
Oak Ridge National Laboratory, Oak Ridge, TN, 37831\\
\email{selesonpd@ornl.gov}
\and X. Tian \at
Department of Mathematics, The University of Texas at Austin, Austin, TX, 78712\\
\email{xtian@math.utexas.edu}
\and Y. Yu \at
Department of Mathematics, Lehigh University, Bethlehem, PA, 07922\\
\email{yuy214@lehigh.edu}
}

\date{Received: date / Accepted: }

\maketitle

\begin{abstract}
Local-to-Nonlocal (LtN) coupling refers to a class of methods aimed at combining nonlocal and local modeling descriptions of a given system into a unified coupled representation. This allows to consolidate the accuracy of nonlocal models with the computational expediency of their local counterparts, while often simultaneously removing additional nonlocal modeling issues such as surface effects. The number and variety of proposed LtN coupling approaches have significantly grown in recent year, yet the field of LtN coupling continues to grow and still has open challenges. 
%
This review provides an overview of the state-of-the-art of LtN coupling in the context of nonlocal diffusion and nonlocal mechanics, specifically peridynamics. We present a classification of LtN coupling methods and discuss common features and challenges. The goal of this review is not to provide a preferred way to address LtN coupling but to present a broad perspective of the field, which would serve as guidance for practitioners in the   selection of appropriate LtN coupling methods based on the characteristics and needs of the problem under consideration. 
\end{abstract}

\section{Introduction}\label{sec:intro}

\subsection{Nonlocal models and the need of coupling methods}\label{sec:nonlocal-models}

Nonlocal models such as nonlocal diffusion and peridynamics can describe phenomena not well represented by classical Partial Differential Equations (PDEs). These include problems characterized by long-range interactions and discontinuities~\cite[Chapter 1]{DElia2015handbook}. For instance, in the context of diffusion, long-range interactions effectively describe anomalous diffusion, whereas in the context of mechanics, cracks formation results in material discontinuities. We refer to these phenomena, in a general sense, as \textit{nonlocal effects}. 
Even though this work is mostly focused on nonlocal models for diffusion and mechanics applications, nonlocal models can characterize a wide range of scientific and engineering problems, including 
subsurface transport 
\cite{Benson2000,katiyar2019general,katiyar2014peridynamic,Schumer2003,Schumer2001},
phase transitions \cite{Bates1999,Chen_nonlocalmodels,dayal2006kinetics},
image processing \cite{Buades2010,Gilboa2007,Lou2010},
multiscale and multiphysics systems \cite{Alali2012,Askari2008,DET19,seleson2010peridynamic}, 
magnetohydrodynamic turbulence \cite{Schekochihin2008}, 
and stochastic processes \cite{Burch2014,DElia2017,Meerschaert2012,MeKl00}.

The fundamental difference between nonlocal models and classical local PDE-based models is the fact that the latter only involve differential operators, whereas the former also rely on integral operators.\footnote{
Nonlocal models can be based on \textit{weakly nonlocal} or \textit{strongly nonlocal} formulations. The former enrich classical PDEs by explicitly including higher-gradients of field variables, whereas the latter is based on integral-type formulations with a weighted  averaging of field variables \cite{bavzant2002nonlocal,Paola2013Review,jirasek2004nonlocal}. 
In this paper, we only concern ourselves with strongly nonlocal formulations.}
The integral form allows for the description of long-range interactions (spanning either small regions or the whole space) and reduces the regularity requirements on problem solutions. Here, we consider nonlocal models, based on integro-differential formulations, characterized by integral operators in space that lack spatial derivatives. This enhances the accuracy of their modeling representations by generalizing the space of admissible solutions, which can feature singularities and discontinuities.

Despite their improved accuracy, the usability of nonlocal equations could be compromised by several modeling and numerical challenges such as the unconventional prescription of nonlocal boundary conditions, the calibration of nonlocal model parameters, often unknown or subject to uncertainty, 
and the expensive numerical solution. In fact, the computational cost of solving a nonlocal problem is significantly higher than that corresponding to PDEs. Specifically, the associated computational expense of a nonlocal problem depends on the ratio between the characteristic nonlocal interaction length (the so-called interaction radius or {\it horizon}) and the chosen simulation grid or mesh size. When this ratio becomes large, simulations can be unfeasible \cite[Chapter 14]{DElia2015handbook}.

Nevertheless, it is often the case that nonlocal effects are concentrated only in some parts of the domain, whereas, in the remaining parts, the system can be accurately described by a PDE. 
The goal of Local-to-Nonlocal (LtN) coupling is to combine the computational efficiency of PDEs with the accuracy of nonlocal models, under the assumption that the location of nonlocal effects can be identified. In this context, the main challenge of a coupling method is to accurately merge substantially different local and nonlocal descriptions of a single system into a mathematically and physically consistent coupled formulation.

As an added value, LtN coupling can provide a viable way to circumvent the non-trivial task of prescribing nonlocal boundary conditions. Such conditions have to be prescribed in a layer surrounding the domain where data are not available or are hard to access; however, surface (local) data are normally available. 
In practice, this requires extending surface (local) boundary conditions to volumetric (nonlocal) boundary conditions in a way that is not always clear or well-defined. 
An \textit{ad hoc} treatment of nonlocal boundaries often results in unphysical surface effects (see dicussions in \cite{le2018surface}). 
Using a local model adjacent to the boundary of the domain allows for prescription of classical boundary conditions, provided the solution is regular enough. Throughout this paper, we mention coupling approaches that have been used for this task. There are other potential benefits of LtN coupling, such as controlling undesired wave dispersion or leveraging available computational tools based on classical PDEs. However, these are beyond the scope of this paper.

The purpose of this paper is to present the state-of-the-art of LtN coupling in the context of nonlocal diffusion and nonlocal mechanics, specifically peridynamics. While the list of proposed LtN coupling methods is extensive, this paper focuses on a select number of approaches. The description provided should, however, give the reader a broad enough perspective on the diversity of LtN coupling techniques and the variety of ways to approach LtN coupling problems, both in terms of mathematical formulation and practical implementation. 
We stress that our goal is strictly to give an overview of available coupling strategies and highlight their properties, while describing common features and challenges. This work does not intend to present a relative assessment of LtN coupling methods. Yet the reader can use this review as a guide for selecting the most appropriate method for the problem at hand.

\subsection{Overview of classes of coupling methods}\label{sec: overview of methods}
\begin{figure}[t]
\centering
\includegraphics[width=1\linewidth]{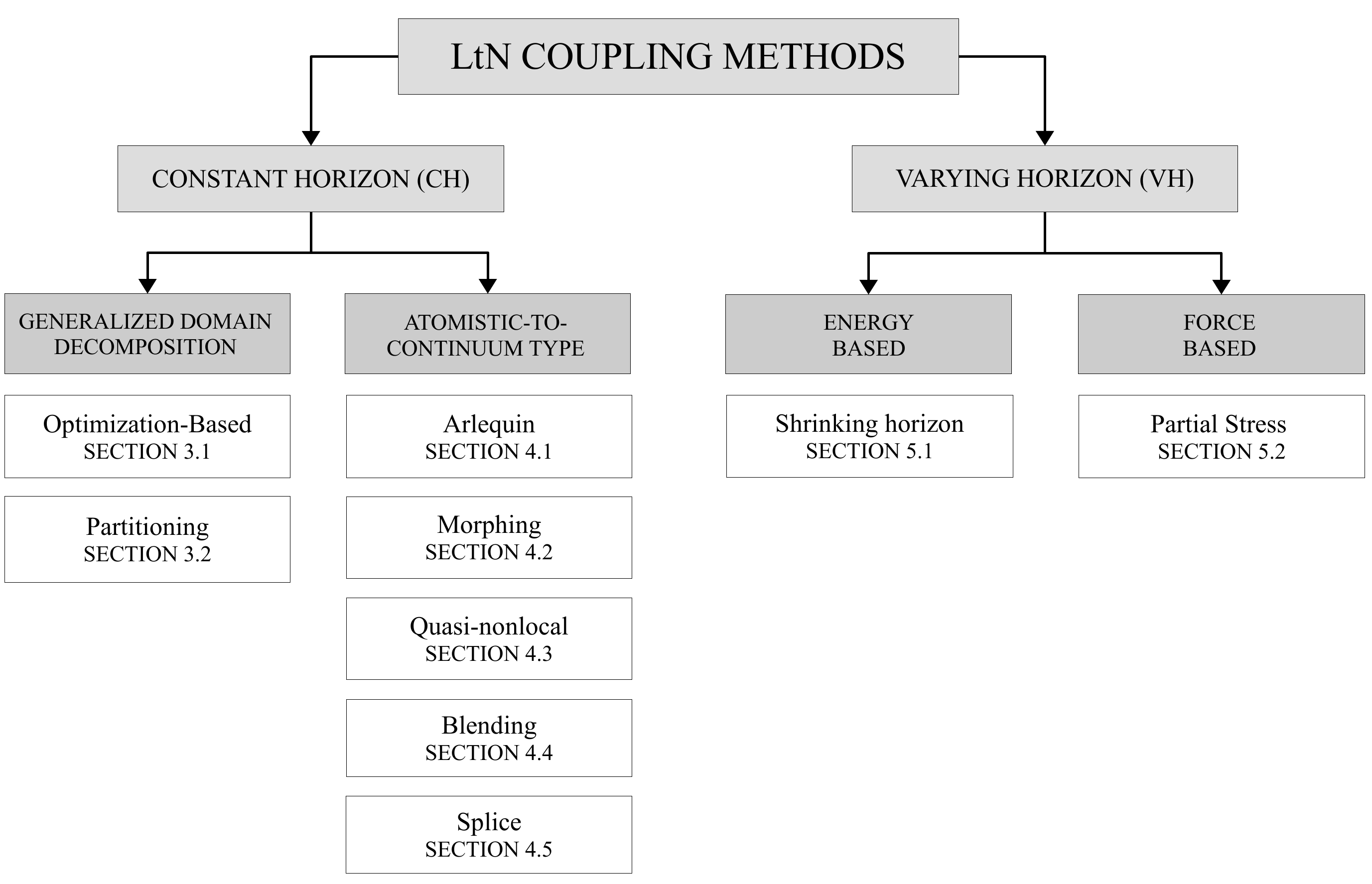}
\vspace{-.5cm}
\caption{Overview of classes of LtN coupling approaches and corresponding methods.}
\label{fig:overview-chart}
\end{figure}

We divide coupling approaches in two main classes based on how the transition from a nonlocal description to a local description is carried out. A schematic overview of such classification and of the methods belonging to each class is reported in Figure~\ref{fig:overview-chart}; that figure also indicates the corresponding sections where the methods are described in detail. We refer to the first class as the {\it Constant Horizon} (CH) class, which contains approaches characterized by an abrupt change in the horizon as we move from the nonlocal region to the local region; whereas we refer to the second class as the {\it Varying Horizon} (VH) class, which contains approaches that, in contrast, feature a smooth transition of the horizon. In the next paragraphs we briefly mention several methods belonging to each class. A detailed description of those methods will be reported in the next sections.

In the CH class, we first identify a group of approaches that resemble 
generalized domain decomposition (GDD) methods. A first example is the optimization-based coupling method introduced in \cite{Bochev_14_INPROC}, analyzed in \cite{Delia2019}, and extended to three-dimensional settings in \cite{Bochev_16b_CAMWA}. This strategy treats the coupling condition as an optimization objective, which is minimized subject to the model equations acting independently in their respective sub-domains. As opposed to other LtN coupling approaches, this method reverses the roles of the coupling conditions and the governing equations, keeping the latter separate. In particular, the coupling of local and nonlocal models is effected by couching the LtN coupling into an optimization problem. The objective is to minimize the mismatch of the local and nonlocal solutions on the overlap of their sub-domains, the constraints are the associated governing equations, and the controls are the nonlocal volume constraints and the local boundary conditions on the virtual boundaries generated by the decomposition.

A second GDD example relies on the partitioned procedure as a general coupling strategy for heterogeneous systems, such as multiscale and multiphysics problems \cite{badia2008fluid,mathew2008domain,quarteroni1999domain,toselli2006domain}.
%
%
In the partitioned procedure, the system is divided into sub-problems in their respective sub-domains. Different models are then employed independently in each sub-problem,  which communicates with other sub-problems only via transmission conditions on the sub-domain interfaces. 
The coupled problem is then solved based on sequential 
%
solutions of sub-problems, and proper transmission conditions are required on the sub-domain interfaces to impose solution continuity on those interfaces and to enforce the energy balance of the whole system. Among the possible transmission conditions, the Robin transmission condition, which is a linear combination of the Dirichlet and Neumann transmission conditions, has been proven to be very efficient (see, e.g., \cite{badia2008fluid,chen2011parallel,discacciati2007robin,douglas1997accelerated}).
%
%
In \cite{yu2018partitioned}, a partitioned procedure with Robin transmission conditions was applied to LtN coupling of mechanics models with an overlapping region, and the method was later extended to LtN coupling without an overlapping region in \cite{you2019coupling}.

A second group of approaches in the CH class is also based on a decomposition of the domain, which resembles an overlapping sub-domain decomposition. 
However, as opposed to GDD methods, these approaches rely on hybrid 
descriptions that combine local and nonlocal models in a transition region between local and nonlocal sub-domains. 
We refer to this group of approaches as Atomistic-to-Continuum (AtC) type coupling approaches due to their resemblance to such coupling methods (see, e.g., the review
articles \cite{Luskin2013a,Tadmor2009a}). We observe that, in some sense, AtC coupling is a special case of LtN coupling, where the nonlocal model is a discrete atomistic representation and the local model is given by a classical PDE.

In the group of AtC type coupling approaches, a first example is the Arlequin method. This is a general coupling technique
introduced in \cite{dhia1998multiscale,ArlequinDhia1999} and further studied in \cite{ArlequinDhia2001,ArlequinDhia2005}. 
This technique 
was applied to AtC coupling in \cite{Bauman2008} with various follow-on works (see, e.g., \cite{Chamoin2010,Prudhomme2008,Prudhomme2009}). 
Application of the Arlequin method
for LtN coupling was proposed in \cite{HanLubineau2012} in the context of static problems and later applied to dynamic settings in \cite{ArlequinWang2019}.
The Arlequin method 
is an example of an energy-based blending approach, where the energy of the system in the transition 
region is defined as a weighted average of the local and nonlocal energies. 
In the local and nonlocal sub-domains, 
the energy is defined according to the models operating in those sub-domains. A Lagrange multiplier enforces compatibility of the kinematics of both models.

A second example of AtC type coupling is the morphing approach proposed in \cite{lubineau2012morphing} and then extended in \cite{azdoud2013morphing,azdoud2014morphing,han2016morphing}. This method is based on blending, or morphing, the material properties of the local and nonlocal models. The method consists of a single model defined over the entire domain with an equilibrium equation that contains both local and nonlocal contributions. The transition between the nonlocal sub-domain to the local sub-domain is achieved through a gradual change in the material
properties characterizing the two models in the transition region or   ``morphing'' zone.
In this region, local and nonlocal material properties are suitably weighted under the constraint of energy equivalence for homogeneous deformations. 

A third example of AtC type coupling is the quasi-nonlocal (QNL) coupling method. This method was originally proposed in the context of AtC coupling and is based on energy minimization with the aim of eliminating linear spurious forces. The QNL method redefines 
the nonlocal energy of the system, where nonlocal interactions contained within the local sub-domain are reformulated. 
The application of the QNL coupling method to LtN coupling appears in \cite{DuLiLuTian2018,XHLiLu2017}. 

A fourth example of AtC type coupling is the force-based blending method proposed in \cite{Seleson2013CMS} and extended in \cite{seleson2015concurrent}. Force-based blending has been studied in the context of AtC coupling (see, e.g., \cite{badia2007force,Fish2007,li2012positive,li2014theory}). 
This approach employs a blending function to create a weighted average of the local and nonlocal governing equations, and that function is chosen in a way that the blended model reduces to the local and nonlocal models in their respective sub-domains. 
Similarly to the morphing approach, a single blended model is defined over the entire domain; however, the force-based blending method does not enforce energy equivalence. 
In this paper, we employ the term ``blending'' to refer to  force-based blending, since other methods which could be classified within a blending category, such as the Arlequin and morphing, are called by their specific names. As opposed to AtC blending methods, which normally seek means to blend given atomistic and continuum models, the peculiarity of the LtN blending method from \cite{Seleson2013CMS,seleson2015concurrent} is that a reference nonlocal model is first postulated over the entire domain, and then the blended model is attained by simply combining the use of a blending function with assumptions on the material response. An underlying connection between the local and nonlocal models is leveraged in the derivation of the blended model.

A final example of AtC type coupling is the splice method. This method was first proposed in~\cite{silling2015variable} to couple two nonlocal models with different horizons and then applied to LtN coupling. A similar approach for LtN coupling, although introduced instead for discretized models, was presented in a series of publications (see \cite{galvanetto2016effective,ni2019coupling,shojaei2016coupled,shojaei2016coupling,shojaei2017coupling,zaccariotto2018coupling,zaccariotto2017enhanced}). In the splice method, the governing equation at each point is given by either the local or nonlocal model. There is no particular coupling enforced, except that points described by the nonlocal model may interact with points described by the local model and vice versa. When such a situation occurs, points with a given description (either local or nonlocal) treat its environment with the same modeling representation as its own. For instance, a point in a nonlocal region interacting with some points in a local region would treat those points as if they were also described by a nonlocal model.

In the VH class, we can also find energy-based and force-based approaches. The idea behind VH is based on the fact that a nonlocal model converges to a local model, under suitable regularity assumptions, when the horizon approaches zero (see Sections~\ref{subsubsec:diffusion} and~\ref{subsubsec:mechanics}). Naturally, allowing the horizon to vary spatially in a domain, so that it approaches zero in certain sub-domains with enough regularity, provides a transition from a nonlocal to a local representation. The first discussions on a variable horizon in peridynamics appeared, in the context of adaptive refinement, in~\cite{bobaru2011adaptive,bobaru2009convergence}. 
Similar adaptivity ideas actually appeared in the context of atomistic systems in~\cite{seleson2010bridging} to remove surface effects. 
A nonlocal diffusion formulation with spatially varying horizon applied to interface problems was presented in~\cite{seleson2013interface}, where two-horizon systems were treated and specialized to the case in which one of the two horizons is taken to zero. A similar formulation applied to two-horizon systems in peridynamics was discussed in \cite{ren2016dual,seleson2010peridynamic}. 
The case of a smoothly varying horizon was presented and analyzed in~\cite{silling2015variable}, and it was used  to couple two peridynamic models with different horizons. In \cite{TTD19,TD17trace}, the validity of nonlocal diffusion models with a shrinking horizon applied to LtN coupling was discussed.

The first group of energy-based approaches in the VH class is related to the formulations presented in~\cite{ren2016dual,seleson2013interface,seleson2010peridynamic,silling2015variable,TTD19,TD17trace}. Such energy-based formulations, unfortunately, introduce linear spurious forces. 
In order to allow varying the horizon in space while 
preventing this coupling artifact, the  partial stress method was presented in~\cite{silling2015variable}. This method belongs to the group of force-based approaches and introduces a new tensor field referred to as the partial stress, which is used to describe the material response in the transition region between the nonlocal and local sub-domains.

There are other approaches in the literature concerning LtN coupling. Specifically, while we focus on LtN coupling of continuum models, some proposed methods address the coupling at the discrete level (we actually mentioned a few above). In fact, the first peridynamic work in this context appear in \cite{macek2007peridynamics}, where a coupling is performed by implementing a peridynamic model in a conventional finite element analysis code using truss elements. 
In this context, an overlapping-based approach between peridynamics and classical finite elements, which relies on interfacial elements, is described in \cite{kilic2010coupling,liu2012coupling}, whereas a sub-modeling approach is presented in \cite{oterkus2012combined}. 
A force-based coupling approach, which resembles blending due to a partition of unity of local and nonlocal displacements within an overlapping region, is presented in~\cite{sun2019superposition}. 
%
%
A different type of coupling approaches have been published specifically addressing dynamic coupling artifacts, such as interfacial spurious wave reflections (see, e.g., \cite{kulkarni2018analytical,nicely2018nonlocal}). 
Finally, a coupling based on coarsening a peridynamic discretizations was proposed in \cite{seleson2010peridynamic} and later implemented in \cite{xu2016multiscale2,xu2016multiscale}; these implementations resemble the splice method in the context of peridynamic finite elements. 
All these approaches, however, are beyond the discussions held in this review.
%

\paragraph{Desired properties of a LtN coupling method.} 
There are various considerations to account for when designing a LtN coupling method. The first one is ensuring that the reference local and nonlocal models are \textit{physically consistent}. This is guaranteed by the convergence of the nonlocal model to the local model as discussed in Sections~\ref{subsubsec:diffusion} and~\ref{subsubsec:mechanics} for nonlocal diffusion and nonlocal mechanics, respectively. 
Once the local and nonlocal models are determined to be physically consistent, the main challenge becomes removing or minimizing the appearance of coupling artifacts on the interface/in the transition 
region between the local and nonlocal sub-domains. We describe below some desired properties of LtN coupling methods and discuss related coupling artifacts.

A basic desired property of a LtN coupling method is \textit{patch-test consistency}, which is established when such a method passes the so-called patch test or consistency test \cite{SandiaReport2015,Luskin2013a,Tadmor2009a}. The main idea is as follows: if for a certain class of problems the local and nonlocal solutions, $u^l$ and $u^{nl}$, respectively, coincide, then ``patching'' the two problems by coupling the corresponding models should still return the same problem solution. Note that to be comparable, the local and nonlocal problems are augmented with consistent boundary conditions and forcing terms. 

As an example, we define a linear patch test in one dimension as follows: 

\begin{definition} [Linear patch test]\label{def:patch-test}
Given a linear function $u^{\rm lin}(x)=a_0+a_1x$ with $a_0$ and $a_1$ constants and local and nonlocal operators $\mcL^{\rm L}$ and $\mcL^{\rm{NL}}$, respectively, such that
$$ 
\mcL^{\rm L} u^{\rm lin}=0 \quad {\rm and} \quad \mcL^{\rm NL} u^{\rm lin}=0,
$$
a coupling method passes the linear patch test if, in the absence of forcing terms and with consistent boundary conditions, $u^{\rm lin}$ is also the solution of the coupled problem.
\end{definition}
Similarly to Definition \ref{def:patch-test}, one can define a higher-order patch test when, given the $p$th degree polynomial $u^{\rm poly}(x)=\sum_{i=0}^p a_i x^i$ with $\{a_i\}_{i=0}^p$ constants, $\mcL^{\rm L} u^{\rm poly}=\mcL^{\rm NL} u^{\rm poly}$. Examples of linear ($p=1$) and quadratic ($p=2$) patch tests are provided in the following sections. 

A related concept to the patch test is that of ``ghost'' forces. Specifically, when a coupling method does not pass the linear patch test, non-physical forces or fluxes normally arise in the transition region between the local and nonlocal sub-domains.
These non-physical forces are often referred to as ``ghost'' forces \cite{Tadmor2009a}.


Another desired property of a LtN coupling method is \textit{asymptotic compatibility}. This ensures that the method preserves the physical consistency of the local and nonlocal models. Specifically, the solution corresponding to a LtN coupling method should be such that it coincides with the local solution everywhere when the nonlocal effects vanish. More details in this regard are provided in the following sections. 

A third desired property of a LtN coupling method is \textit{energy equivalence}. Due to the physical consistency of the local and nonlocal models, it is expected that these models would have equivalent energy descriptions for a certain class of problems. Consequently, an appealing property of a LtN coupling method is to also preserve the energy description for such class of problems. 
This consideration becomes a natural requirement for energy-based LtN coupling methods, since they possess an associated energy functional. In contrast, force-based LtN coupling methods do not possess, in general, an associated energy.

While the above properties concern common coupling artifacts in static problems, other spurious effects can emerge in dynamic scenarios. A prime example of these spurious effects is wave reflections on the interface/in the transition region between the local and nonlocal sub-domains. While most LtN coupling methods retain at least some of the above static properties to a certain degree, controlling spurious wave reflections appears to be a much harder task and thus related discussions are excluded from this review.

We point out that the majority of the methods described in this review are not tied to a particular discretization. However, the choice of discretization used for local and nonlocal models may affect some of the properties mentioned above. As an example, when the nonlocal discretization does not guarantee asymptotic compatibility, the coupling method will inherit that limitation.
 
\subsection{Outline of the paper}
This review is organized as follows. In Section \ref{sec:models-results}, we introduce the notation, discuss coupling configurations, and describe the 
 mathematical models and recall relevant results, including a description of the nonlocal vector calculus \cite{Du2013}. 
 In Section \ref{sec:GDD}, we provide a description of GDD approaches, i.e., Optimization-Based methods and the Partitioned procedure. 
 Section \ref{sec:AtC} describes several methods belonging to the group of AtC type coupling approaches, including the Arlequin, Morphing, Quasi-Nonlocal, Blending, and Splice methods. In Section \ref{sec:VH}, we present two VH approaches: the energy-based Shrinking Horizon method and the force-based Partial Stress method. 
 %
 For each method, we provide a mathematical formulation and describe its properties. We also report relevant numerical results based on available literature with the purpose of illustrating theoretical properties and showing applicability to realistic settings.
In Section \ref{sec:conclusion}, we draw conclusions and present guidelines for an appropriate choice of LtN coupling methods based on necessities and constraints. The chart in Figure \ref{fig:overview-chart} summarizes the classes of LtN coupling approaches and corresponding methods. 

\section{
Notation, coupling configurations, and  
 mathematical models
}\label{sec:models-results}
In this section, we introduce the notation and coupling configurations used in this paper (see Section~\ref{sec: notation and configurations}) and describe the 
nonlocal models and recall relevant results (see  Section~\ref{sec: nonlocal models})  that will be useful throughout the following sections.

\subsection{Nonlocal variables and nonlocal domains} 
\label{sec: notation and configurations}
Let $\omg\in\Rn$, $n=1,2,3$, be a bounded open domain. We are interested in functions $u:\omg\to\mbR$ and $\ub:\Omega\to\Rn$, $n=1,2,3$, solutions of nonlocal diffusion and nonlocal mechanics problems. 
%
Specifically, in diffusion, $u$ represents the concentration of a diffusive quantity and, in mechanics, $\ub$ represents the displacement in $n$ dimensions.

\begin{figure}[H]
\centering
\includegraphics[width=0.41\linewidth]{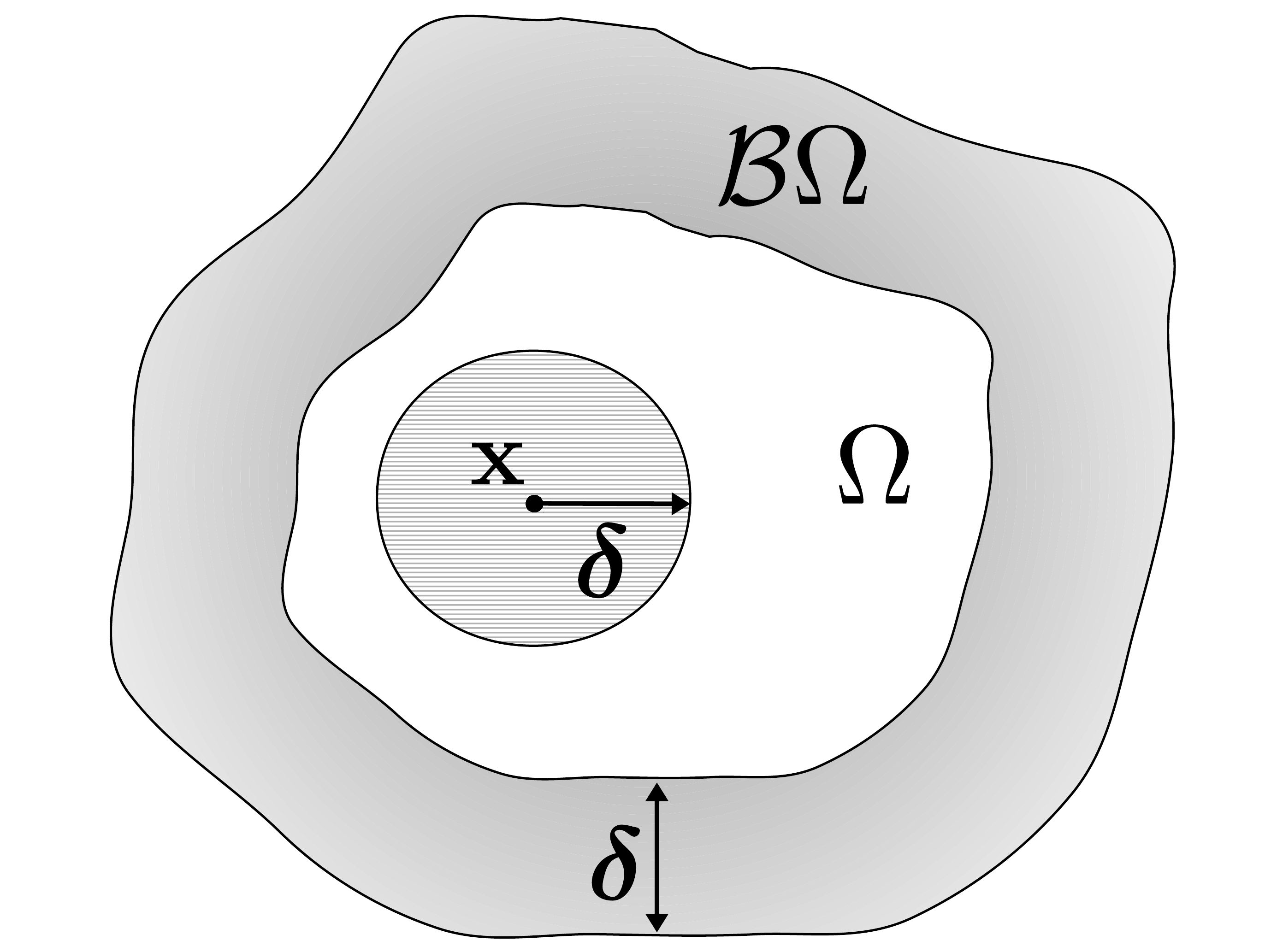}
\caption{Domain $\omg$, nonlocal boundary $\omgb$, and nonlocal neighborhood $B_\horizon(\xb)$. 
}
\label{fig:standard-domain}
\setlength{\unitlength}{1cm}
\begin{picture}(10,6)
\put(4.3,8.7){\tiny $\boldsymbol{B_\horizon(\xb)}$}
\end{picture}
\vspace*{-2in}
\end{figure}

In nonlocal settings, every point in a domain interacts with a neighborhood of points. Usually, such neighborhood is an Euclidean ball surrounding points in the domain, i.e.,
\begin{equation}\label{eq:Euclidean-ball}
B_\horizon(\xb)=\{\xbp\in\mbRn: \|\xbp - \xb\|\leqslant\horizon\},
\end{equation}
where $\horizon$ is referred to as the {\it horizon}.  
This fact has implications on the concept of boundary conditions that are no longer prescribed on $\partial\omg$, but on a collar of thickness of at least $\horizon$ surrounding the domain that we refer to as the {\it nonlocal volumetric boundary domain}, $\omgb$, or simply {\it nonlocal boundary}.\footnote{In many nonlocal models, such as bond-based peridynamic models, a collar of thickness $\horizon$ is normally sufficient as nonlocal boundary. However, more general nonlocal models, such as state-based perdiynamic models, may require a  collar of thickness $2\horizon$ (see, e.g., \cite{seleson2016convergence}).} This set, by definition, consists of all points outside the domain that interact with points inside the domain. Here, prescription of volume constraints guarantees the well-posedness of the problem \cite{Du2013}. Figure~\ref{fig:standard-domain} provides a two-dimensional configuration. We denote the union of $\oomg$ and $\overline{\omgb}$ by $\ooomg:=\oomg\cup\overline{\omgb}$, where the overline notation indicates closure in a mathematical sense. 

\subsubsection{Coupling configurations}
\begin{figure}[H]
\centering
\includegraphics[width=0.7\linewidth]{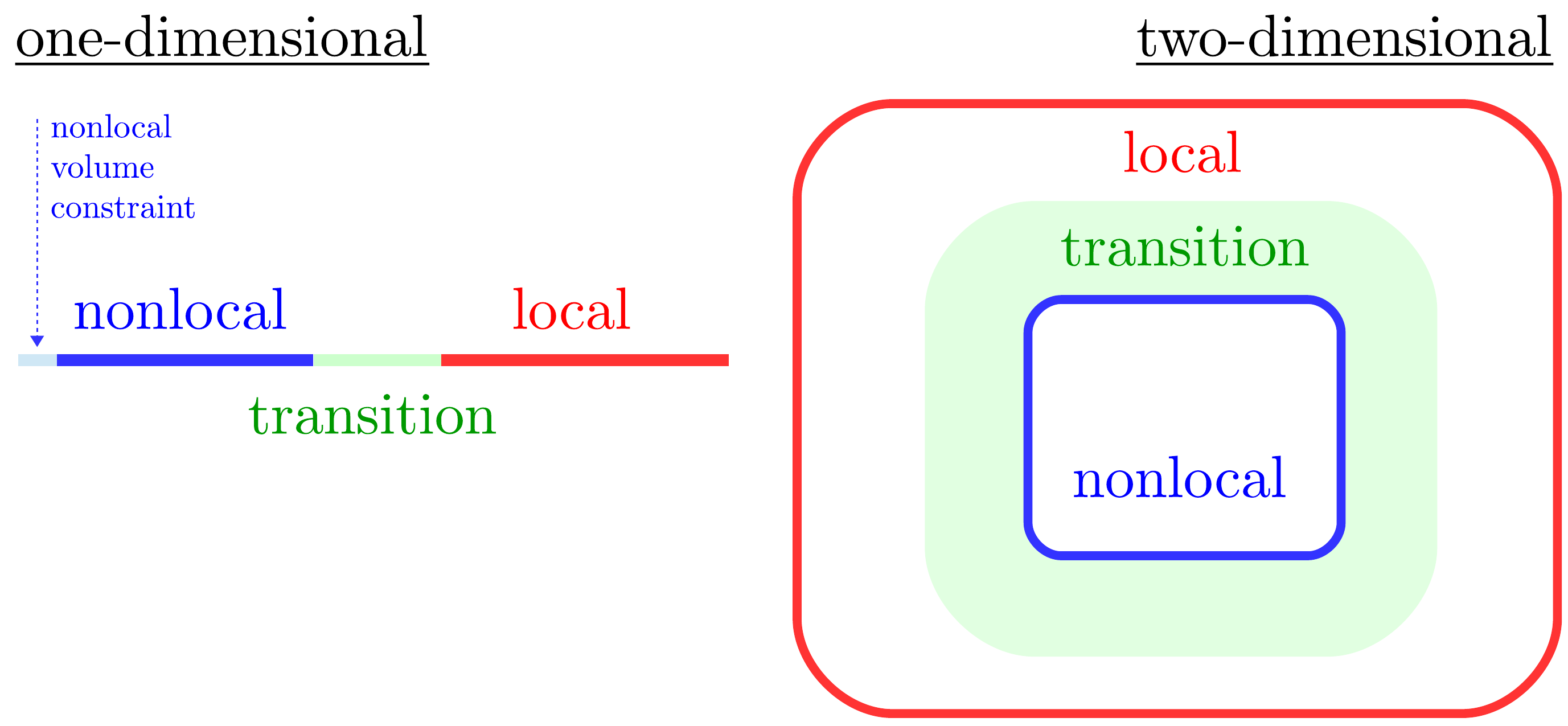}
\caption{General coupling configuration: the domain is decomposed into a purely nonlocal sub-domain, a transition region, and a purely local sub-domain.}
\label{fig:general-domains}
\end{figure}

\begin{figure}[H]
\centering
\includegraphics[width=0.65\linewidth]{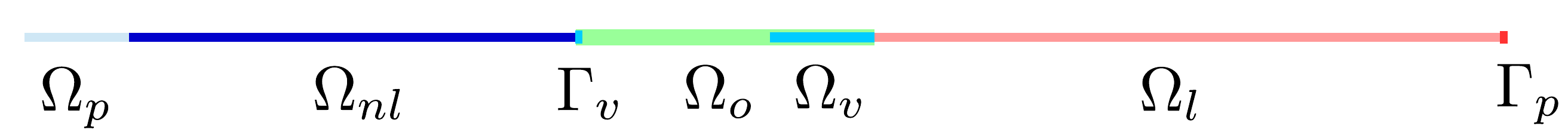}\\[4mm]
\includegraphics[width=0.65\linewidth]{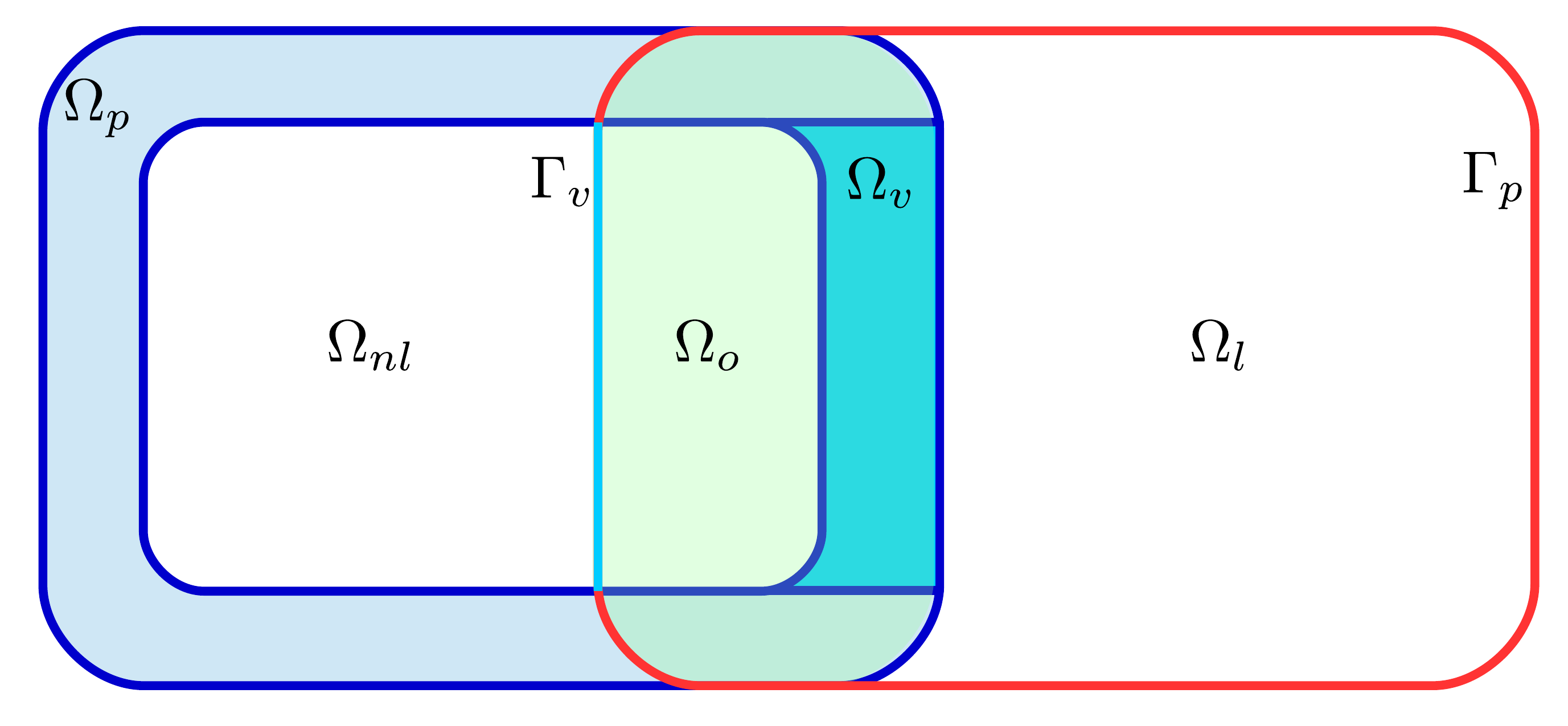}\\[4mm]
\includegraphics[width=0.43\linewidth]{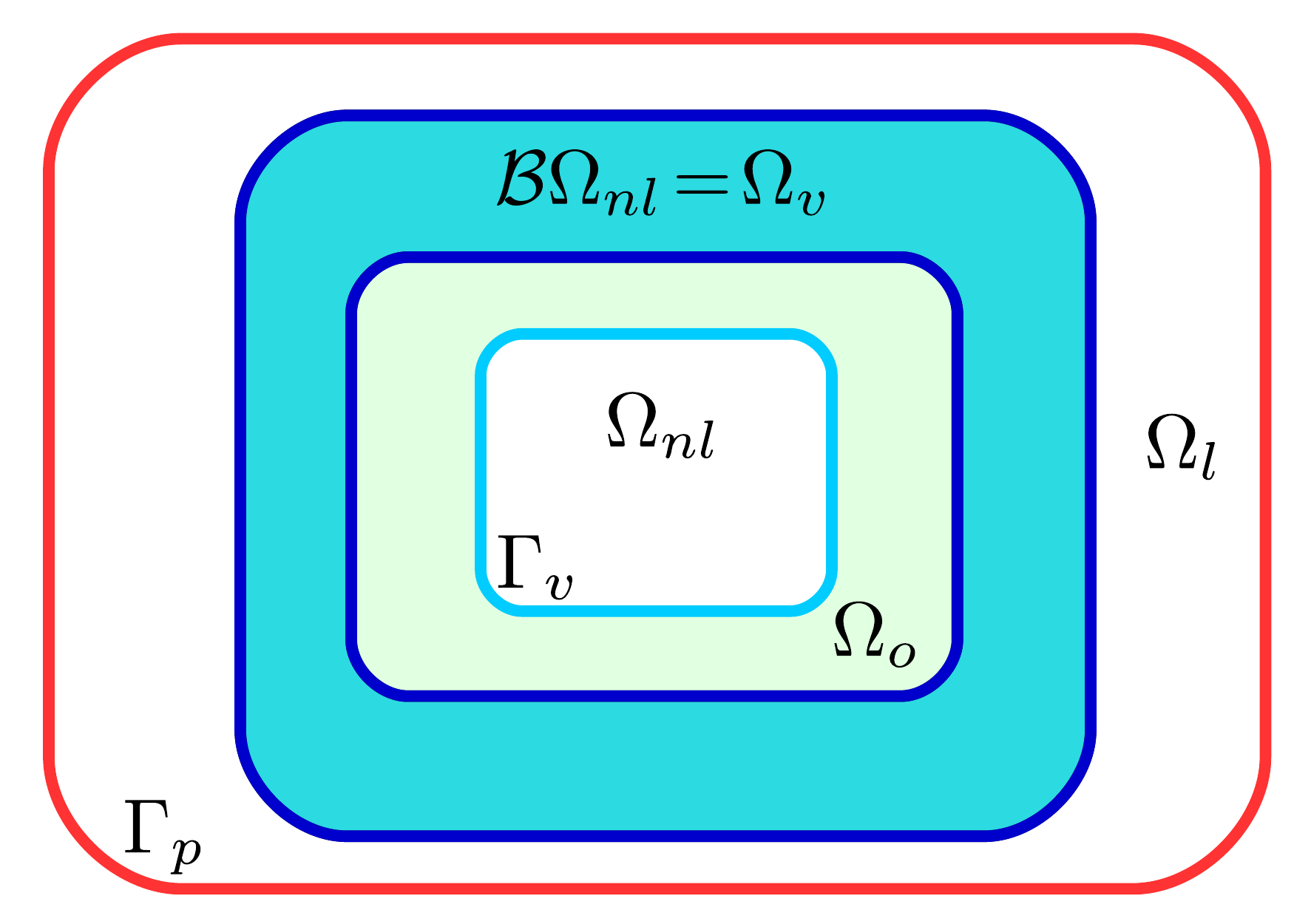}
\caption{Coupling configuration in one dimension (top), two dimensions with adjacent sub-domains (center), and two dimensions with embedded sub-domains (bottom). 
%
%
This is a domain-decomposition setting with overlap (the green region), where both nonlocal and local solutions co-exist.}
\label{fig:decomposed-domains}
\end{figure}

\begin{figure}[H]
\centering
\vspace{-0.5cm}
\includegraphics[width=0.7\linewidth]{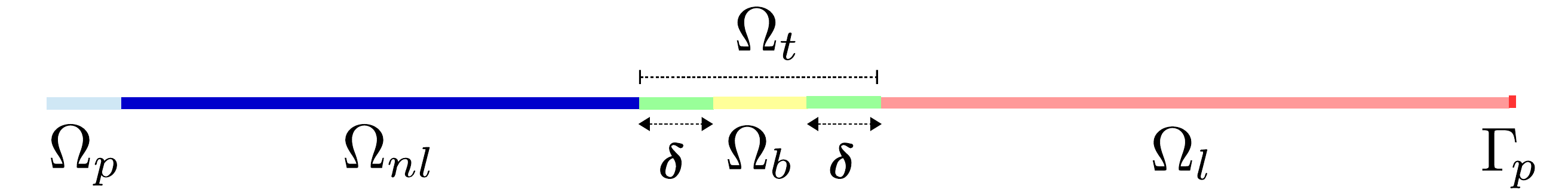}\\[.3cm]
\includegraphics[width=0.5\linewidth]{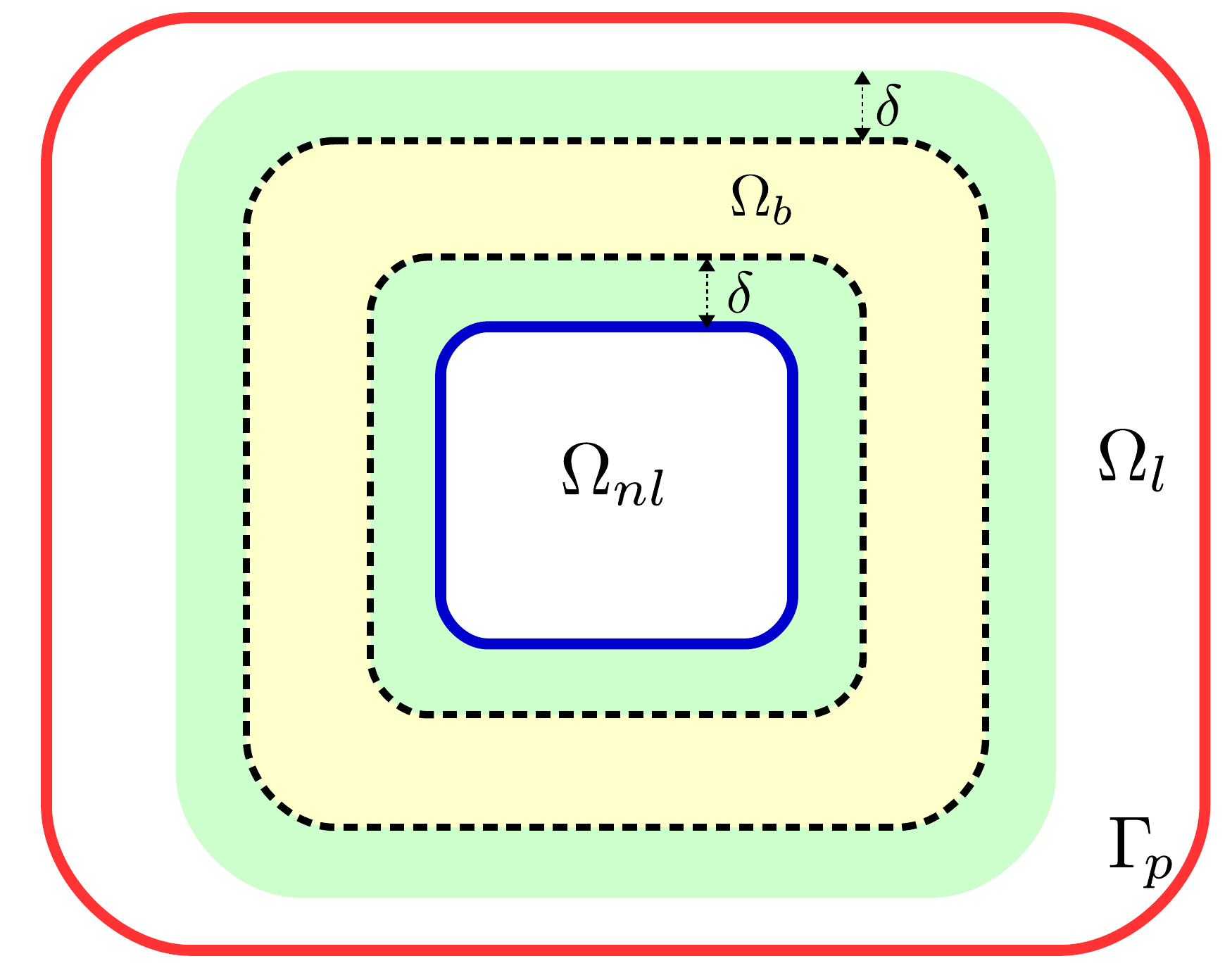}
\caption{Coupling configuration in one dimension (top) and two dimensions with embedded sub-domains (bottom)  for a setting with a transition region $\Omega_t$ (green and yellow) and a blending region $\Omega_b$ (yellow). The blending region is part of the transition region and it is $\horizon$-far from the local and nonlocal sub-domains. 
%
}
\label{fig:Omegab-domain}
\end{figure}

\begin{figure}[H]
\centering
\includegraphics[width=0.7\linewidth]{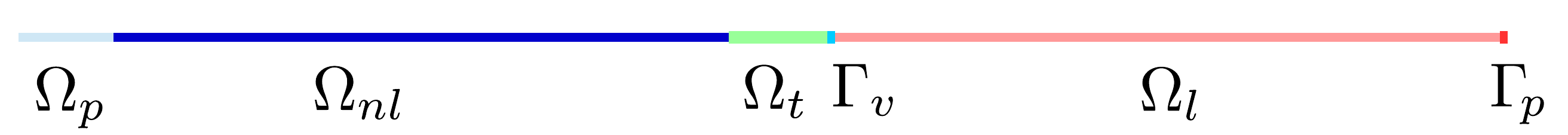}\\[4mm]
\includegraphics[width=0.45\linewidth]{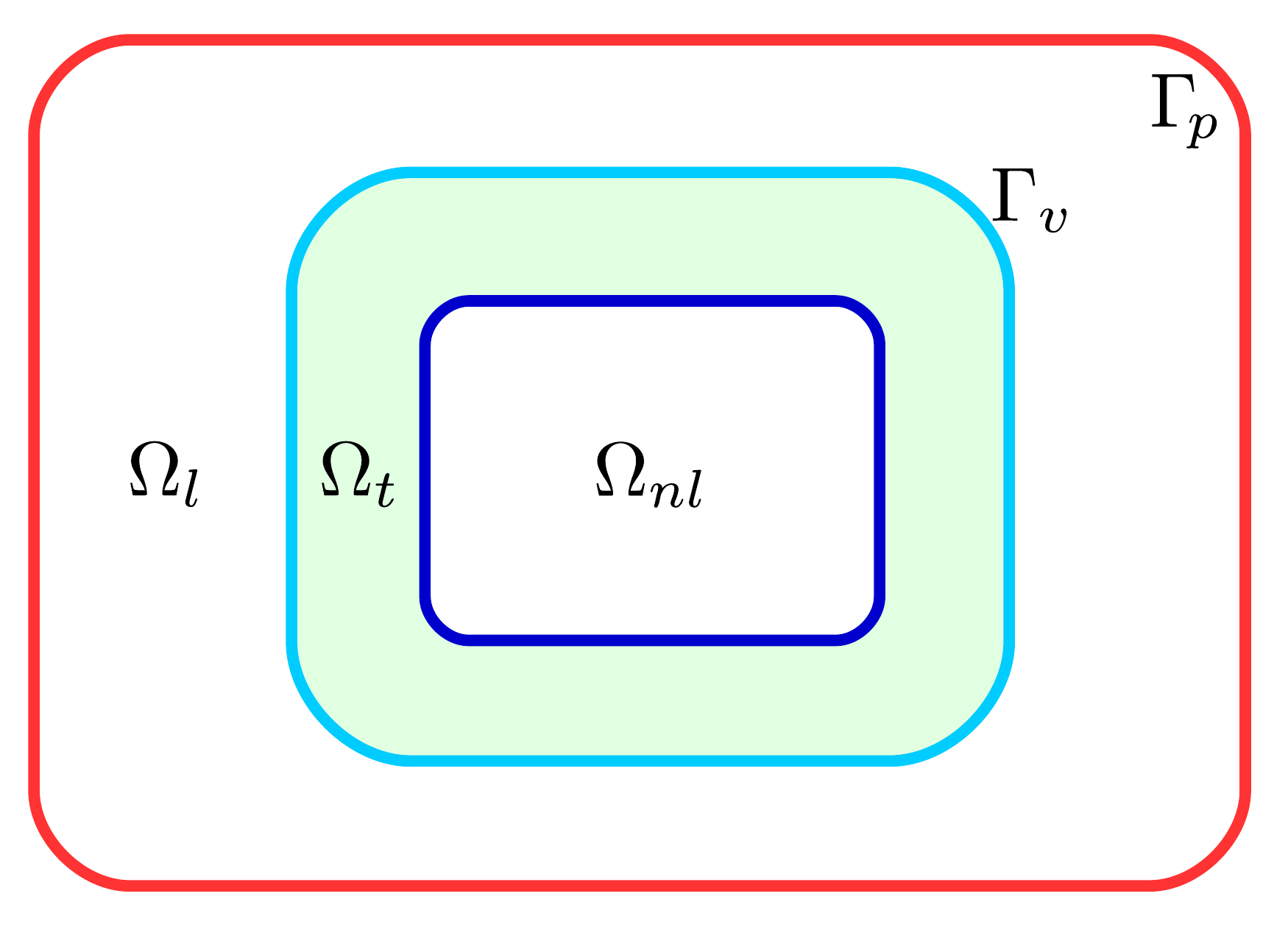}
\caption{Coupling configuration in one dimension (top) and two dimensions with embedded sub-domains (bottom)  for a setting with a transition region.}
\label{fig:blended-domains}
\end{figure}

\begin{figure}[H]  
\centering
\includegraphics[width=0.7\linewidth]{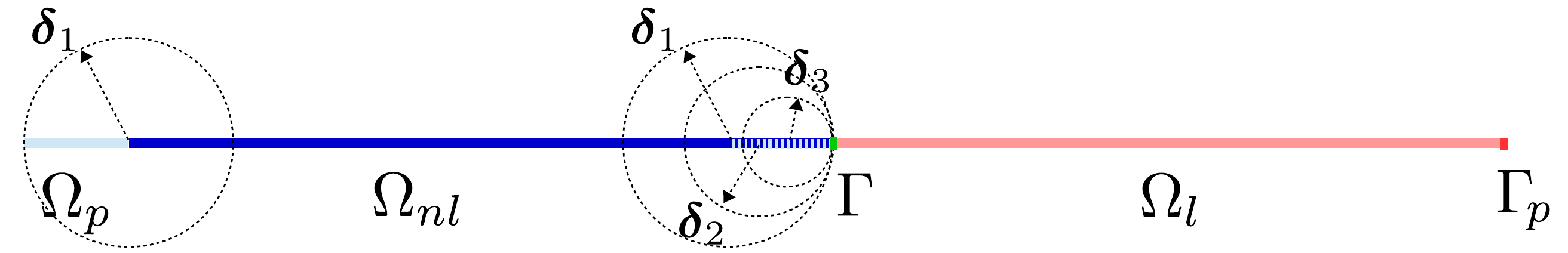}\\[4mm]
\includegraphics[width=0.7\linewidth]{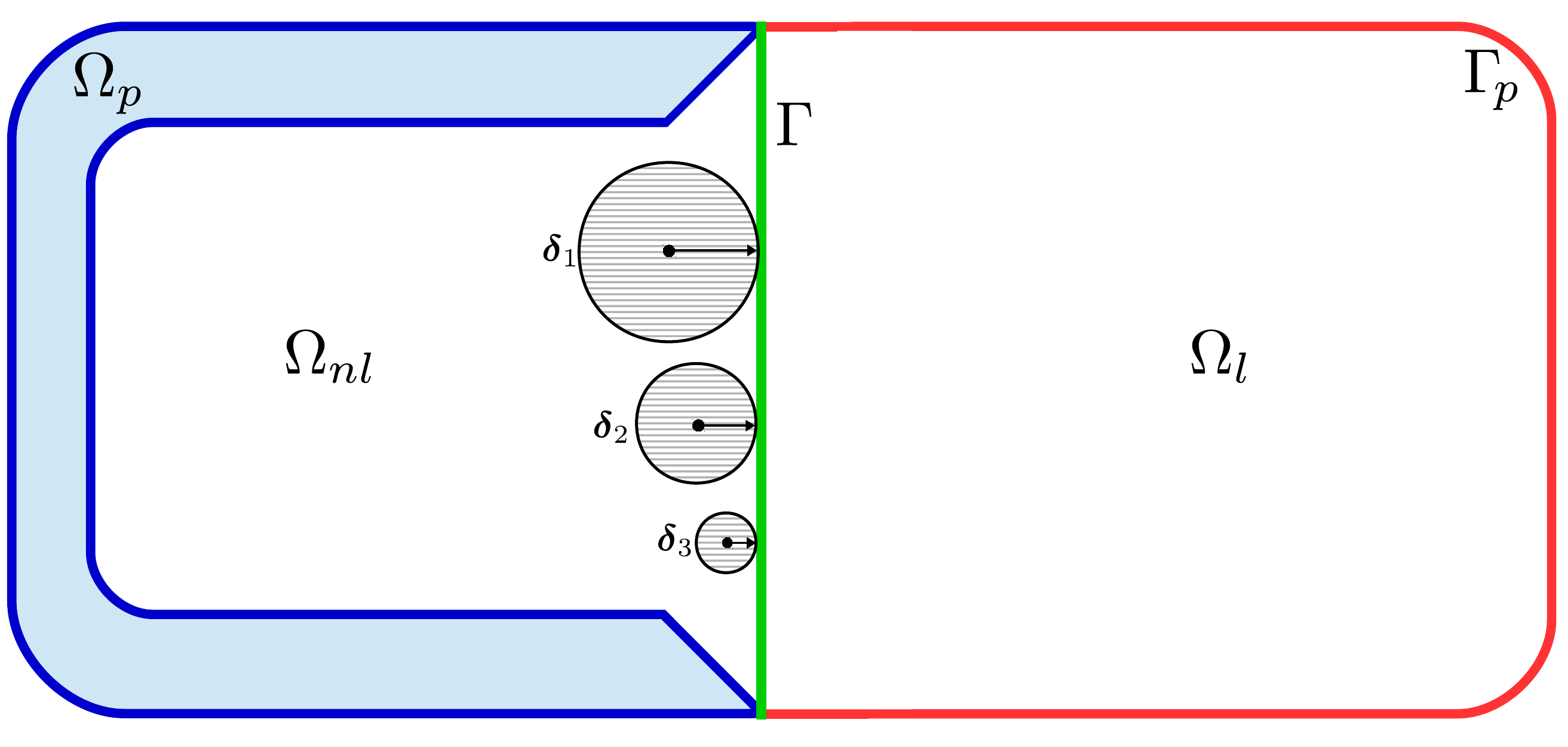}
\caption{Coupling configuration in one (top) and two (bottom) dimensions for a variable horizon setting: as $\xb\in \omgn$ approaches $\Gamma$ the size of the 
nonlocal neighborhood 
linearly shrinks and never crosses the interface.}
\label{fig:variable-domain}
\end{figure}

In a general LtN coupling scenario, the domain $\ooomg$ is decomposed into a purely local sub-domain (described by the local model), a purely nonlocal sub-domain (described by the nonlocal model), and a transition region connecting those two sub-domains. An illustration is provided in Figure \ref{fig:general-domains}. 
In the figure, we report both a one-dimensional (left) and a two-dimensional (right) configuration. In the former, the sub-domains are adjacent, whereas, in the latter, the nonlocal sub-domain and transition region are embedded into the local sub-domain.
The way the transition between the purely nonlocal sub-domain and the purely local sub-domain is performed is method-dependent. For instance, in the CH class, we may have either co-existing local and nonlocal models or a hybrid model in the transition region. On the other hand, in the VH class, the transition region could be identified with the region of varying horizon. To provide a better understanding of the differences between the coupling configurations corresponding to the different methods, we illustrate various configurations in Figures~\ref{fig:decomposed-domains}--\ref{fig:variable-domain}. While any LtN coupling method can be implemented, in general, both on adjacent and embedded configurations, we provide below only the specific configurations used in the following sections.

\paragraph{Domain decomposition with overlap.
}
We refer to Figure \ref{fig:decomposed-domains}: we report the one-dimensional configuration (top),  two-dimensional configuration with adjacent sub-domains (center), 
and two-dimensional configuration with embedded sub-domains (bottom). 
Here, the transition region is the overlap of nonlocal and local sub-domains, where both operators are defined and the solutions co-exist. The domain $\ooomg$ is decomposed as $\ooomg=\overline{\omg}_{l}\cup\ooomg_{nl}$, where $\ooomg_{nl}=\oomg_{nl}\cup\overline{\omgb}_{nl}$ is the union of the nonlocal sub-domain and its nonlocal boundary. The overlapping region 
is defined as $\overline{\omg}_o: =\oomg_{l}\cap\ooomg_{nl}$ and includes the local virtual boundary $\gammav$ and the nonlocal virtual boundary $\omgv$. We also introduce the ``physical'' local and nonlocal boundaries $\gammad=\partial\omgl\!\setminus\!\gammav$ and $\omg_p=\omgb_{nl}\!\setminus\!\omgv$, where we assume that conditions coming from the physics of the problem are provided. 

This configuration is used in Sections \ref{subsec:OBM}, \ref{subsec:Partition_Robin} (for the overlapping case), and 
\ref{subsec:Arlequin}.
%

\paragraph{Domain decomposition without overlap and with blending.}
We refer to Figure \ref{fig:Omegab-domain}: we report 
the one-dimensional configuration (top) and two-dimensional configuration with embedded sub-domains (bottom). 
%
%
In this decomposition, the nonlocal and local sub-domains do not overlap, but are separated by the transition region, $\omgt$. 
A blending function is used as a partition-of-unity function, which changes in the {\it blending} region, $\Omega_b$. 
%
Due to nonlocal contributions, $\omgt = \Omega_b\cup\omgb_{b}$.
In the top figure, the domain is decomposed into four disjoint sub-domains: $\ooomg=\omgd\cup\omg_{nl}\cup\Omega_t\cup\omgl$, i.e.,  the physical nonlocal boundary, the nonlocal sub-domain, the transition region, and the local sub-domain. In the bottom figure, the nonlocal sub-domain is fully embedded in the local sub-domain and $\ooomg$ is decomposed into three disjoint sub-domains: $\ooomg=\omg_{nl}\cup\Omega_t\cup\omgl$. 

This configuration is used in Sections \ref{subsec:morphing} and  \ref{subsec:blending}.

\paragraph{Domain decomposition without overlap and no blending.} 
We refer to Figure \ref{fig:blended-domains}: we report the one-dimensional configuration (top) 
and two-dimensional configuration with embedded sub-domains (bottom). 
In the top figure, the domain is decomposed into {four} disjoint sub-domains: {$\ooomg=\omgd\cup\omg_{nl}\cup\Omega_t\cup\omgl$}, i.e.,  the physical nonlocal boundary, the nonlocal sub-domain, the transition region and the local sub-domain.
In the bottom figure, the nonlocal sub-domain is fully embedded in the local sub-domain and $\ooomg$ is decomposed into three disjoint sub-domains: $\ooomg=\omg_{nl}\cup\Omega_t\cup\omgl$. 
%
This configuration is used in Section \ref{subsec:Partition_Robin} (for the non-overlapping case)
and \ref{subsec:quasinonlocal}.


\paragraph{Domain decomposition with smooth transition from a nonlocal to a local sub-domain.}
In the case of VH coupling methods, we consider the configurations in Figure \ref{fig:variable-domain} for one-dimensional (top) and two-dimensional (bottom) settings. In these configurations, the domain is decomposed by means of a sharp interface $\Gamma$ between a nonlocal, $\ooomg_{nl}=\oomg_{nl}\cup\omgb_{nl}$, and local, $\omgl$, sub-domains. The extent of the nonlocal interactions decreases as points in the nonlocal sub-domain approach $\Gamma$. The same holds true for the nonlocal boundary, that approaches $\partial\ooomg_{nl}$ in the area surrounding the interface. Here, $\omg_p$ coincides with the nonlocal boundary, $\omgb_{nl}$. 
 
This configuration is used in Section \ref{sec:VH}.

\subsubsection{Illustration of blending functions}

\begin{figure}[H]
\begin{center}
\subfigure[Arlequin]{
\includegraphics[width=0.45\textwidth]{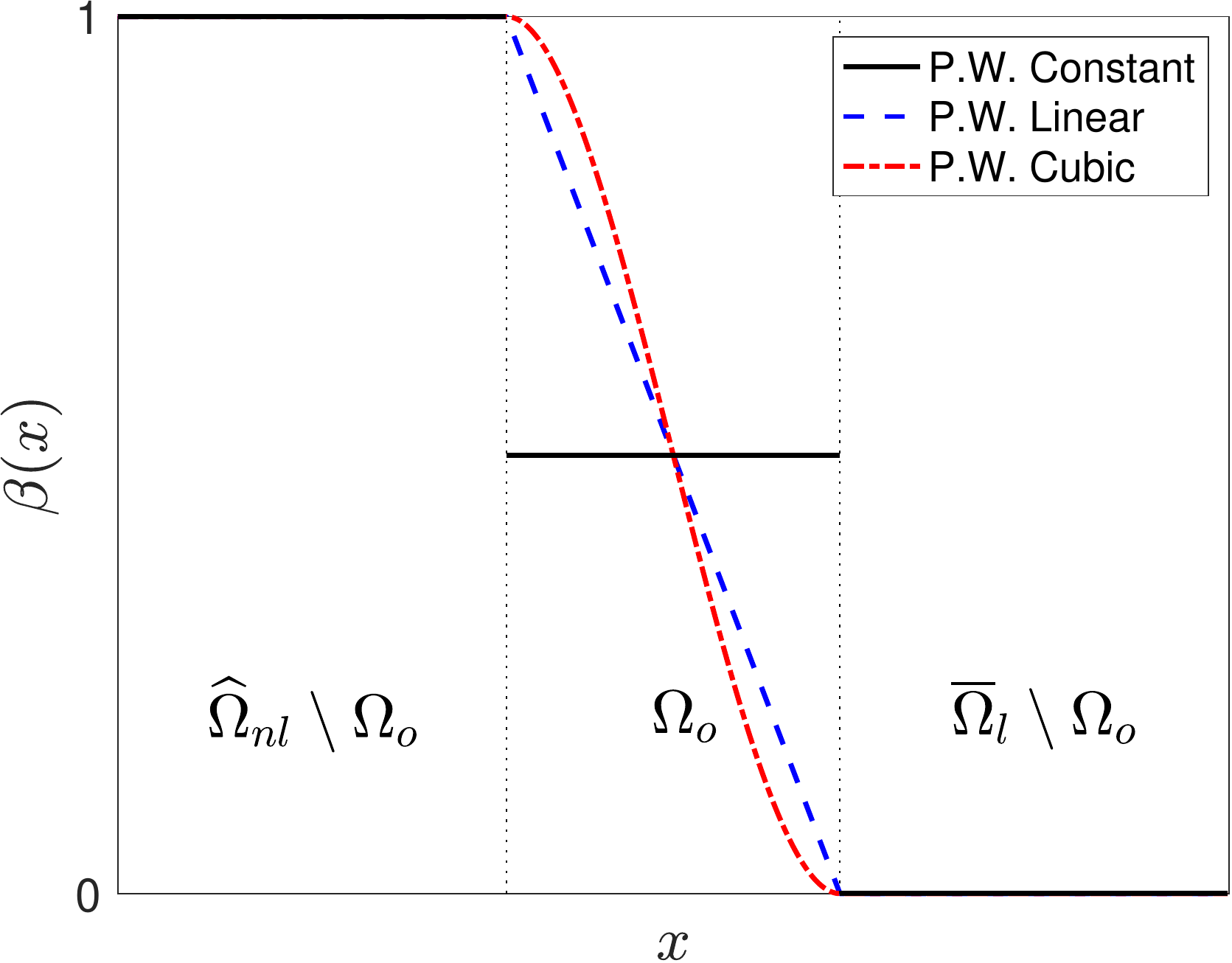}
}
\subfigure[Morphing/Blending]{
\includegraphics[width=0.45\textwidth]{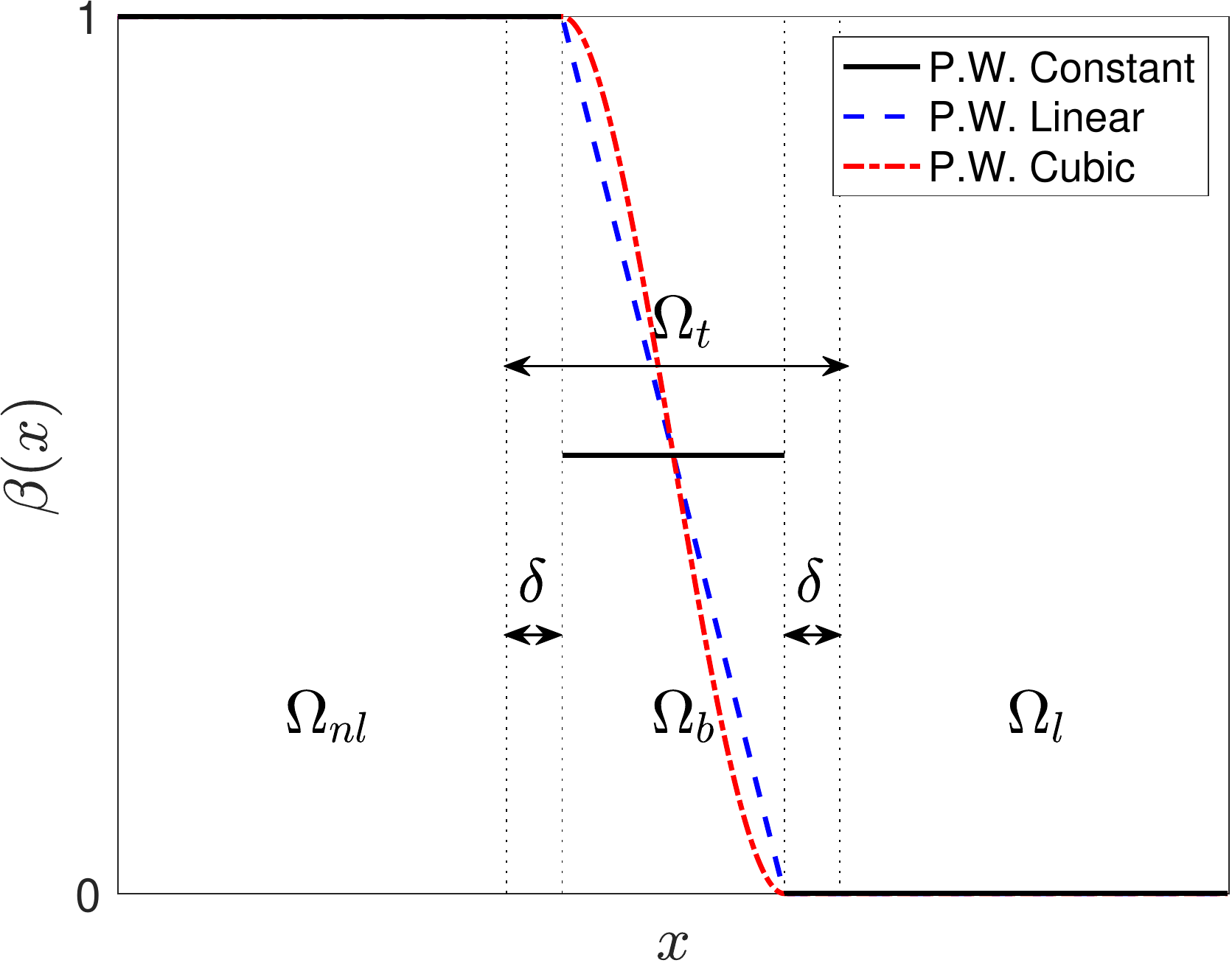}
}
\caption{Illustration of blending functions in one dimension. (a) Domain decomposition with overlap for the Arlequin method. (b) Domain decomposition without overlap and with blending for the morphing and blending methods. Three choices of blending function are shown:  piece-wise (P.W.) constant (black solid line), linear (blue dashed line), and cubic (red dashed-dotted line).}
\label{fig: blending function2}
\end{center}
\end{figure}

In various LtN coupling methods, 
such as the Arlequin, morphing, and blending, the idea of a partition of unity is used by means of  a blending function, $\beta(\xb)$, such that
%
\begin{equation}\label{general_def_blend_fnc}
   \beta(\xb)=\begin{cases}
   & 1 \quad \mbox{in purely nonlocal sub-domain},\\
   & 0 \quad \mbox{in purely local sub-domain},
   \end{cases} 
\end{equation}
which takes values between $0$ and $1$ in the transition region 
(see  Figure \ref{fig:general-domains}).
In the Arlequin method, the transition region coincides with the overlapping region, $\omgo$ (see Figure \ref{fig:decomposed-domains}), and the blending function is chosen as a polynomial, 
normally constant, linear, or cubic, in that region, as illustrated in one dimension in Figure \ref{fig: blending function2}(a). 
In contrast, in the morphing and blending methods, the polynomial choice of the blending function normally occurs within 
a sub-region of the transition region referred to as the blending region, $\omg_b$ (see Figure \ref{fig:Omegab-domain}), resulting in a piece-wise form as illustrated  in one dimension in Figure \ref{fig: blending function2}(b).
%
%
%
%

\subsection{Nonlocal models}\label{sec: nonlocal models}
We describe two nonlocal models: nonlocal diffusion and nonlocal mechanics, specifically peridynamics. We recall that in the case the of diffusion the unknown is the scalar-valued function $u$, whereas for mechanics the unknown is the vector-valued function $\ub$.

\subsubsection{Nonlocal diffusion}\label{subsubsec:diffusion}
Nonlocal diffusion models have been used in many applications such as describing complex turbulence  \cite{Bakunin2008} and nonlocal heat conduction  \cite{bobaru2010peridynamic} when the classical Fick's first law or standard Brownian motion fail to delineate the underlying phenomena \cite{DGLZ,Klafter2004a,Neuman2009,Valdinoci2016}.

Given a domain $\omg \in \Rn$, the time-dependent nonlocal diffusion 
for a point $\xb\in \omg$ at time $t\geqslant 0$ is 
\begin{equation}\label{eq:laplacian}
\frac{\partial u}{\partial t}(\xb,t)=\mcL^{\rm ND} u(\xb,t) +f(\xb,t) = \int_\Rn \gamma(\xbp,\xb)(u(\xbp,t)-u(\xb,t))\,d\xbp +f(\xb,t), 
\end{equation}
where $\gamma:\Rn\times\Rn\to\mbR$ is a symmetric kernel\footnote{This discussion could be generalized to non-symmetric \cite{DElia2017} or sign-changing \cite{Mengesha2013sign-changing} kernels, as well as more general nonlinear kernels~\cite{seleson2013interface}.} and $f$ is a source. In accordance with the definition of nonlocal neighborhood,  
we consider localized kernels of the form
\begin{equation}\label{eq:loc-kernel}
\gamma(\xbp,\xb)=  \mcX_{B_\delta} (\xbp-\xb) k(\xbp,\xb),
\end{equation}
where 
\begin{equation}\label{eq: characteristic function}
\mcX_{\Omega}(\xb) 
:= \left\{
\begin{array}{cl}
  1   & \qquad  \xb \in \Omega \\
  0   & \qquad \mbox{else}\\
\end{array}
\right.
\end{equation}
is the characteristic function of a domain $\Omega$, and we use the notation $B_\delta:=B_\delta (\bf 0) $ (see \eqref{eq:Euclidean-ball}). 
%
%
Since, by definition, the nonlocal boundary consists of those points outside of $\omg$ that interact with points inside of $\omg$, \eqref{eq:loc-kernel} implies that the nonlocal boundary is defined as
\begin{equation}\label{eq:int-domain-diff}
\omgb:=\{\xbp\in\Rn\setminus\omg: \|\xbp- \xb\| \leqslant  \horizon, \;{\rm for \, any}\; \xb\in \omg\}.
\end{equation}
As a consequence of \eqref{eq:loc-kernel}, for $\xb\in\Omega$ the nonlocal operator in \eqref{eq:laplacian} reads
\begin{equation}\label{eq:laplacian-ball}
\mcL^{\rm ND} u(\xb,t) = \int_{B_\horizon(\xb)}k(\xbp,\xb) (u(\xbp,t)-u(\xb,t))\,d\xbp.
\end{equation}
In this work, we mainly focus on static problems, where $u$ depends on $\xb$ only. Also, we consider nonlocal diffusion operators that, in the limit of vanishing nonlocality, i.e., as $\horizon\to 0$, under sutable regularity assumptions, converge to the classical Laplacian or div-grad operator, $\Delta$ or $\nabla\cdot(\nabla)$, respectively. It is possible to show that, when the kernel function $k(\xbp,\xb)$ is properly scaled, we have the following property (see, e.g. \cite{Tao2017,YuNeumann2019}):
\begin{equation}\label{eq:local-limit}
\mcL^{\rm ND} u(\xb) = \Delta u(\xb) + \mcO(\horizon^2) {\rm D}^{(4)}u(\xb),
\end{equation}
where ${\rm D}^{(4)}$ is a combination of the fourth-order derivatives of $u$.\\
\begin{remark}\label{rem:equivalence}
Property \eqref{eq:local-limit} implies that $\mcL^{\rm ND}$ and $\Delta$ are equivalent for polynomials up to the third order, i.e., $\mcL^{\rm ND}p(\xb)=\Delta p(\xb)$, for $p\!\in\!\mathbb P^3(\mbRn)$.
Examples of properly scaled kernels, in one dimension, include
\begin{displaymath}
 k(x',x)=\dfrac{3}{\horizon^3} \quad {\rm and} \quad
k(x',x)=\dfrac{2}{\horizon^2}\frac{1}{|x'-x|}.
\end{displaymath}
%
Note that \eqref{eq:local-limit} also implies that the classical Poisson problem is a $\horizon^2$-approximation of nonlocal diffusion and, as such, it should be used as the reference local problem in the design of LtN coupling methods. Furthermore, the equivalence property implies that polynomials up to order three are good candidates for patch tests.
\end{remark}

We consider nonlocal diffusion problems in the configuration of Figure \ref{fig:standard-domain} and, without loss of generality, we limit our description to the time-independent case. Specifically, given $f:\omg\to\mbR$ and $g:\omgb\to\mbR$, we solve
\begin{equation}\label{eq:ND static problem}
\left\{
\begin{aligned}
-\mcL^{\rm ND} u = f & \quad\xb\in\omg, \\[2mm]
u = g               & \quad\xb\in\omgb,
\end{aligned}\right.
\end{equation}
where the second condition is the nonlocal counterpart of a Dirichlet boundary condition for PDEs. For this reason, we refer to it as Dirichlet volume constraint\footnote{Nonlocal Neumann conditions are also an option. However, for the sake of clarity and without loss of generality, we do not discuss them.}.

\subsubsection{Nonlocal mechanics}\label{subsubsec:mechanics}


In this paper, the nonlocal mechanics models are provided by the peridynamic theory of solid mechanics. 
%
Peridynamics has been applied for material failure and damage simulation (see, e.g., the review paper \cite{diehl2019review}) and seems to provide robust modeling capabilities for analysis of complex crack propagation phenomena, such as  branching~\cite{bobaru2015cracks}.

Given a domain $\omg \in \Rn$, the {\it state-based} peridynamic equation of motion for a point $\xb\in \omg$ at time $t\geqslant 0$ is~\cite{Silling2007}
%
%
\begin{equation}\label{eq: state-based-PD-EOM}
\rho(\xb)\frac{\partial^2 \ub}{\partial t^2}(\xb,t)
= \mcL^{\rm PD} \ub (\xb,t)+ \bb(\xb,t),
\end{equation}
where the nonlocal operator in \eqref{eq: state-based-PD-EOM} is the peridynamic internal force density, 
\begin{equation}\label{eq: PD operator state-based}
\mcL^{\rm PD} \ub (\xb,t) = \int_{ B_\horizon(\xb)} \left\{
\statevecT[\xb,t]\langle \xbp - \xb\rangle  - \statevecT[\xbp,t]\langle \xb - \xbp \rangle
  \right\}d\xbp, 
\end{equation}
$\rho$ is the mass density, 
%
$\statevecT$ is a force vector state,
 and $\bb$ is a prescribed body force density.  
Note that, as for the nonlocal diffusion operator~\eqref{eq:laplacian-ball}, the nonlocal interactions in~\eqref{eq: PD operator state-based} are restricted to the nonlocal neighborhood.
The force vector state $\statevecT[\xb,t]\langle \xib\rangle$ is an operator defined at a given point $\xb$ at time $t$ that maps a peridynamic bond $\xib:=\xbp - \xb$ to force per unit volume squared;  
this operator contains the constitutive relation characterizing the specific material under consideration. 
%
%
Similarly to the case of  nonlocal diffusion, in this work, we mainly consider static problems, where $\ub$ depends only on~$\xb$. In addition,
it has been shown for elastic materials that, under suitable regularity assumptions, the peridynamic 
internal force density convergences, in the limit as $\horizon \to 0$, to the classical elasticity operator \cite{silling2008convergence}
i.e., 
\begin{equation}\label{eq: PD operator state-based limit}
\lim_{\horizon \to 0} \mcL^{\rm PD} \ub (\xb) = \nabla\cdot \nub^0(\xb), 
\end{equation}
%
%
where $\nub^0$ is the collapse peridynamic stress tensor, which is an admissible Piola-Kirchhoff stress tensor.
%

A special case of state-based peridynamics is when the material response of any bond is independent of other bonds. This is referred to as {\it bond-based} peridynamics~\cite{Silling2000} for which the force vector state is given by  \cite{Silling2007}
\begin{equation}\label{eq: bond-based-PD-force-state}
\statevecT[\xb]\langle \xib\rangle = \frac{1}{2} \fb(\etab,\xib),
\end{equation}
where $\fb$ is the pairwise force function and $\etab := \ub(\xb + \xib) -  \ub(\xb)$ is the relative displacement. 
In this case,
\begin{equation}\label{eq: PD operator bond-based bulk}
\mcL^{\rm PD} \ub (\xb) = \int_{B_\horizon({\bf 0})} \fb(\etab,\xib) d\xib.
\end{equation}
%
 %
We observe that \eqref{eq: PD operator bond-based bulk} implies that a nonlocal boundary, as in \eqref{eq:int-domain-diff}, is required to inform points near the domain boundary $\partial\omg$. 

For small deformations, i.e. $\|\etab\|\ll 1$, we can linearize the pairwise force function as\footnote{We assume a pairwise equilibrated reference configuration, i.e., $\fb(\zerob,\xib)= \zerob$ for all $\xib \in \Rn$ and neglect terms of order $\mcO(\|\etab\|^2)$.
}
\begin{equation}\label{eq: bond-based-PD-force-state-linear}
\fb(\etab,\xib) = \lambda(\xib) (\xib\otimes \xib ) \etab,
\end{equation}
%
where $\lambda$ is a scalar-valued micromodulus function, which in the case of isotropy 
$ \lambda(\xib)  =  \lambda(\| \xib \|)$. 
%
%
%
In this case, under a smooth deformation, in the limit of $\horizon \to 0$, we have the following property for~\eqref{eq: PD operator bond-based bulk} (similarly to~\eqref{eq:local-limit})~\cite{seleson2015concurrent}:
\begin{equation}\label{eq: PD operator bond-based bulk limit}
\mcL^{\rm PD} \ub (\xb) = \mathcal{C}_{ijkl}\frac{\partial^2 u_j}{\partial x_k\partial x_l}(\xb) {\bf e}_i + \mcO(\horizon^2) {\rm D}^{(4)}\ub(\xb),
\end{equation}
where we used Einstein summation convention for repeated indices and the orthonormal basis $\{{\bf e}_1,{\bf e}_2,{\bf e}_3\}$, and 
$\mathcal{C}_{ijkl}$ are the components of 
 a fully-symmetric fourth-order elasticity tensor, which is related to the micromodulus function by 
\begin{equation}\label{eq: fourth-order elasticity tensor bond-based}
\mathcal{C}= \frac{1}{2}  \int_{B_\horizon({\bf 0})}\lambda(\|\xib\|) \xib \otimes \xib \otimes \xib \otimes \xib d\xib.
%
\end{equation}
In the limit as $\horizon \to 0$, \eqref{eq: PD operator bond-based bulk limit} reduces to the classical linear elasticity operator
%
\begin{equation}\label{eq: CE-Navier-Cauchy}
\mcL^{\rm CE} \ub (\xb) = 
\left\{
\begin{array}{cc}
\displaystyle
\frac{4E}{5} \left[ \nabla^2\ub(\xb)  + \frac{1}{2}\nabla \otimes \nabla \ub(\xb)\right]
    & \qquad \text{in three dimensions,}\\
\displaystyle
\frac{3E}{5} \left[ \nabla^2\ub(\xb)  + \frac{1}{2}\nabla \otimes \nabla \ub(\xb)\right]
    & \qquad \text{in two dimensions,}\\
\displaystyle
E\frac{\partial^2 u}{\partial x^2}(x)
   & \qquad \text{in one dimension},
\end{array}
\right.
\end{equation}
where $E$ is Young's modulus. Note that the resulting Poisson's ratio is constrained to $\nu = 1/4$ (in three dimensions) and $\nu = 1/3$ (in two dimensions)~\cite{TrageserSeleson2019}; this restriction is characteristic of bond-based peridynamics and is removed in state-based peridynamics.

We consider peridynamic problems in the configuration of Figure \ref{fig:standard-domain} and, without loss of generality, we limit our description to the static case. Specifically, given $\bb:\omg\to\Rn$  and $\gb:\omgb\to\Rn$, we solve
\begin{equation} \label{eq: PD static problem}
\left\{
\begin{aligned}
-\mcL^{\rm PD} \ub(\xb) = \bb(\xb) & \quad\xb\in\omg, \\[2mm]
\ub(\xb) = \gb(\xb)       & \quad \xb\in\omgb,
\end{aligned}\right.
\end{equation}
where the second condition is, as in the nonlocal diffusion case, the nonlocal counterpart of a Dirichlet boundary condition for PDEs. 

\subsubsection{Nonlocal vector calculus}\label{subsub:nonlocal_calculus}

In this section, we reformulate nonlocal diffusion and peridynamic operators in terms of the {\it  nonlocal vector calculus} developed in \cite{Du2013}. This theory is the nonlocal counterpart of the classical calculus for differential operators and, by providing a variational setting, allows one to analyze nonlocal problems in a very similar way as for PDEs. Such calculus is necessary to analyze several LtN coupling methods described in this review; as such, we only report relevant results and refer the reader to \cite{DGLZJoE,Du2013} for more details.

Following the notations in \cite{Du2013}, the {\itshape nonlocal divergence operator} $\mcD$ acting on a two-point function $\nub(\xb,\xbp): \R^n\times\R^n\to \R^n$ is defined as
\[
\mcD(\nub)(\xb) := \int_{\R^n} (\nub(\xb,\xbp) +\nub(\xbp,\xb))\cdot  \alphab(\xb, \xbp) d\xbp \,,
\]
where $\alphab(\xb,\xbp)$ is anti-symmetric, i.e., $\alphab(\xb,\xbp)=-\alphab(\xbp,\xb)$. 
Without loss of generality, we assume
\begin{equation} \label{kernel:alpha}
\alphab(\xb,\xbp) =\frac{\xbp-\xb}{\|\xbp-\xb\|}.
\end{equation}
The {\itshape nonlocal gradient operator} $\mcD^\ast$ acting on a one-point function $u(\xb):\R^n\to \R$ is defined by
\[
\mcD^\ast(u) (\xb,\xbp):= - (u(\xbp)-u(\xb)) \alphab(\xb, \xbp) \,.
\]
Using the notations of nonlocal divergence and gradient operators, the nonlocal diffusion operator $\mcL^{\rm{ND}}$ in \eqref{eq:laplacian} is then written as 
\[
\mcL^{\rm{ND}}  u(\xb) = - \frac{1}{2}\mcD(\gamma \mcD^\ast u) (\xb)\,,
\]
where $\gamma=\gamma(\xbp, \xb)$ is the kernel function defined in Section \ref{subsubsec:diffusion}. We next define the energy functional associated with the problem \eqref{eq:ND static problem}\\
\[
\begin{split}
E_\delta^{\rm{ND}}(u, f,  g):=& \frac{1}{4}\int_{\ooomg}\int_{\ooomg} \gamma(\xbp,\xb) ((\mcD^\ast u)(\xb,\xbp))^2 d\xbp d\xb  \\
&-\int_{\Omega} f(\xb) u(\xb) d\xb - \int_{\omgb} g(\xb) u(\xb) d\xb.
\end{split}
\]
The corresponding nonlocal energy norm, nonlocal energy space, and constrained energy space are then given, respectively, by
\begin{align*}
   \| u\|_{\rm{ND}} &  = (E_\delta^{\rm{ND}}(u, 0, 0))^{1/2}\,, \\
   \mcS^{\rm{ND}}(\ooomg) & =\{ u\in L^2(\ooomg) : \| u\|_{\rm{ND}} <\infty \} \,,\\
   \mcS^{\rm{ND}}_c(\ooomg) & =  \{ u\in \mcS^{\rm{ND}}(\ooomg) : u|_{\omgb} =0 \}\,.
\end{align*}

Given the tensor two-point function $\Psib:\R^n\times \R^n \to \R^{n\times n}$
and the one-point function $\ub: \R^n \to \R^n$, we can similarly define the nonlocal divergence operator $\mcD \Psib$ for tensors and its adjoint $\mcD^\ast \ub$ by 
\begin{align*}
\mcD (\Psib) (\xb) & =\int_{\R^n} (\Psib(\xb,\xbp) +\Psib(\xbp,\xb))\cdot  \alphab(\xb, \xbp) d\xbp\,, \\
\mcD^\ast(\ub)(\xb, \xbp) & =  - (\ub(\xbp)-\ub(\xb))\otimes\alphab(\xb, \xbp)\,,
\end{align*}
where $\alphab$ is given by \eqref{kernel:alpha}. The linear bond-based peridynamic operator $\mcL^{\rm{PD}}$ in \eqref{eq: PD operator bond-based bulk} with \eqref{eq: bond-based-PD-force-state-linear} is then written as 
\[
\mcL^{\rm{PD}}(\xb) = -\frac{1}{2}\mcD(\gamma (\mcD^\ast \ub)^T )(\xb)\,,
\]
where $\gamma=\gamma(\xbp, \xb):=\lambda(\|\xbp-\xb\|) \|\xbp-\xb\|^2$. We can also define the energy functional for problem \eqref{eq: PD static problem}:
\[
\begin{split}
E_\delta^{\rm{PD}}(\ub, \bb, \gb):=& \frac{1}{4}\int_{\ooomg}\int_{\ooomg} \gamma(\xbp,\xb) (\rm{Tr}( \mcD^\ast \ub)(\xb,\xbp))^2 d\xbp d\xb \\
& -\int_{\Omega}  \bb(\xb) \ub(\xb) d\xb - \int_{\omgb} \gb(\xb) \ub(\xb) d\xb  \,, 
\end{split}
\]
where $\rm{Tr}( \mcD^\ast \ub)$ is the trace of $\mcD^\ast \ub$. 
The corresponding energy norm, energy space, and constrained energy space are given, respectively, by
\begin{align*}
  \| \ub\|_{\rm{PD}} & = (E_\delta^{\rm{PD}}(\ub, {\bf 0}, {\bf 0}))^{1/2}\,, \\
   \mcS^{\rm{PD}}(\ooomg) & =\{ {\ub}\in L^2(\ooomg,\R^n) : \| {\ub}\|_{\rm{PD}} <\infty \} \,,\\
  \mcS^{\rm{PD}}_c(\ooomg) & =  \{ {\ub}\in \mcS^{\rm{PD}}(\ooomg) : \ub|_{\omgb} =0 \}\,.
\end{align*}
Note that $L^2(\ooomg)= L^2(\ooomg, \R)$ under our notation.

The nonlocal energy spaces $\mcS^{\rm{ND}}(\ooomg)$ and $\mcS^{\rm{PD}}(\ooomg)$ defined in the above are Hilbert spaces, as shown in \cite{Mengesha2013sign-changing,mengesha2014}. Moreover, the nonlocal diffusion problem \eqref{eq:ND static problem} and the nonlocal mechanics problem \eqref{eq: PD static problem} are well-posed as a result of Poincar\'e type inequalities. Here we will state the inequalities without proof, and more details can be found in \cite{Du2012,Mengesha2013sign-changing,mengesha2014}.
\begin{lemma}[Nonlocal Poincar\'e inequalities]
There exists a positive constant $C$ such that the following nonlocal Poincar\'e type inequalities hold: 
\begin{align*}
  & \| u \|_{L^2(\ooomg)}\leqslant C\| u\|_{\rm{ND}}, \quad \forall u \in\mcS^{\rm{ND}}_c(\ooomg) \,, \\
  & \| \ub \|_{L^2(\ooomg, \R^n)}\leqslant C\| \ub\|_{\rm{PD}}, \quad \forall \ub \in\mcS^{\rm{PD}}_c(\ooomg)\,.   \\
 \end{align*}
\end{lemma}

\section{Generalized domain-decomposition methods}\label{sec:GDD}

\subsection{Optimization-based methods}\label{subsec:OBM}
In Optimization-Based Methods (OBMs) the coupling of local and nonlocal models is effected by couching the LtN coupling into an optimization problem. This approach is inspired by GDD 
methods for PDEs 
\cite{Discacciati_13_SIOPT,Du_01_SINUM,Du_00_SINUM,Gervasio_01_NM,Gunzburger_00_AMC,Gunzburger_00_SINUM,Gunzburger_99_CMA}
and it has also been applied to AtC coupling 
in \cite{Bochev_14_LSSC,Bochev_14_SINUM}. 
%
A main feature of OBMs is that numerical solutions only require the implementation of the optimization strategy as the local and nonlocal solvers can be used as black boxes. For this reason, OB couplings (OBCs) are considered {\it non-intrusive} as opposed to other methods, whose implementation requires modification of the basic governing equations in the transition region between local and nonlocal sub-domains. 
 
\subsubsection{Mathematical formulation}

We provide a very general formulation of an OBM. Let $\mcL^{\rm NL}$ be a nonlocal operator that accurately describes the system in a bounded domain $\omg$ and let $\mcL^{\rm L}$ be a corresponding local operator that describes the system well enough where nonlocal effects are negligible. OBMs tackle the LtN coupling by solving a minimization problem where the difference between the local and nonlocal solutions is minimized in the overlap between local and nonlocal sub-domains, tuning their values on the virtual boundaries 
induced by the domain decomposition (see Figure \ref{fig:decomposed-domains}). Here, we consider a time-independent problem and we state the LtN OBC in a very general form for vector functions:
\begin{equation}\label{eq:formal-OBM}
\begin{aligned}
&\hspace{-.3cm} \min\limits_{\ub^{nl},\ub^l,\nub^{nl},\nub^l} \frac12 \,
                \|\ub^{nl}-\ub^l\|^2_{*,\omg_o} \quad\mbox{such that}\quad\\[2mm]
& 
\;\left\{
\begin{array}{rll}
- \mcL^{\rm NL} \ub^{nl} =& \!\bb & \xb\in\omgn \\[1mm]
\ub^{nl} = & \!{\bf g} &  \xb\in\omgd\\[1mm]
\ub^{nl} = & \!\nub^{nl} & \xb\in\omgv
\end{array}\right.
\quad\mbox{and}\quad
\left\{\begin{array}{rll}
- \mcL^{\rm L} \ub^l =& \!\bb & \xb\in\omgl \\[1mm]
\ub^l = & \!{\bf g} &  \xb\in\gammad\\[1mm]
\ub^l = & \!\nub^l & \xb\in\gammav
\end{array}\right.,
\end{aligned}
\end{equation}
where $\bb$ is a forcing term over $\omgn\cup\omgl$, ${\bf g}$ is a boundary data over $\omg_p$ and  $\gammad$, $\|\cdot\|_{*,\omg_o}$ is a suitable norm in the 
overlapping region 
$\omg_o$, and $(\nub^{nl},\nub^l)\in\mcC$ (an appropriate control space) are the control variables. 
%
The goal of OBC is to find optimal values for the virtual controls $\nub^{nl}$ and $\nub^l$ such that $\ub^{nl}$ and $\ub^l$ are as close as possible in the 
overlapping region,
while still satisfying their corresponding governing equations. 
We denote the \textit{optimal} controls and corresponding nonlocal and local solutions by $\nub^{nl*},\nub^{l*}$ and $\ub^{nl*},\ub^{l*}$, respectively. The global coupled solution is then defined as 
\begin{equation}\label{eq:OBMsolution}
\ub^*(\xb) = \ub^{nl*}(\xb)\, \mcX_{\omgn}(\xb) + \ub^{l*}(\xb)\, \mcX_{\omgl}(\xb) \qquad \mbox{for all }\xb \in \ooomg,
\end{equation}
where $\mcX_{\omgn}(\cdot)$ and $\mcX_{\omgl}(\cdot)$ are the characteristic functions of $\omgn$ and $\omgl$, respectively (see \eqref{eq: characteristic function}).

Applications of this technique to nonlocal diffusion can be found in \cite{Bochev_14_INPROC,Delia2019} in one dimension and in \cite{Bochev_16b_CAMWA} in three dimensions, whereas applications to nonlocal mechanics can be found in \cite{DElia2015handbook}.

We now present a more rigorous formulation of OBC for a nonlocal diffusion problem, following \cite{Delia2019}. In this case, the unknowns are the scalar-valued functions $u^{nl}\!\in\! \mcS^{\rm{ND}}_c(\ooomg_{nl})$ and $u^l\!\in \!H^1_\gammad(\omgl)$\footnote{We denote by $H^1_\gammad(\omgl)$ the space of functions in $H^1(\omgl)$ that vanish on $\gammad$.} such that
\begin{equation}\label{eq:subprob}
\left\{\begin{array}{rlll}
-\mcL^{\rm ND}u^{nl} & = & f & \xb\in\omgn\\[0.3ex]
u^{nl} & = & 0 & \xb\in\omg_p\\[0.3ex]
u^{nl} & = & \nu^{nl} & \xb\in\omgv
\end{array}\right.
\quad\mbox{and}\quad
\left\{\begin{array}{rlll}
-\Delta u^l & = & f & \xb\in\omgl\\[0.3ex]
u^l & = & 0 & \xb\in\gammad\\[0.3ex]
u^l & = & \nu^l & \xb\in\gammav
\end{array}\right.
\end{equation}
where $\mcL^{\rm ND}$ is the nonlocal diffusion operator introduced in Section \ref{subsubsec:diffusion} and $(\nu^{nl},\nu^l)\in \mcC:=\widetilde \mcS^{\rm{ND}}_c \times H^\frac12(\gammav)$, where $\widetilde \mcS^{\rm{ND}}_c = \{ v|_\omgv\!:\! v\in \mcS^{\rm{ND}}_c(\ooomg_{nl})\}$ 
and $v|_\omgv$ denotes a restriction of $v$ to $\omgv$,  
are  undetermined nonlocal and local Dirichlet boundary conditions. 
Thus, the OBC can be formulated as the following constrained optimization problem:
\begin{equation}\label{eq:OBM-minimization}
\begin{array}{c}
\displaystyle\min\limits_{u^{nl},u^l,\nu^{nl},\nu^l}  \frac12 \| u^{nl}-u^l \|^2_{L^2(\omg_o)}
\quad
\mbox{subject to \eqref{eq:subprob}.}
\end{array}
\end{equation}
Given the optimal controls $(\nu^{nl*},\nu^{l*})$, the global coupled solution is  defined as in \eqref{eq:OBMsolution}.

\paragraph{Properties.} 
OBC approaches such as \eqref{eq:OBM-minimization} have several desirable properties:

\begin{itemize}
\item Provided that the nonlocal and local equations are well-posed (see Section \ref{subsub:nonlocal_calculus}), the optimal control problem is well-posed, i.e., \eqref{eq:OBM-minimization} has a unique solution. In particular, it is possible to show that the Euler-Lagrange equations associated with the reduced functional 
$\mathcal J(\nu^{nl},\nu^l)=\frac12 \| u^{nl}(\nu^{nl})-u^l(\nu^l) \|^2_{L^2(\omg_o)}$ define a coercive variational form in the control space $\mathcal C$ \cite{Delia2019,Bochev_16b_CAMWA}, where the notations $u^{nl}(\nu^{nl})$ and $u^l(\nu^l)$ denote the solutions $u^{nl}$ and $u^l$ given the virtual controls $\nu^{nl}$ and $\nu^l$, respectively.
\item When $\mcL^{\rm ND}$ is properly scaled (see Remark \ref{rem:equivalence}) the method passes the linear, quadratic and cubic patch tests. Furthermore,  \cite{Delia2019} shows that as $\horizon\to 0$ the coupled solution approaches the solution of a global local problem, i.e., the method is asymptotically compatible. However, this result can be achieved as long as the overlapping region $\omgo$ is ``large enough'', as convergence estimates depend on $|\omgo|^{-1}$, where $|\omgo|$ denotes the size of $\omgo$.
\item 
The solution of nonlocal and local problems is completely uncoupled. As a consequence, nonlocal and local governing equations can be discretized on separate computational sub-domains with different discretization methods (e.g., meshfree method for the nonlocal problem and the finite element method (FEM) for the local problem). This implies that available software for the solution of each equation can be used as a black box.
\item
The mathematical formulation presented in this section can be easily extended to the peridynamic model presented in Section \ref{subsubsec:mechanics} by simply substituting $\mcL^{\rm NL}$ with $\mcL^{\rm PD}$, either in the state-based or bond-based forms, and $\mcL^{\rm L}$ with the corresponding classical elasticity operator (e.g., $\mcL^{\rm CE}$ in \eqref{eq: CE-Navier-Cauchy} for bond-based peridynamics), as demonstrated in Section \ref{sec:OBMresults}. 
%
\item
The optimization does not affect the rate of convergence of the discretization method used for the governing equations (e.g., a piece-wise linear FEM discretization preserves second order convergence of the solution in the $L^2$-norm).
\end{itemize}

\paragraph{The time-dependent problem.}
The extension of this method to time-dependent problems is straightforward. However, due to the high computational cost, still has limited applicability. Such extension consists in performing the optimization at every time step of the 
time integration. 
%
Specifically, at every time step the objective functional is the same as in \eqref{eq:formal-OBM}, but the constraints are the semi-discrete (in time) nonlocal and local equations. An even more expensive option is to formulate the time-dependent coupling as a continuous global (in space and time) optimization problem, where the objective functional is the sum of the misfits over the whole domain and at any time instant and the constraints are the time-dependent nonlocal and local equations. None of these approaches has been rigorously analyzed nor implemented.

\subsubsection{Applications and results}\label{sec:OBMresults}

\begin{figure}[H]
\centering
{\includegraphics[width=0.5\textwidth]{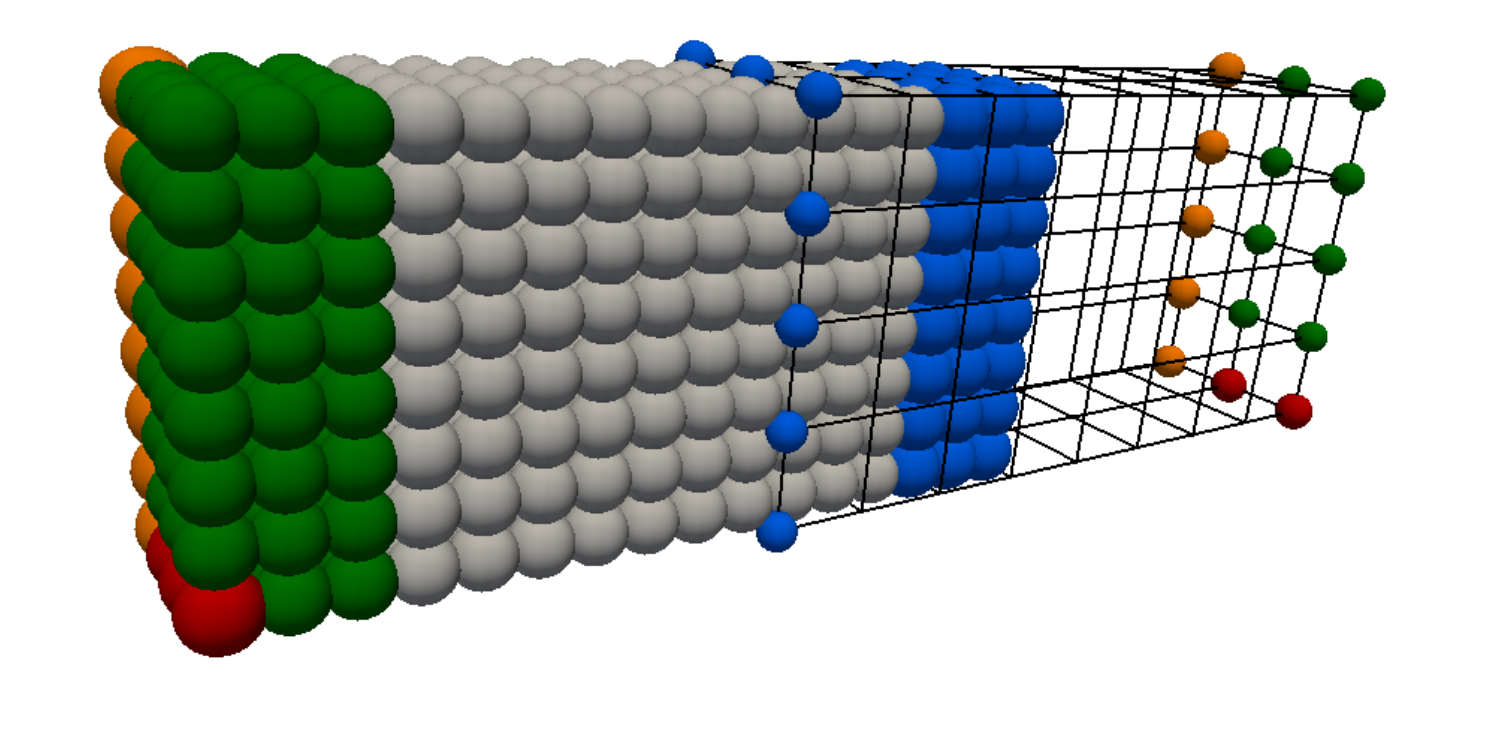}}
\vspace{-5ex}
\caption{Computational domain for linear and quadratic patch tests. On the right, a meshfree discretization; on the left, an FEM discretization. 
Blue nodes indicate the control variables degrees of freedom. 
}
\label{fig:patch-tests-setup}
\end{figure}

\begin{figure}[H]
\begin{center}
\subfigure[Linear Patch Test]{
\includegraphics[width=0.45\textwidth]{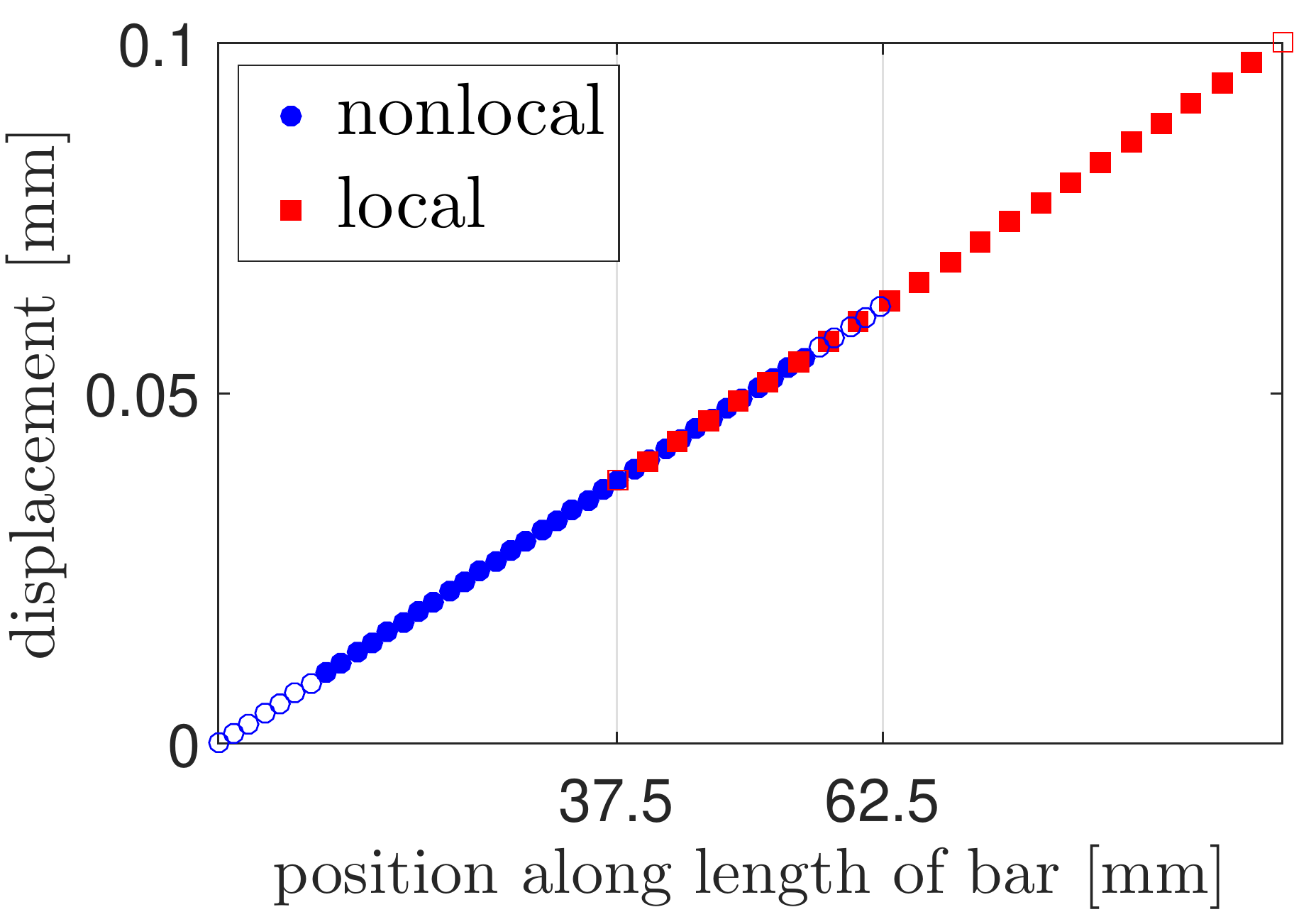}
}
\subfigure[Quadratic Patch Test]{
\includegraphics[width=0.45\textwidth]{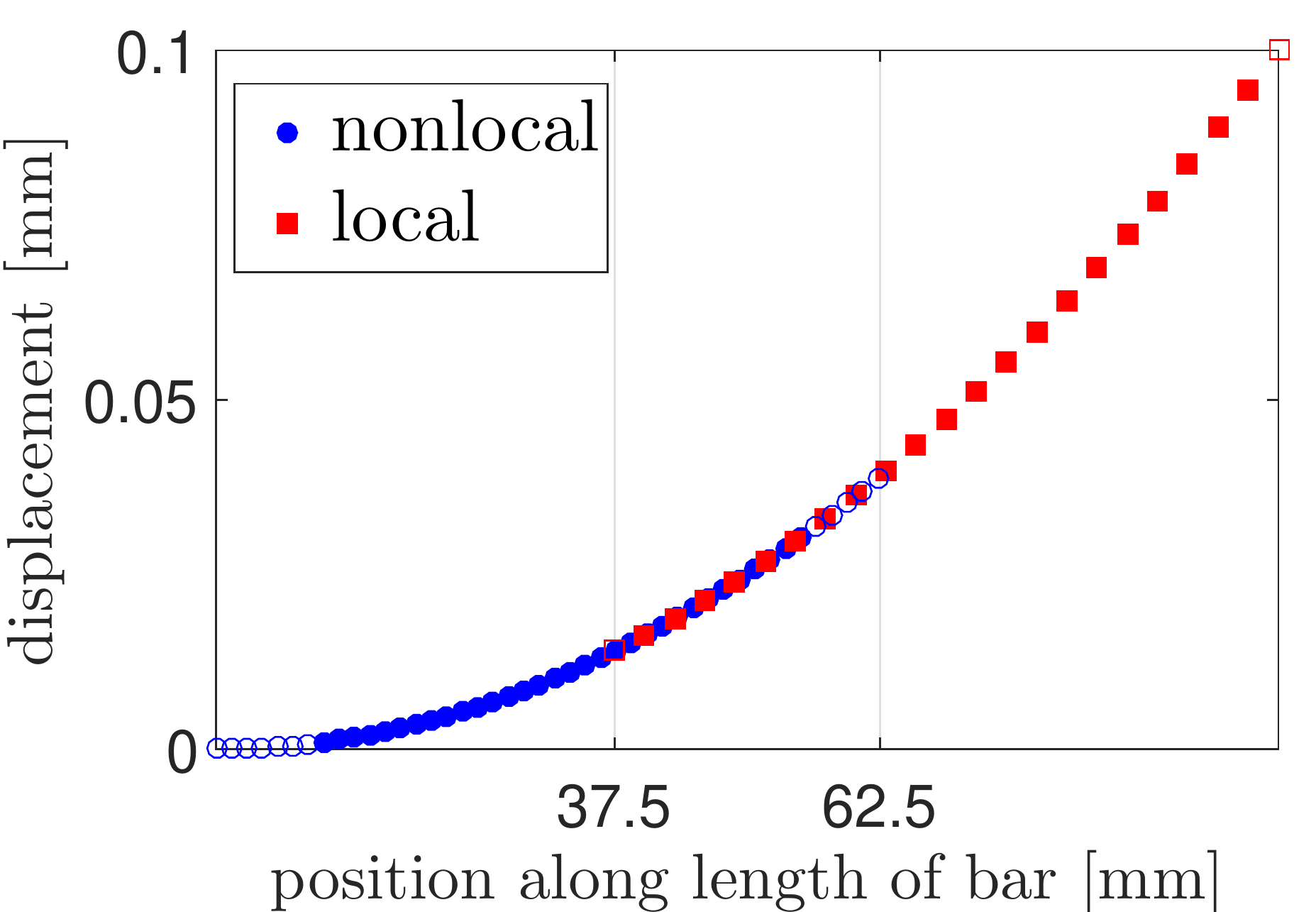}
}
\caption{Patch tests for the optimization-based method: linear (left) and quadratic (right). The figures show displacements along the horizontal line passing through the center of the bar in Figure \ref{fig:patch-tests-setup}. Nodes in the nonlocal sub-domain are represented by blue filled circles, whereas nodes in the (physical and virtual) nonlocal boundaries are empty blue circles. Nodes in the local sub-domain are represented by red filled squares, whereas nodes in the (physical and virtual) local boundaries are empty red squares.}
\label{fig:patch-tests-1D}
\end{center}
\end{figure}

As anticipated in the previous section, the abstract formulation in \eqref{eq:formal-OBM} can be easily extended to state-based peridynamics. 
Here, 
we report results for the linearized version (see \cite{Silling2010}) of the linear peridynamic solid model introduced 
in \cite{Silling2007}. 
This is an isotropic state-based peridynamic model that 
converges, in the limit as $\horizon\!\to\! 0$ and under suitable regularity assumptions, to classical linear elasticity 
without the Poisson's ratio restriction discussed in Section \ref{subsubsec:mechanics}.  
We report linear and quadratic patch tests for a stainless steel bar as obtained in \cite{DElia2015handbook} (where all the specifics are listed). 
In both tests, the nonlocal problem is discretized using a meshfree method \cite{silling2005meshfree} (see also \cite{Littlewood2015Peridigm}) 
and the local problem is discretized using piece-wise linear FEM. 
The computational domain is reproduced in Figure \eqref{fig:patch-tests-setup}. 
On the left, the nonlocal problem is discretized with a meshfree method; 
on the right, the local problem is discretized with FEM .
Results are reproduced in Figure \ref{fig:patch-tests-1D}, where displacement solutions in the horizontal direction are reported along a horizontal line passing through the center of the bar. The patch test results are in good agreement with the expected linear and quadratic solutions for both the nonlocal and local problems. 
%

Further results in \cite{DElia2015handbook} on stainless steel bars with cracks demonstrate the ability of OBC to capture cracks in nonlocal regions as well as the ease of dealing with boundary conditions that can be applied to local regions circumventing the non-trivial task of prescribing nonlocal boundary conditions.

\subsection{Partitioned procedure}
\label{subsec:Partition_Robin}

In this section we introduce another type of GDD methods, known as partitioned procedures, which can be applied to heterogeneous domains such as those characterized by multiscale and multiphysics dynamics. In the partitioned procedure, the coupled problem is solved based on subsequent solutions of sub-problems. The sub-problems are coupled through transmission conditions on interfaces, which enforce the continuity of solutions and the energy balance of the whole system. One of the main features of this approach is software modularity. In fact, since each sub-problem is solved separately, possibly with different solvers, the partitioned procedure allows for the reuse of existing codes/discretization methods. Moreover, because only  transmission conditions on the coupling interfaces are needed and provided for each solver, the partitioned procedure is also generally non-intrusive, similarly to OBMs.

While partitioned procedures can be applied to static and quasi-static problems, they are more commonly used for dynamic problems. Therefore, we describe the partitioned formulation for dynamic LtN coupling problems. The partitioned procedure can be broadly classified as either explicit or implicit (in the literature these two approaches are sometimes referred to as loosely coupled strategy and strongly coupled strategy, respectively). In explicit coupling strategies, the solution of each sub-problem and the exchange of interface data are performed only once (or a small number of times) per time step. For reasons pointed out below, this causes an information mismatch on the interface. As a consequence, the coupled system does not satisfy exactly the coupling transmission conditions and the energy exchanged between sub-problems is not perfectly balanced.
In contrast, in implicit coupling strategies, each sub-problem is solved in a partitioned way via sub-iterations until convergence is reached. Thus, transmission conditions are always satisfied and the energy of the system is balanced. This generally allows one to achieve numerical stability, even for  
highly nonlinear coupled systems.

The partitioned procedure is introduced in \cite{you2019coupling} for nonlocal diffusion problems and in \cite{yu2018partitioned} for nonlocal mechanics problems. In both works, the authors use Robin transmission conditions since, compared to Dirichlet and Neumann transmission conditions, are more flexible and often exhibit better robustness.

\subsubsection{Mathematical formulation}
To introduce the partitioned procedure formulation for LtN coupling  with Robin transmission conditions, we consider a general time-dependent LtN coupling problem. At time step $k$ and sub-iteration step $j$, we solve for $\ub^{nl}_{k,j}$ and $\ub^{l}_{k,j}$ using the solutions at the previous time step, $\ub^{nl}_{k-1}$ and $\ub^{l}_{k-1}$, and the solutions at the previous sub-iteration, $\ub^{nl}_{k,j-1}$ and $\ub^{l}_{k,j-1}$. We first describe the partitioned procedure for the overlapping case, i.e., with a configuration corresponding to Figure \ref{fig:decomposed-domains}:
\begin{equation}\label{eq:dynamic-robin}
\begin{aligned}
& \left\{
\begin{array}{rll}
D_t{\ub^{nl}_{k,j}}- \mcL^{\rm NL} \ub^{nl}_{k,j} =& \!\sbp_k & \qquad \xb\in\omgn \\[2mm]
\ub^{nl}_{k,j} = & \!{\bf g}_k & \qquad \xb\in\omgd\\[2mm]
R_{1}\ub^{nl}_{k,j}+T^{\rm NL}(\ub^{nl}_{k,j}) = & \!R_{1}\ub^{l}_{k,j-1}+T^{\rm L}(\ub^{l}_{k,j-1}) & \qquad \xb\in\omgv
 \end{array}\right.\\[2mm]
&\left\{\begin{array}{rll}
D_t{\ub^l_{k,j}}- \mcL^{\rm L} \ub^l_{k,j} =& \!\sbp_k & \qquad\xb_n\in\omgl \\[2mm]
\ub^l_{k,j} = & \!{\bf g}_k &  \qquad\xb\in\gammad\\[2mm]
R_2\ub^{l}_{k,j}+T^{\rm L}(\ub^{l}_{k,j}) = & \!R_2\ub^{nl}_{k,j}+T^{\rm NL}(\ub^{nl}_{k,j}) &\qquad \xb\in\gammav
\end{array}\right.,
\end{aligned}
\end{equation}
where $\sbp$ is a forcing term over $\omgn\cup\omgl$ and ${\bf g}$ is a boundary data over $\omg_p$ and $\gammad$. For diffusion problems, $\ub$ is a scalar-valued function, $\mcL^{\rm NL}=\mcL^{\rm ND}$, $\mcL^{\rm L}=\Delta$, and $\sbp=f$ (see Section \ref{subsubsec:diffusion}). 
%
For mechanics problems,
 $\ub$ is a vector-valued function, $\mcL^{\rm NL}\ub(\mathbf{x},t)=\dfrac{1}{\rho(\mathbf{x})}\mcL^{\rm PD}\ub(\mathbf{x},t)$, $\mcL^{\rm L}\ub(\mathbf{x},t)=\dfrac{1}{\rho(\mathbf{x})}\nabla\cdot \nub^0(\mathbf{x},t)$, and $\sbp(\mathbf{x},t)=\dfrac{\bb(\mathbf{x},t)}{\rho(\mathbf{x})}$ (see Section \ref{subsubsec:mechanics}). $T^{\rm L}$ is a flux operator in diffusion problems or a normal stress operator in mechanics problems, and $T^{\rm NL}$ is a corresponding nonlocal operator to be defined. $D_t$ is the appropriate numerical time-difference operator to approximate a first-order derivative (in diffusion problems) or a second-order derivative (in mechanics problems) of $\ub$. $R_1$ and $R_2$ are constant Robin coefficients. 
Note that when $R_i=0, i=1,2$, the Robin transmission condition is equivalent to the Neumann-type transmission condition, whereas when $R_i\rightarrow\infty$, $i=1,2$, the Dirichlet transmission condition is obtained. In explicit coupling procedures, the results are updated as $\ub^{nl}_{k}=\ub^{nl}_{k,j}$ and $\ub^{l}_{k}=\ub^{l}_{k,j}$ after one (or a few) sub-iteration. Note that in explicit coupling procedures, the Robin transmission condition in the nonlocal sub-problem becomes $R_{1}\ub^{nl}_{k}+T^{\rm NL}(\ub^{nl}_{k}) =  \!R_{1}\ub^{l}_{k-1}+T^{\rm L}(\ub^{l}_{k-1}) $, so the transmission condition $R_{1}\ub^{nl}+T^{\rm NL}(\ub^{nl}) =  \!R_{1}\ub^{l}+T^{\rm L}(\ub^{l}) $ is not exactly satisfied, and therefore the energy of the system is not balanced. In implicit coupling procedures, the iterations stop only when a prescribed stopping criterion is reached, and the energy exchanged on the interface is balanced. 
To apply Robin transmission conditions in LtN coupling problems, one of the challenges is to provide a well-defined nonlocal operator $T^{\rm NL}$. In the  overlapping case, both local and nonlocal solutions exist in ${\Omega}_o$. 
In \cite{yu2018partitioned}, the authors proposed to employ the local operator $T^{\rm L}(\ub^{nl})$ as an approximation for the nonlocal operator $T^{\rm NL}(\ub^{nl})$. 

In \cite{you2019coupling}, the authors considered an LtN coupling configuration for the non-overlapping case  
as in Figure \ref{fig:blended-domains}, except that the nonlocal sub-domain includes both $\Omega_{nl}$ and $\Omega_t$, where the local and nonlocal solvers interact on a sharp 
interface 
that coincides with the local interface $\Gamma_v$. {On the local side, a Dirichlet transmission condition $u^{l}_{k,j}=u^{nl}_{k,j}$ is employed on $\Gamma_v$. On the nonlocal side, instead of explicitly defining the nonlocal Neumann-type operator $T^{\rm NL}$, the authors developed a nonlocal formulation which converts the local flux to a correction term in the nonlocal model and provides an estimate for the nonlocal interactions across
$\Gamma_v$. For $\xb\in\omg_t$, where $\omg_t$ is a collar of thickness $\delta$ adjacent to the LtN coupling interface (see Figure \ref{fig:blended-domains}), a modified nonlocal formulation is employed: }
\begin{align*}
&Q_\delta(\xb) D_t \ub^{nl}_{k,j}(\xb,t)-\mcL_{N\delta}^{\rm NL} \ub^{nl}_{k,j}(\xb,t)+R_1V_\delta(\xb)\ub^{nl}_{k,j}(P(\xb),t)\\
&~~~~~=Q_\delta(\xb)\sbp_k(\xb,t)+V_\delta(\xb) T^{\rm L}(\ub^{l}_{k,j-1}(\xb,t))+R_1V_\delta(\xb) \ub^{l}_{k,j-1}(P(\xb),t),
\end{align*}
where $Q_\delta$ and $V_\delta$ are functions of $\xb$, $P$ is the projection operator onto $\Gamma_v$, and $\mcL_{N\delta}^{\rm NL}$ is a modified nonlocal operator as proposed in \cite{you2019coupling}. This formulation provides a generalization of the classical local Robin transmission condition, such that the nonlocal problem converges to the corresponding local problem with local Robin transmission condition, as $\delta\rightarrow 0$, under suitable regularity assumptions.

\paragraph{Properties.} To achieve a stable and/or fast convergent LtN  coupling algorithm, one needs to choose optimal Robin coefficients, $R_1$ and $R_2$; see next section. With proper Robin coefficients, the partitioned procedure has several properties:
\begin{itemize}
\item The partitioned procedure is not tied to any particular discretization. Actually, it provides a flexible coupling framework which enables software modularity: different discretization methods, including different spatial and time resolutions, 
can be employed for the two sub-problems.
\item When the nonlocal sub-problem with Robin transmission condition is asymptotically compatible to the local limit, the coupling framework preserves the asymptotic compatibility.
%
\item The partitioned procedure passes the patch tests 
as long as the solution 
satisfies the chosen Robin transmission condition. For example, when taking $T^{\rm NL}$ as the corresponding local operator $T^{\rm L}$ in the partitioned procedure for the overlapping case, the coupling framework passes up to cubic patch tests.  
On the contrary, in the partitioned procedure for the non-overlapping case, as in \cite{you2019coupling}, the coupling framework simply passes the linear and quadratic patch tests 
only when the interface is flat.
\item The partitioned procedure with the implicit coupling approach is energy preserving. However, the explicit coupling approach does not necessarily guarantee energy preservation.
\item The partitioned procedure for the non-overlapping case can also be applied to general heterogeneous LtN coupling problems, i.e., where the local and nonlocal problems have 
different physical properties.
\end{itemize}

\subsubsection{Applications and results}

In applications, the optimal Robin coefficients can be estimated either theoretically or numerically. In problems with relatively simple and/or structured domain settings, one can {perform a Fourier decomposition of the solution and define an analytical reduction factor as the convergence ratio of the semi-discretized solution in the frequency space. 
Then, the optimal Robin coefficient can be obtained by minimizing the analytical reduction factor, as shown in \cite{yu2018partitioned}}. There, the authors consider a LtN coupling for two-dimensional static/quasi-static mechanics problems, involving a nonlinear bond-based peridynamic model given by the prototype microelastic brittle (PMB) model \cite{silling2005meshfree} discretized with a meshfree method,  
coupled to a classical linear elasticity model (see \eqref{eq: CE-Navier-Cauchy}) discretized with the FEM.
%
%
%
On a simple problem setting representing a plate under uniaxial tension, the optimal Robin coefficient is provided by simplifying the coupling problem to a one-dimensional model problem and obtaining expressions for the reduction factor of the solution. The developed optimal Robin transmission condition is also applied to capture crack initiation and growth for a plate with a hole loaded with increasing tension.

\begin{figure}[H]\centering
\subfigure{\includegraphics[width=0.365\textwidth]{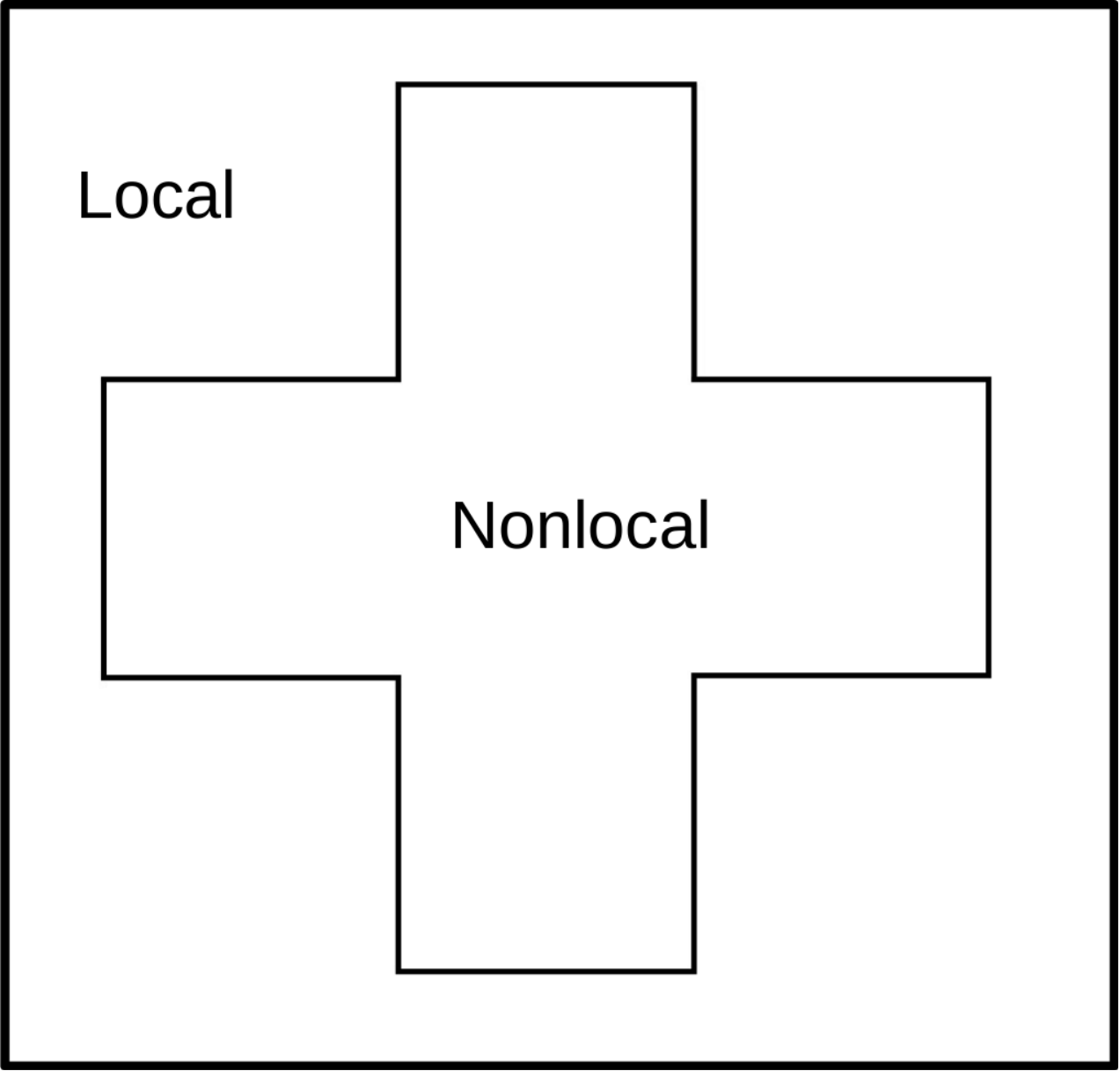}}
\subfigure{\includegraphics[width=0.6\textwidth]{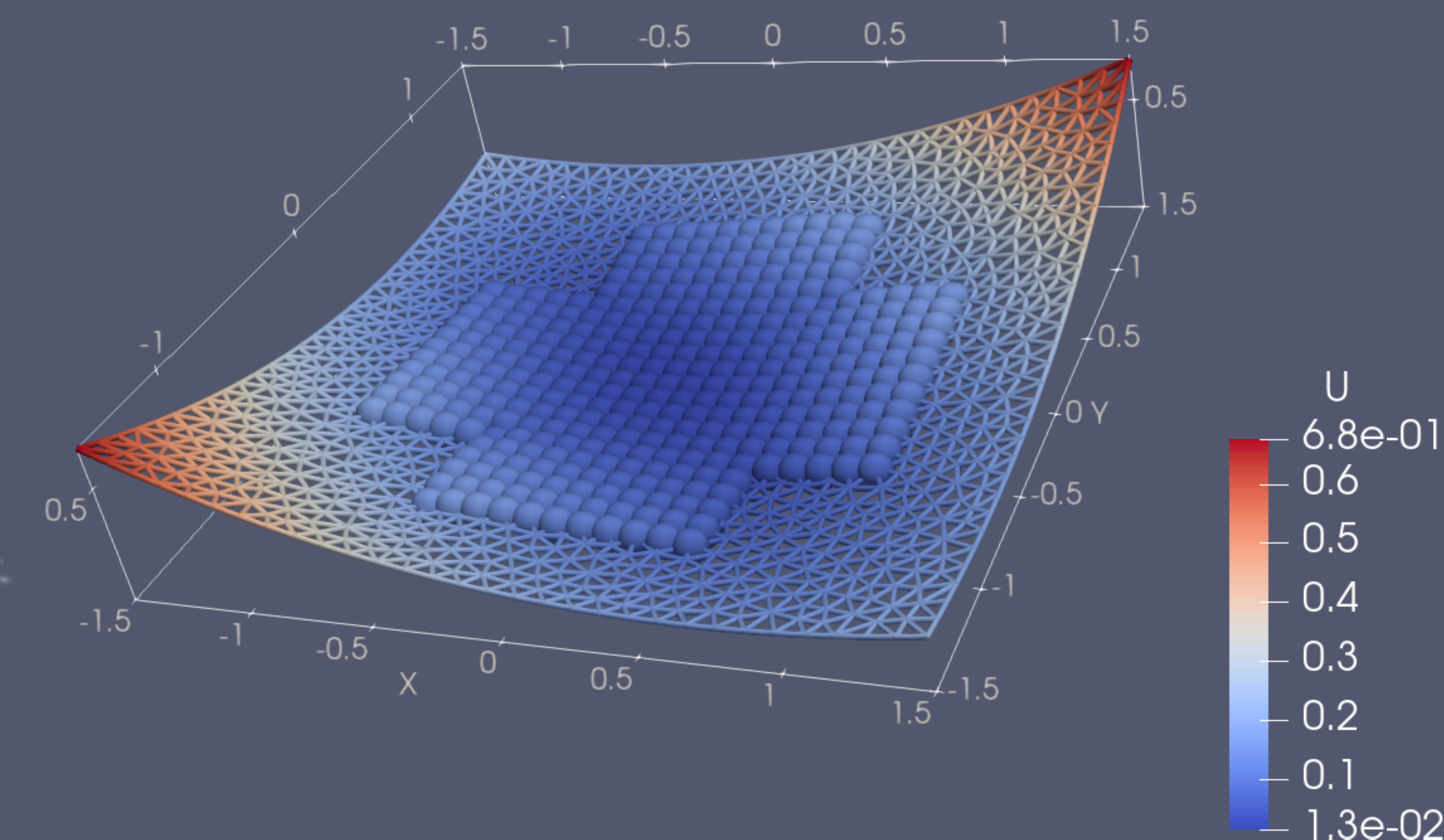}}
\caption{{Non-overlapping partitioned procedure for LtN coupling. 
Left: problem setting. Right: simulation results where the blue spheres represent the nonlocal solution with a meshfree solver and the triangular mesh represents the local solution obtained via FEM.}}
\label{Figure:cross}
\end{figure}

\begin{figure}[H]
\begin{center}
\subfigure[Linear Patch Test]{
\includegraphics[width=0.49\textwidth]{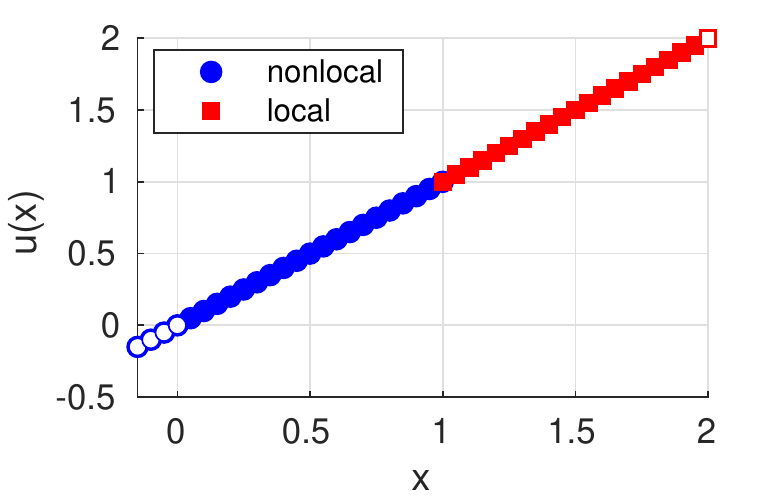}}
\subfigure[Quadratic Patch Test]{\includegraphics[width=0.49\textwidth]{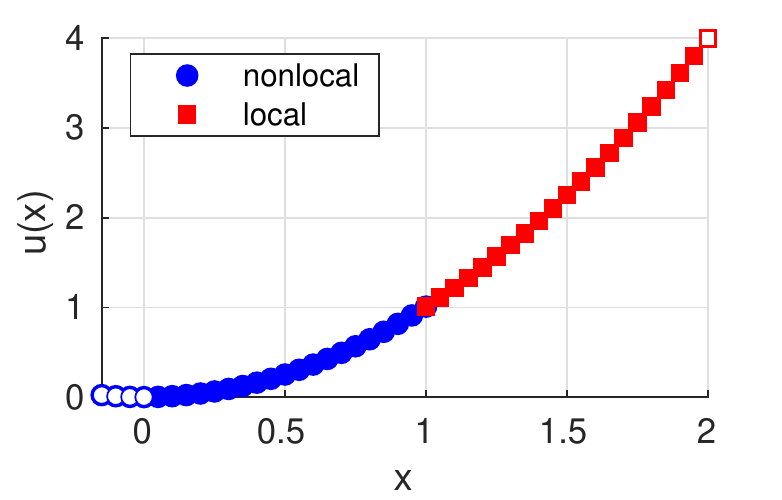}}
\caption{Patch tests for the non-overlapping partitioned procedure: linear (left) and quadratic (right). Nodes in the nonlocal sub-domain are represented by blue filled circles, whereas nodes in the left nonlocal boundary  are empty blue circles. Nodes in the local sub-domain are represented by red filled squares, whereas the node on the right local boundary is an empty red square.}
\label{Fig:Patch Tests PP}
\end{center}
\end{figure}

For LtN coupling problems with general geometry and heterogeneous material properties, deriving the analytical expression of the reduction factor is generally not straightforward, and typically a numerical approximation has to be considered. In \cite{you2019coupling}, LtN coupling in a diffusion problem for the non-overlapping case is studied. Using Robin transmission conditions, the authors developed a stable explicit partitioned procedure, where the optimal Robin coefficient is obtained by minimizing the magnitude of the maximum eigenvalue of the discretized coupled system. When the time step  satisfies the CFL condition $\Delta t\leqslant C h^2$, where $h$ is the spatial discretization size and $C$ is a constant, the resulting partitioned procedure with optimal Robin coefficient is robust and capable to handle sub-domains with complicated geometries, as shown in Figure \ref{Figure:cross}. The numerical results on patch tests are shown in Figure \ref{Fig:Patch Tests PP}, where both the local and nonlocal sub-domains are two-dimensional squares and the interface $\Gamma_v$ is a straight line segment. 
In Figure \ref{Fig:Patch Tests PP}, solutions along the domain center line are reported. For both linear and quadratic solutions, the numerical solution is machine accurate.

\section{Atomistic-to-Continuum type coupling}\label{sec:AtC}

\subsection{Arlequin method}\label{subsec:Arlequin}
The Arlequin method 
\cite{dhia1998multiscale,ArlequinDhia1999} 
is a general energy-based engineering design tool for multiscale modeling which superimposes multiple models and glues them to each other while partitioning and weighting the energies. In this section, we review the Arlequin approach for LtN coupling, following \cite{HanLubineau2012}. In that work, the authors employ the nonlocal elasticity model developed in \cite{Paola2009,Paola2010}. Here, for consistency with the rest of the paper, we adapt the formulation from \cite{HanLubineau2012} to couple the linear bond-based peridynamic model \eqref{eq: bond-based-PD-force-state-linear} and classical linear elasticity (see \eqref{eq: CE-Navier-Cauchy}). 
%
%
For state-based peridynamics, the Arlequin method has been applied in 
\cite{ArlequinWang2019}.

\subsubsection{Mathematical formulation}
We refer to Figure~\ref{fig:decomposed-domains} (bottom): the 
domain $\ooomg$ is decomposed into $\oomg_{l}$ and $\ooomg_{nl}$ with an overlapping region $\overline{\omg}_o: =\oomg_{l}\cap\ooomg_{nl}$,  
whose volume has to be sufficiently large to ensure stability \cite{Bauman2008,HanLubineau2012}. 
The strain energy density of the system is defined in the following way: it equals the nonlocal energy density in $\omgn\setminus\omgo$ and the local energy density in $\omgl\setminus\omgo$. 
%
%
In the overlapping region, $\Omega_{o}$, the strain energy density is defined as a weighted combination of the strain energy densities of the two models
using complementary weight functions. We introduce below, in a three-dimensional setting, the weak form of the system for the Arlequin formulation.

The functional spaces for the local model are defined as
\[
\begin{split}
&W_1:=\left\{\ub\in H^1\left(\oomg_l\right): \ub(\xb)=\gb(\xb) \text{ on }\gammad\right\},\\
&W^0_1:=\left\{\vb\in H^1\left(\oomg_l\right): \vb(\xb)=\bm{0} \text{ on }\gammad\right\}.
\end{split}
\]
Given a function pair $\left(\ub_1,\vb_1\right)\in W_1\times W_1^0$, we define a local weighted bilinear form: 
\begin{equation}\label{weight_local_bilinear}
\begin{split}
    a_1(\ub_1,\vb_1):=&\int_{\omgl}\alpha_1(\xb)\frac{4E}{5}\big\{
   (\nabla\cdot \ub_1(\xb))\,(\nabla\cdot \vb_1(\xb))\\
    &\qquad\qquad\quad+\frac{1}{2}\left(\nabla \ub_1(\xb)+\nabla \ub_1^T(\xb)\right):
    \left(\nabla \vb_1(\xb)+\nabla \vb_1^T(\xb)\right)
    \big\}d\xb.
\end{split}
\end{equation}
%
%
For the nonlocal model, we define the functional space $W_2 : = \mcS^{\rm{PD}}(\ooomg_{nl})$. 
%
%
Given a function pair $\left(\ub_2,\vb_2\right)\in W_2\times{W_2}$,
%
%
%
we define a nonlocal weighted bilinear form: 
\begin{equation}\label{weight_nonlocal_bilinear}
  \begin{split}
    a_2(\ub_2, \vb_2)
    :=&\frac{1}{4}\int_{\omgn} 
    \int_{B_{\horizon}(\xb)}\alpha_2(\xb,\xbp) \lambda(\|\xbp-\xb\|) \left(\ub_2(\xbp) - \ub_2(\xb)\right)^T\\
    &\qquad \qquad
   \ \big((\xbp-\xb)\otimes(\xbp-\xb) \big)  \left(\vb_2(\xbp) - \vb_2(\xb)\right) d\xbp d\xb.
\end{split}
\end{equation}
Here, $\alpha_1(\xb)$ and $\alpha_2(\xb, \xbp)$ are complementary (see below) scalar-valued weight functions, which can be defined via a blending function $\beta(\xb)$, such as the one introduced in \eqref{general_def_blend_fnc} and illustrated in Figure~\ref{fig: blending function2}(a).  
%
%
We take $\alpha_1(\xb)=1-\beta(\xb)$. In contrast to $\alpha_1(\xb)$, which is a one-point function,  $\alpha_2(\xb, \xbp)$ is a two-point function so that its definition is not unique. Nevertheless, as discussed in \cite{Bauman2008,Chamoin2010}, to reduce spurious effects  $\alpha_2(\xb, \xbp)$ can chosen as a symmetric function (see \cite{HanLubineau2012}):
%
%
  \begin{equation}\label{weight_nonlocal}
 \alpha_2(\xb,\xbp) =\beta\left(\frac{\xb+\xbp}{2}\right)\quad \forall \xb, \xbp\in\Omega.
 \end{equation}
 Another standard choice of $\alpha_2$ can be defined in an average sense \cite{Chamoin2010,Prudhomme2009}.

Next, we define the weak formulation of local and nonlocal body forces, respectively. For $\vb_1\in W_1^0$, the local term is defined as
\begin{equation}\label{weight_local_load}
l_1(\vb_1):=\int_{\Omega_{l}}\alpha_1(\xb) \bb(\xb)\cdot \vb_1(\xb)\,d\xb,
\end{equation}
and, for $\vb_2\in W_2$, the nonlocal term is defined as
\begin{equation}\label{weight_nonlocal_load}
l_2(\vb_2):=\int_{\Omega_{nl}}(1-\alpha_1(\xb)) \bb(\xb)\cdot \vb_2(\xb)\,d\xb.
\end{equation}
In order to couple the displacement fields of the local and nonlocal models 
in $\Omega_o$, a weak compatibility between the kinematics of both models using Lagrange multipliers is enforced. The result is the following saddle point problem:
\begin{align}
&\text{Find }(\ub_1,\ub_2,\phib)\in W_1\times W_2\times H^1(\Omega_o) \quad \mbox{such that} \nonumber\\
& a_1(\ub_1,\vb_1)+C(\phib,\vb_1)=l_1(\vb_1),\label{weak_formula_local}\\
& a_2(\ub_2,\vb_2)-C(\phib,\vb_2)=l_2(\vb_2),\label{weak_formula_nonlocal}\\
& C(\psib,\ub_1-\ub_2)=0,\label{weak_formula_Lagrange}\\
&\text{for all }
(\vb_1,\vb_2,\psib)\in W^0_1\times W_2\times H^1(\Omega_o).\nonumber
\end{align}
The bilinear form $C(\cdot,\cdot)$ describes the coupling and is defined, for any $(\psib,\vb)\in H^1(\Omega_o)\times H^1(\Omega_o)$, as
\begin{equation}\label{weak_coupling}
C(\psib,\vb):=\int_{\Omega_o}\left(\kappa_0\psib\cdot \vb
+\kappa_1\epsb(\psib):\epsb(\vb)\right)d\xb,
\end{equation}
where $\epsb(\cdot)$ denotes the \textit{infinitesimal strain tensor} 
$
\epsb(\vb):=\frac{1}{2}\left(\nabla \vb+\nabla\vb^{T}\right).
$
The coefficients $0\leqslant \kappa_0, \kappa_1 \leqslant 1$ are non-negative coupling parameters. For example, $(\kappa_0, \kappa_1)=(1,0)$ defines the $L^2$-norm coupling; $(\kappa_0,\kappa_1)=(0,1)$ defines the $H^1$-seminorm coupling; and $(\kappa_0, \kappa_1)=(1,1)$ is the $H^1$-norm coupling  \cite{ArlequinDhia2005,HanLubineau2012,Belytschko2007a,Belytschko2009a}.

Once the mixed formulation \eqref{weak_formula_local}-\eqref{weak_formula_Lagrange} is solved, we need to reconstruct 
a displacement field $\ub$ for the entire domain $\Omega$, 
because $\ub_1$ and $\ub_2$ are only defined in $\Omega_{l}$ and $\ooomg_{nl}$, respectively. One convenient reconstruction option is given as follows:
\begin{equation}\label{Arlequin_solution}
\ub(\xb)=\begin{cases}
\ub_1(\xb)\qquad
\xb\in\Omega_{l}\setminus\Omega_o,\\
\ub_2(\xb)\qquad
\xb\in\Omega_{nl}\setminus\Omega_o,\\
\alpha_1(\xb)\ub_1(\xb)+(1-\alpha_1(\xb))\ub_2(\xb)\qquad
\xb\in\Omega_o.
\end{cases}
\end{equation}

\paragraph{Properties.}
Even though a mathematical analysis and a rigorous study of the
influence of parameters (such as weight functions or size of the overlapping region) of the Arlequin method for mechanics are not available, several references in the literature address properties of 
this 
technique when applied to other coupling problems
(see, e.g., \cite{ArlequinDhia2001,Belytschko2007a,Bauman2008,Chamoin2010}). Based on these references, we list here some general properties: 
\begin{itemize}
\item If the coupling parameters defined in \eqref{weak_coupling} satisfy $\kappa_0\geqslant 0$ and $\kappa_1>0$, and the overlapping region $\Omega_o$ is sufficiently large, then the corresponding saddle point problem \eqref{weak_formula_local}-\eqref{weak_coupling} is well-posed. However, for the $L^2$-norm coupling 
$(\kappa_0,\kappa_1)=(1,0)$, the well-posedness is unclear.
\item Since the linear patch test is only weakly imposed via \eqref{weak_formula_Lagrange}, ghost forces and spurious effects cannot be completely removed in the overlapping region.
\item The energy equivalence of this formulation depends on the choice of the weight functions $\alpha_1$ and $\alpha_2$, the size of overlapping region $\omgo$, and the choice of the Lagrange multiplier space \cite{Chamoin2010,Belytschko2007a}. In some cases, the Arlequin method is not equivalent neither to the local model nor the nonlocal model, even when homogeneous deformations are assumed. In fact, as pointed out in \cite{Chamoin2010}, with inappropriate choices of $\alpha_1$, $\omgo$, and $(\kappa_0,\kappa_1)$, the Arlequin formulation may be not coercive.
\end{itemize}

\paragraph{The time-dependent problem.}
The Arlequin method has been used in many time-dependent applications, and the extension is straightforward, see, e.g., \cite{Bauman2008,Winker2013,ArlequinWang2019}. 
In particular, applications to dynamic LtN coupling mechanical problems appear in \cite{ArlequinWang2019}. 
However, such extension consists in solving the saddle point problem \eqref{weak_formula_local}-\eqref{weak_formula_Lagrange} at every time step, 
which limits its applicability due to high computational cost.

\subsubsection{Applications and results}
In this section, as done in others below, we provide references of applications of the Arlequin method as the reproduction of numerical results is non-trivial.

In \cite{HanLubineau2012}, the authors first apply the $H^1$-norm coupling 
and piece-wise linear weight functions to study a two-dimensional cantilever beam 
of isotorpic homogenous material. 
%
%
%
Their results show that the accuracy of Arlequin solutions is comparable to that of the fully nonlocal elasticity model.
Next, they test a static cracked square plate using various options of $(\kappa_0, \kappa_1)$. When an $H^1$-norm coupling 
with piece-wise linear weight functions is used, the strain distribution from the Arlequin approach agrees with the strain field, especially near the crack-tip, computed with the fully nonlocal model. 

In \cite{Prudhomme2008}, the authors use the Arlequin method to investigate a one-dimensional problem that consists of a collection of springs that exhibit a localized defect, resulting in a sudden change in the spring properties. They test the method and study its accuracy for several choices of coupling parameters. They prove that the Arlequin formulation is well-posed with both $H^1$-seminorm and $H^1$-norm couplings. Their numerical results also indicate that the method is sensitive to the location and size of the overlapping region, and they propose to utilize adaptive strategies, based on a posteriori error estimates, to identify them.

\subsection{Morphing method}\label{subsec:morphing}

The morphing method for LtN coupling was developed in \cite{azdoud2013morphing,azdoud2014morphing,han2016morphing,lubineau2012morphing} based on a blending approach to morph the material properties of local and nonlocal sub-domains.
The coupling formulation consists of a single unified model obtained by a transition (morphing) from local to nonlocal descriptions. More specifically, a hybrid model is introduced in the transition region or morphing zone, $\Omega_t$ (see Figure~\ref{fig:Omegab-domain}), whose constitutive law changes gradually from a local to a nonlocal response. As a result, a single model with evolving material properties is defined on the whole domain and the equivalence of the energy of the system with the fully nonlocal energy is enforced under homogeneous deformations in the morphing zone \cite{lubineau2012morphing}.

\subsubsection{Mathematical formulation}
We describe the morphing method for the coupling of the linear isotropic bond-based peridynamic model~\eqref{eq: bond-based-PD-force-state-linear} and the corresponding classical linear elasticity model \eqref{eq: CE-Navier-Cauchy}, following \cite{lubineau2012morphing}. 
We note, however, that this technique has also been applied to more complex material models, including linear anisotropic bond-based peridynamic models \cite{azdoud2013morphing} and a state-based peridynamic model  \cite{han2016morphing}.
%

The strain energy density 
of the linear bond-based peridynamic model is given by
\begin{equation}\label{morphing_nonlocal}
\begin{split}
 W^{nl}(\xb)=&\dfrac{1}{4}\int_{B_\delta(\mathbf{0})}\lambda(\|\xib\|) (\ub(\xb + \xib) - \ub(\xb))^T (\xib\otimes \xib ) (\ub(\xb + \xib) - \ub(\xb)) d\xib,
\end{split} 
\end{equation}
whereas the corresponding local strain energy density 
is given by
%
\begin{equation}\label{morphing_local}
\begin{split}
W^l(\xb)=&
{\frac{4E}{5}\big\{
   (\nabla\cdot \ub(\xb))\,(\nabla\cdot \ub(\xb))}\\
    & {\qquad
    +\frac{1}{2}\left(\nabla \ub(\xb)+\nabla \ub^T(\xb)\right):
    \left(\nabla \ub(\xb)+\nabla \ub^T(\xb)\right)\big\} }\\
=&\dfrac{1}{2}\epsb(\xb):\mathcal{C}^l:\epsb(\xb),
\end{split}
\end{equation}
where $\mathcal{C}^l$ is the 
fourth-order elasticity tensor and $\epsb = \dfrac{1}{2}(\nabla \ub+\nabla\ub ^T)$ is the infinitesimal strain tensor. When considering an infinitesimal homogeneous deformation, we can define a local stiffness tensor $\mathcal{C}^0$ from the nonlocal model \eqref{morphing_nonlocal} such that the strain energy density of the resultant local model is equal to the strain energy density of the nonlocal model, i.e.,
$$
W^{nl}(\xb)\approx \dfrac{1}{2}\epsilon(\xb):\mathcal{C}^0:\epsilon(\xb),$$
where $\mathcal{C}^0$ is given by \eqref{eq: fourth-order elasticity tensor bond-based}.
%
To define the morphing model, due to consistency requirements on the energy densities, we assume $\mathcal{C}^0 = \mathcal{C}^{l}$. 
%

Given a blending or morphing function $\beta$, such as the one introduced in \eqref{general_def_blend_fnc} and illustrated in Figure~\ref{fig: blending function2}(b), 
%
the morphing model is fully defined by the following strain energy density:
\begin{align}
\nonumber &W^m(\xb)=\dfrac{1}{2}\epsilon(\xb):\mathcal{C}(\xb):\epsilon(\xb)\\
\nonumber &+\dfrac{1}{4}\int_{B_\delta(\mathbf{0})}\lambda(\|\xib\|)\dfrac{\beta(\xb + \xib)+\beta(\xb)}{2} (\ub(\xb + \xib) - \ub(\xb))^T (\xib\otimes \xib ) (\ub(\xb + \xib) - \ub(\xb)) d\xib,   
\end{align}
where
$$
\mathcal{C}(\xb):=(1-\beta(\xb))\mathcal{C}^l+\int_{B_\delta(\mathbf{0})}\lambda(\|\xib\|)\dfrac{\beta(\xb)-\beta(\xb + \xib)}{4} \xib\otimes\xib\otimes\xib\otimes\xib  d\xib.
$$
Note that, similarly to \eqref{eq: fourth-order elasticity tensor bond-based}, $\mathcal{C}(\xb)$ is a fully-symmetric fourth-order tensor. 
%
%
%

\paragraph{Properties.} The morphing method has the following properties:
\begin{itemize}

\item For $\xb\in \omgl$,  $\mathcal{C}(\xb)=\mathcal{C}^l$ and $W^m(\xb)=W^l(\xb)$.

\item For $\xb\in \omgn$, $\mathcal{C}(\xb)=0$ and $W^m(\xb)=W^{nl}(\xb)$.

\item This method does not pass the {linear patch test}. In fact, $u^{\rm lin}$ (see Definition \ref{def:patch-test}) does not satisfy the equilibrium equation throughout the morphing zone and a nonzero ghost force density arises \cite{lubineau2012morphing}. However, these ghost forces can be approximately corrected using exactly the same deadload correction approach used in AtC coupling methods \cite{Shenoy1999b}. Besides, the ghost force intensity decreases when using smoother morphing functions $\beta$ or sufficiently large morphing zones. Moreover, ghost forces are localized to the morphing zone and vanish when $\delta\rightarrow 0$. 
\item For homogeneous deformations, the {strain energy density} is equivalent to {both strain energy densities of the local and nonlocal models}.  Therefore, this method is considered { energy preserving} under homogeneous deformations \cite{lubineau2012morphing}.
\item Even though there are a few theoretic studies regarding the morphing  method for LtN coupling, this method has been studied as a type of blending for AtC coupling \cite{li2012positive}. The corresponding operator is coercive with respect to the nonlocal energy norm with smooth morphing function $\beta$ and sufficiently large morphing zone $\omg_t$ \cite{li2012positive}.
\end{itemize}
\paragraph{The time-dependent problem.} 

The extension of the morphing method to time-dependent problems is straightforward, even though the implementation of the model has been only demonstrated in static/quasi-static  problems. 



\subsubsection{Applications and results}
In this section, we provide references to applications of the morphing method as the reproduction of numerical results is non-trivial. In \cite{lubineau2012morphing} this strategy is applied to couple linear
 bond-based peridynamics with classical linear elasticity for isotropic materials. 
%
The authors investigate the ghost force intensity when using different 
morphing 
functions. They perform a one-dimensional analysis followed by numerical studies in one and two dimensions. The results suggest that the ghost forces are localized to the morphing zone and a smoother morphing function $\beta$ reduces the maximum relative ghost forces \cite{lubineau2012morphing}. Follow-on two-dimensional simulations for a cracked plate under both traction and shear 
demonstrate the effectiveness of the method, compared to a fully peridynamic simulation. 
%
%
%
%
Later on, in \cite{azdoud2014morphing}, the morphing method is combined with an adaptive algorithm that updates the nonlocal sub-domain based on damage progression. The resulting coupling framework is applied to three-dimensional quasi-static problems. 
In \cite{han2016morphing}, the method is further extended to couple a linearized state-based peridynamic model and the corresponding classical linear elasticity model. For anisotropic materials, in \cite{azdoud2013morphing}, the authors introduce anisotropic nonlocal models based on spherical harmonic descriptions, and present three-dimensional results.

\subsection{Quasi-nonlocal method}\label{subsec:quasinonlocal}
The quasi-nonlocal (QNL) method is an energy-based coupling approach introduced in the context of AtC coupling. 
This method redefines the nonlocal energy via a ``geometric reconstruction'' scheme  in the transition region and local sub-domain of a LtN coupling configuration,  
in such a way that the method is {linearly} patch-test consistent \cite{Luskin2013a,Shimokawa:2004,E:2006}. We point out that the idea of ``geometric reconstruction'' is not limited to AtC coupling {of solids;} 
similar coupling strategies in the literature have been applied, for example in computational fluid dynamics (see, e.g., the review papers \cite{AidunClausen2010,LiLiu2002}).
%
%
Here, we focus on the QNL method for 
LtN coupling of one-dimensional diffusion models, following \cite{DuLiLuTian2018,XHLiLu2017}.

\subsubsection{Mathematical formulation}
We refer to Figure \ref{fig:blended-domains} (top): without loss of generality, consider the domain $\ooomg=[-1-\delta,1]$. %
The domain is decomposed into four  
disjoint sub-domains: 
{$\ooomg=\omgd\cup\omg_{nl}\cup\Omega_t\cup\omgl$}, which include the left physical nonlocal boundary  $\omgd=(-1-\delta, -1)$, the nonlocal sub-domain  $\omgn=(-1, \;x^*)$, the transition region $\omgt=(x^*,\; x^*+\delta)$, and the local sub-domain $\omgl=(x^*+\delta, \;1)$. 
The right physical local boundary is $\Gamma_p=\{1\}$. 
%
Note that the interface between the nonlocal sub-domain and the transition region occurs at $x^*$ satisfying
$x^*\in(-1+2\delta, 1-2\delta)$, 
and the transition region, $\omg_t$, has thickness $\horizon$. 
%
%

The crucial step in the QNL formulation is the ``geometric reconstruction'' \cite{Shapeev2012a,Shimokawa:2004,E:2006} of the directional distance $u(\xbp)-u(\xb)$ in the definition of the energy. Because of difficulties arising from reconstructing geometries with corners in high dimensions, we limit the discussion to the one-dimensional case.

Recall that the one-dimensional nonlocal diffusion energy density associated with the bond $\xi = x'-x$ is 
\begin{equation}\label{eq:nl-energy-density}
\frac{1}{2}\gamma_\delta(x'-x)\left(u(x')-u(x)\right)^2,
\end{equation}
where we assume a radially symmetric nonlocal diffusion kernel, i.e., $\gamma_{\delta}(\xi) = \gamma_{\delta}(|\xi|)$. 
%
In the QNL model, such bond energy density is modified when the 
bond 
is entirely located in the local sub-domain. Specifically, the nonlocal energy is redefined by substituting the directional distance $\left(u(x')-u(x)\right)$ with a path integral, such that the local energy density, $\frac{1}{2}|u'(x)|^2$, is equivalent to the nonlocal one for sufficiently smooth~$u$. Thus, we have that \eqref{eq:nl-energy-density} is replaced by
%
\begin{equation}\label{eq:geom_reconstruct}
{\frac{1}{2}}\gamma_{\delta}(x'-x)\int_{0}^{1}\left|u'\big(x+t(x'-x)\big)\right|^2|x'-x|^2dt.
\end{equation}
%
The combined total energy of the QNL model with interface at $x^*$ is
\begin{align}\label{eq:qnl_energy}
E^{\rm{QNL}}_\delta(u)
=& \frac{1}{4} 
\iint\limits_{x\leqslant {x^*} \text{ or }x'\leqslant {x^*}
} \gamma_{\delta}(\abs{x'-x})\left({u(x')}-u(x)\right)^2\, dx'dx \\
& +\frac{1}{4}
\iint\limits_{x> {x^*}\text{ and }x'> {x^*}
} \gamma_{\delta}(\abs{x'-x} )\int_{0}^1 {\abs{u'(x+t(x'-x))}}^2 |x'-x|^2 dt \, dx'dx.\nonumber
\end{align}
%
%
The definition of the QNL coupling operator $\mcL^{\rm{QNL}}$ is obtained by taking the negative first variation of the total energy \eqref{eq:qnl_energy}. We split it into three parts \cite{DuLiLuTian2018}:
\begin{itemize}
\item[] {\bf I.} 
For $x\in \omgn$: 
   \begin{align}\label{force_case1}
   \mcL^{\rm{QNL}}u(x)
   =& \int_{x'
   \in \ooomg
   } \gamma_{\delta} (\abs{x'-x})\left(u(x')-u(x)\right) dx'.
   \end{align}
\item[] {\bf II.} 
For $x\in\omgt$: 
   \begin{align}\label{force_case2}
   \mcL^{\rm{QNL}}u(x) =& \int_{{x'<x^*} }
   \gamma_{\delta}(\abs{x'-x})\left(u(x')-u(x)\right)dx'+ \left(\omega_\delta(x) u'(x)\right)',\,
   \end{align}
   where $\omega_{\delta}(x)$ is defined as
   \begin{equation*}
  {
   \omega_\delta(x):=\int_{0}^{1} dt \int_{|\xi|<{\frac{|x-x^*|}{t}} } |\xi|^2\gamma_{\delta}(|\xi|) \, d\xi.
   }
   \end{equation*}
%
\item[] {\bf III.} 
For $x\in\omgl$:
   \begin{align}\label{force_case3}
   \mcL^{\rm{QNL}} u(x)=& \frac{\partial^2 u}{\partial x^2}(x) .
   \end{align}
\end{itemize}

\paragraph{Properties.}
We summarize the properties of the QNL method, following \cite{DuLiLuTian2018}:
\begin{itemize}
\item The QNL operator is {self-adjoint} (i.e., symmetric) as it is derived from a combined total energy \cite{CaffarelliChanVasseur:11}.
 \item The method passes the {linear patch-test}.
\item The operator is positive-definite and, hence, it is energy stable {with respect to the QNL energy defined in \eqref{eq:qnl_energy}} as well as $L^2$ stable {with respect to the  $L^2(\ooomg)$ norm}.
\item The formulation satisfies the weak maximum principle and therefore it is mass-conserving.
\item The method is {asymptotically compatible} and the coupled solution converges as $O(\delta)$ to the corresponding local solution.
\item The total combined energy is equivalent to either the fully nonlocal or the fully local energy up to linear functions.
\end{itemize}

\paragraph{The time-dependent problem.}
The extension of the QNL method to time-dependent problems is simple and straightforward. This follows from the fact that the QNL operator $\mcL^{\rm{QNL}}$ is coercive on the entire domain $\ooomg$ and, hence, it is stable and monotonic.  

\subsubsection{Applications and results}
\begin{figure}[H]
\begin{center}
\subfigure[Linear Patch Test]{
\includegraphics[height=4cm, width=0.4\textwidth]{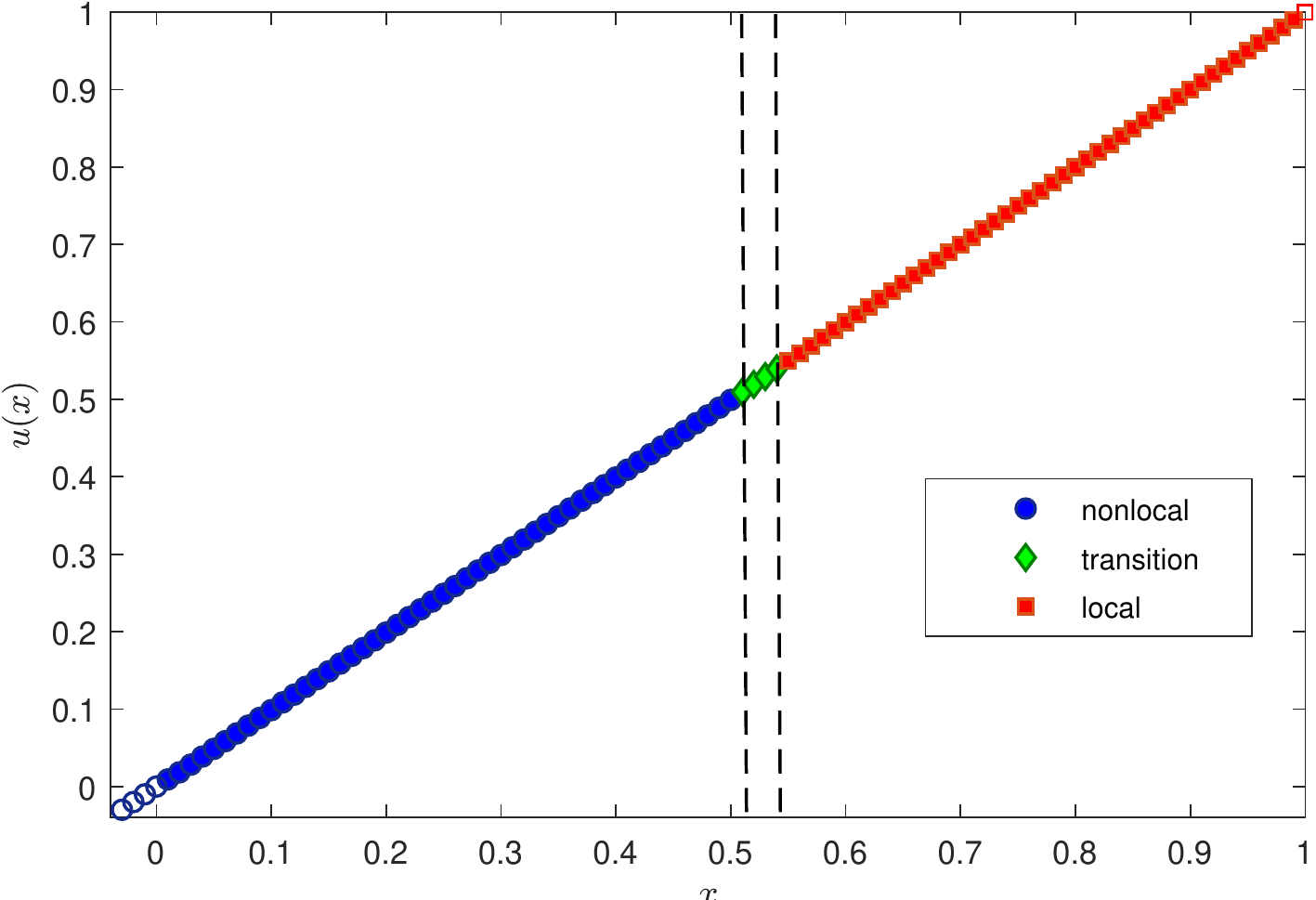}}
\qquad
\subfigure[Quadratic Patch Test]{
\includegraphics[height=4cm, width=0.4\textwidth]{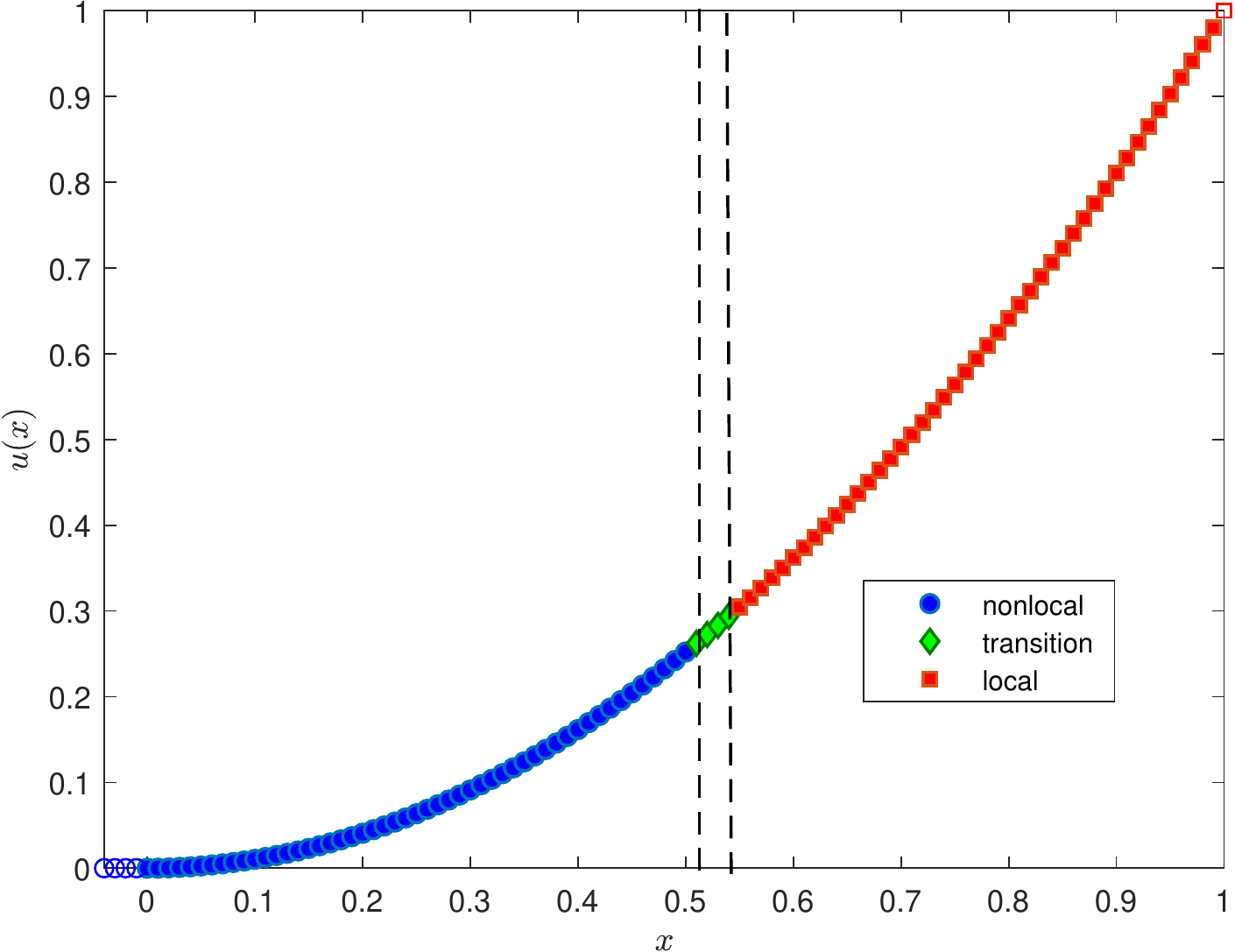}}
\vspace{-0.14 in}
\caption{Patch tests for the quasi-nonlocal (QNL) method: linear (left) and quadratic (right). Nodes in the nonlocal sub-domain are represented by blue filled circles, whereas nodes in the left nonlocal boundary are empty blue circles. Nodes in the local sub-domain are represented by red filled squares, whereas the node on the right local boundary is an empty red square. Green diamonds represent nodes in the transition region described by the QNL model.
%
}
\label{Fig:Patch_Tests_QNL}
\end{center}
\end{figure}

Consider the one-dimensional domain $\ooomg=[-\delta,\;1]$ with interface at $x^*=1/2$.
The domain is decomposed into the subdomains $\omgd=(-\delta, 0)$, $\omgn=(0, \;1/2)$, $\omgt=(1/2,\; 1/2+\delta)$, and $\omgl=(1/2+\delta, \;1)$. 
%
%
Results are obtained with the first-order asymptotically compatible finite difference method introduced in \cite{DuLiLuTian2018}. We report results of patch tests for both linear and quadratic solutions. 
In this case, the horizon is taken as $\horizon = 0.04$ and the grid size 
is $h = 0.01$. Results are presented in Figure~\ref{Fig:Patch_Tests_QNL} and demonstrate that, even though theoretically the QNL method is only linearly path-test consistent, numerically both linear and quadratic patch tests are well passed.
%

\begin{figure}[H]
\begin{center}
 \subfigure[Nonlocal-local coupling]{\includegraphics[height=3cm, width=0.35\textwidth]{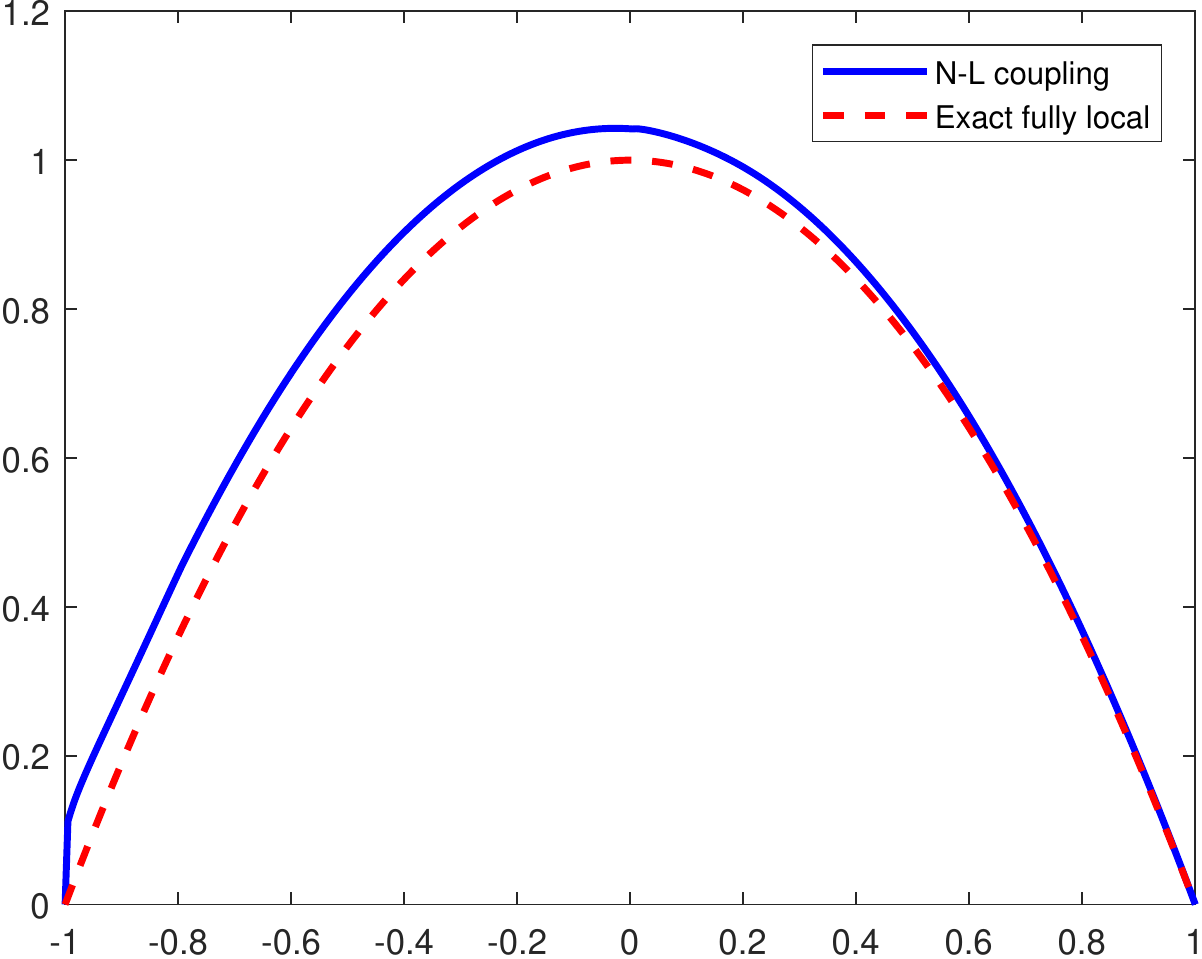}}
\qquad
\subfigure[Local-nonlocal-local coupling]{\includegraphics[height=3cm, width=0.35\textwidth]{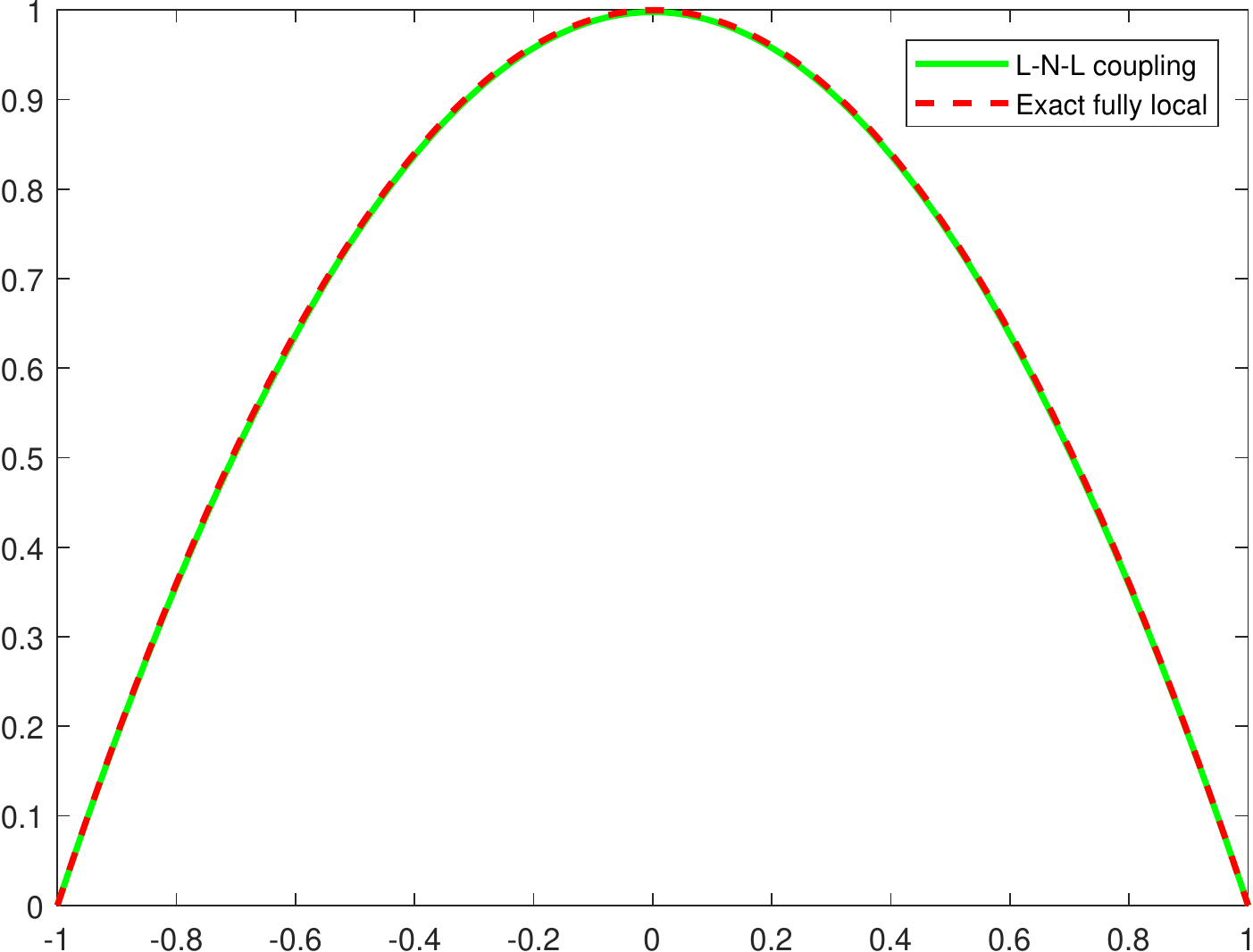}}
\vspace{-0.14 in}
\caption{
Application of the QNL method to impose classical local boundary conditions in a nonlocal problem. Left: nonlocal-local (N-L) coupling with a classical Dirichlet boundary condition extended to the left physical nonlocal boundary and a classical Dirichlet boundary condition imposed on the right local boundary. 
Right: local-nonlocal-local (L-N-L) coupling with classical Dirichlet boundary conditions imposed on both left and right local boundaries.}
%
%
%
\label{Fig:LtoNtoL_demo}
\end{center}
\end{figure}

{We now discuss the application of the QNL method in \cite{DuLiLuTian2018} to remove surface effects.} 
As mentioned in Section~\ref{sec:intro}, volumetric constraints for nonlocal models have to be prescribed in a layer surrounding the domain, where data are not available or are hard to access. Consequently, an ad hoc  
treatment of nonlocal boundaries 
often causes unphysical surface effects. 
%
In \cite{DuLiLuTian2018}, the authors apply the QNL method to allow the prescription of classical local boundary conditions in a nonlocal problem. Specifically, they 
address the nonlocal boundary issue by coupling nonlocal and local models, being the latter placed in the nonlocal boundary region. 
%
%
%
Figure \ref{Fig:LtoNtoL_demo}(b) shows the solution of a local-nonlocal-local (L-N-L) coupling with boundaries being treated in a fully classical way; whereas Figure~\ref{Fig:LtoNtoL_demo}(a) shows the solution of a nonlocal-local (N-L) coupling by simply imposing a nonlocal Dirichlet type boundary condition $u(x)=u(-1)$ in the left physical nonlocal boundary 
$\omgd$.  
In those two examples, the nonlocal diffusion kernel is chosen as  $\gamma_{\delta}(x)=\frac{3}{\delta^3}\chi_{(-\delta,\delta)}(x)$ 
and the source is $f\equiv 1$. 
The L-N-L coupling has two interfaces located at $x^a=-\frac{1}{2}$ and $x^b=\frac{1}{2}$; 
the N-L coupling has one interface at $x^*=0$.
%
The horizon is  $\delta=0.2$ and the mesh size is $h=1/800$. 
%
%
%
We observe that the N-L coupling displays unphysical surface effects 
on the nonlocal side, whereas the L-N-L coupling eliminates those effects. 
%
%
Consequently, L-N-L coupling provides 
a way of imposing 
classical local boundary conditions in a nonlocal problem and helps avoiding 
the need of 
complicated extensions of local boundary conditions to nonlocal boundaries.
\begin{remark}\label{rmk:qnl}
Finding an analytic geometric reconstruction formula becomes more difficult in higher dimensions. The nonlocal  neighborhood, $B_\horizon(\xb)$, becomes
a disk (in two dimensions) or a ball (in three dimensions), making the intersections with the interface highly more complex. As a result, it is difficult to identify an explicit expression for the QNL  
method in a general setting (see discussions from \cite{Shapeev2012a}); even for a simple interface, such as a sphere, the exact formula is not yet available. In \cite{JiangShen2018}, the authors propose to adopt a least-square fitting procedure  to find an approximation of this reconstruction; however, due to the inexactness of the approximation, the properties listed above no longer hold.
\end{remark}

\subsection{Blending method}
\label{subsec:blending}
In this section, we discuss a force-based blending approach originally presented for one-dimensional linear bond-based peridynamic models in~\cite{Seleson2013CMS} and then extended to general  bond-based peridynamic models in higher-dimensions in~\cite{seleson2015concurrent}. Force-based blending is based on weighting governing equations via the introduction of a blending function. As discussed in~\cite{badia2008atomistic}, external or internal blending is possible. The former does not change the definitions of the reference models to be coupled,  
whereas the latter modifies their internal force operators. 
%
In comparison to the morphing method discussed in Section \ref{subsec:morphing},  force-based blending does not necessarily has the physical interpretation of blending material properties. 
Force-based blending has been proposed in AtC coupling (see, e.g.,~\cite{badia2007force,Fish2007,li2012positive,li2014theory}). These AtC coupling methods are based on external blending and simply weight governing equations via the introduction of a blending function. In contrast, the approach from~\cite{Seleson2013CMS} is based on internal blending and arrives at a blended model through derivation from a single reference peridynamic model  by only resorting to assumptions on the deformation. 
Here, as mentioned in Section \ref{sec: overview of methods}, we employ the term  ``blending'' to refer to force-based blending, because related methods, such as  Arlequin and morphing, are called by their specific names.

\subsubsection{Mathematical formulation}

We refer to Figure~\ref{fig:Omegab-domain}: in this decomposition, the nonlocal sub-domain, $\omg_{nl}$, and the local sub-domain, $\omg_{l}$, do not overlap, but are separated by the transition region, $\omgt$. A blending function $\beta (\xb)$, such as the one introduced in \eqref{general_def_blend_fnc} and illustrated in Figure~\ref{fig: blending function2} (b), is used to characterized the different sub-domains. The change in the blending function occurs within the so-called  blending region, $\omg_{b}\subset \omgt$, and normally takes a polynomial shape. Due to nonlocality, however, the influence of the variation in the blending function extends beyond the blending region and affects the entire transition region, which is given by $\omgt = \Omega_b\cup\omgb_{b}$.

The force-based blending approach presented here begins with a single reference bond-based peridynamic model given by \eqref{eq: PD operator bond-based bulk} for all $\xb \in \omg$.  
%
Employing a symmetric combination of the blending function~\eqref{general_def_blend_fnc}, the operator is first split in two contributions, so that one of the two has nonzero support in regions of small smooth deformation thus allowing linearization followed by a Taylor expansion of the displacement at any neighbor point $\xb' \in B_\horizon({\xb})$ around $\xb$ 
up to second order, which leads to an expression connected to the local model. The resulting blended operator is given by:
%
\begin{align}\label{eq: force-based blending operator bulk}
\mcL^{b}\ub(\xb) = &  \int_{B_\horizon({\bf 0})} \left(\frac{\beta(\xb) + \beta(\xb + \xib)}{2} \right) \fb(\etab,\xib) d\xib \nonumber \\
&- \frac{1}{2} \left[ \int_{B_\horizon({\bf 0})} \beta(\xb+\xib)\lambda(\|\xib\|) \xi_i \xi_j\xi_k d\xib\right] \frac{\partial u_j}{\partial x_k}(\xb) {\bf e}_i \nonumber\\
&- \frac{1}{2} \left[ \int_{B_\horizon({\bf 0})} \left( \frac{\beta(\xb) + \beta(\xb+\xib)}{2} \right)\lambda(\|\xib\|) \xi_i \xi_j\xi_k \xi_l d\xib\right] \frac{\partial^2 u_j}{\partial x_k\partial x_l}(\xb) \hat{\bf e}_i\nonumber\\
&+  C_{ijkl}\frac{\partial^2 u_j}{\partial x_k\partial x_l}(\xb) {\bf e}_i,
\end{align}
where $C_{ijkl}$ are the components of the  fourth-order elasticity tensor given by~\eqref{eq: fourth-order elasticity tensor bond-based}. We further used antisymmetry to simplify the term involving first-order derivatives of displacements. 

\paragraph{Properties.} From~\eqref{eq: force-based blending operator bulk}, we have the following properties: 
\begin{itemize}
\item For $\xb \in \omg_{l}$, $\mcL^{b}\ub(\xb) = \mcL^{\rm CE}\ub(\xb)$ (see \eqref{eq: CE-Navier-Cauchy}).
%
\item For $\xb \in \omg_{nl}$, $\mcL^{b}\ub(\xb) = \mcL^{\rm PD}\ub(\xb)$ (see   \eqref{eq: PD operator bond-based bulk}).


%
\item 
By construction, the blended operator~\eqref{eq: force-based blending operator bulk} is equivalent to the reference bond-based peridynamic operator~\eqref{eq: PD operator bond-based bulk} and the classical linear elasticity operator~\eqref{eq: CE-Navier-Cauchy} for quadratic displacements. Therefore, the model automatically passes the linear and quadratic patch tests. Employing a higher-order Taylor expansion up to third order instead would result in a 
blended model which also passes a cubic patch test. 
%
\item 
An error estimate provided in~\cite{seleson2015concurrent} shows that the blended operator~\eqref{eq: force-based blending operator bulk} converges to the classical local operator~\eqref{eq: CE-Navier-Cauchy} with an order of convergence of $\mcO(\horizon)$. Employing a higher-order Taylor expansion up to third order in the derivation of the blended model would increase the order of convergence to $\mcO(\horizon^2)$.
%
\end{itemize}

\paragraph{The time-dependent problem.} The extension of the blending method to time-dependent problems is straightforward. In fact, in \cite{Seleson2013CMS,seleson2015concurrent} the blended model is presented based on equations of motion, even though the implementation of the model is demonstrated in static problems. 

\subsubsection{Applications and results}

\begin{figure}[htbp]
\begin{center}
\subfigure[Linear Patch Test]{
\includegraphics[scale = 0.3]{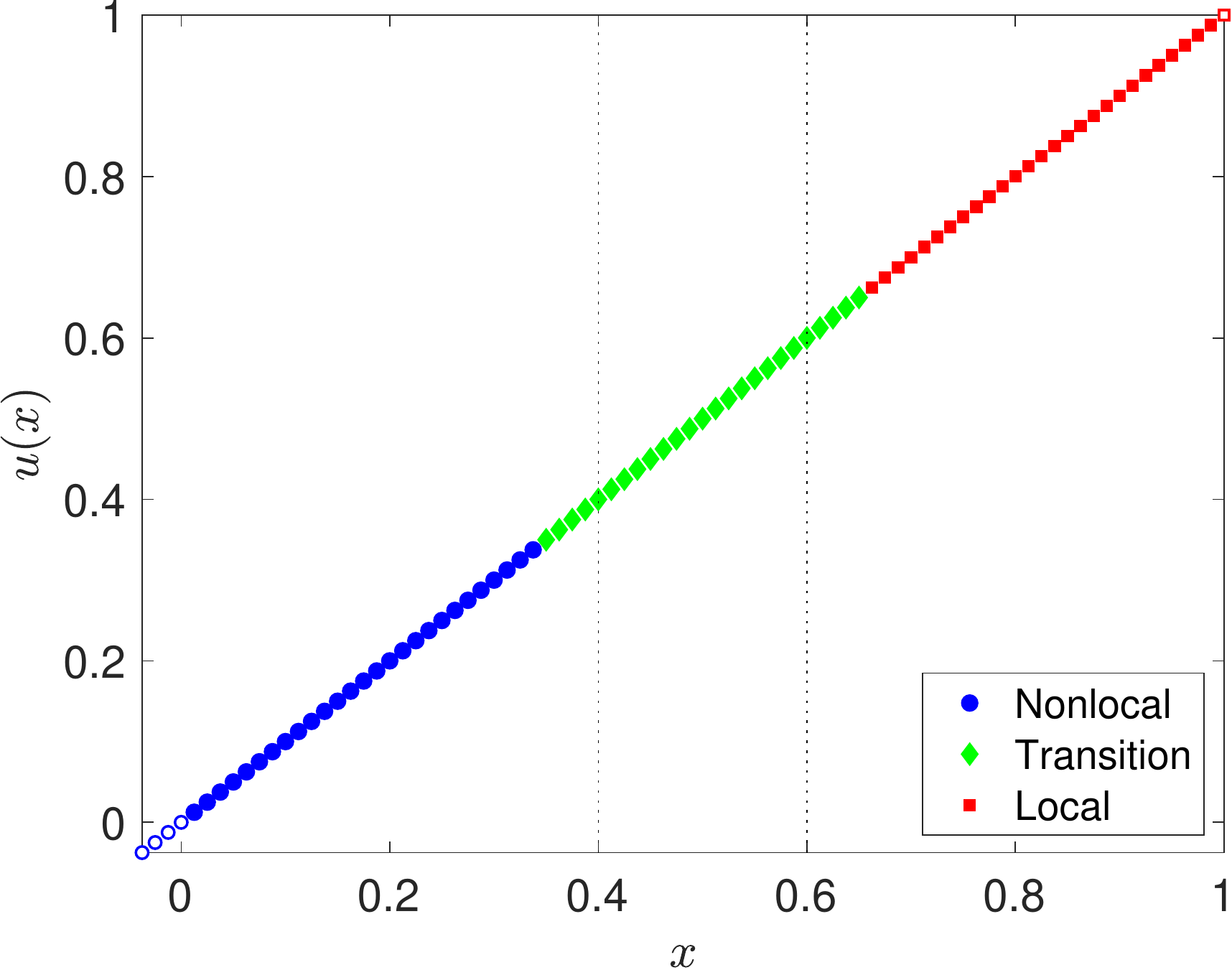}}
\qquad
\subfigure[Quadratic Patch Test]{
\includegraphics[scale = 0.3]{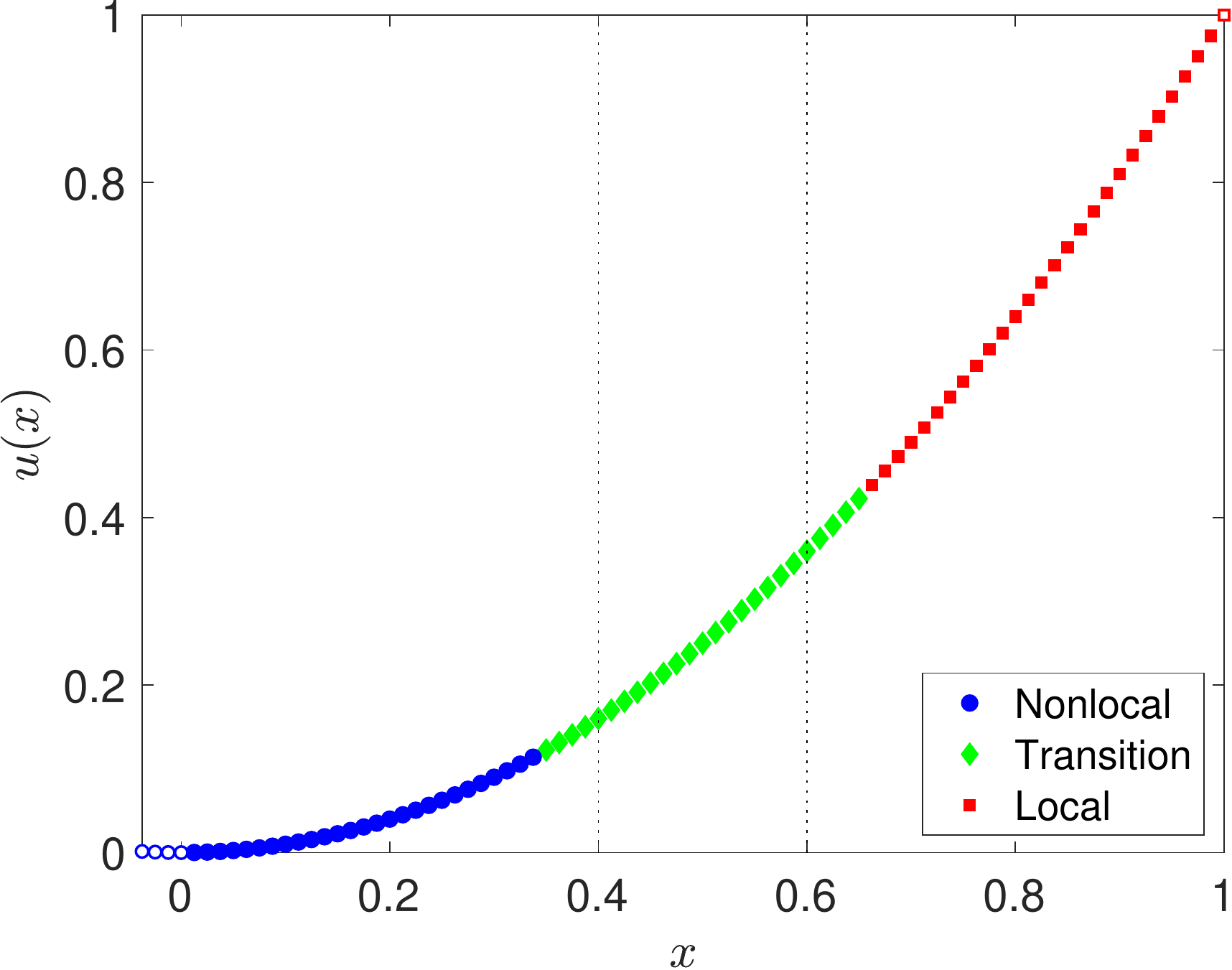}}
\vspace{-0.14 in}
\caption{Patch tests for the blending method: linear (left) and quadratic (right). Nodes in the nonlocal sub-domain 
are represented by blue filled circles, whereas nodes in the left nonlocal boundary 
are empty blue circles. Nodes in the 
local sub-domain 
are represented by red filled squares, whereas the node in the right local boundary is an empty red square. Green filled diamonds represent nodes in the transition 
region described by the blended model.}
\label{Fig:Patch Tests Blended Model}
\end{center}
\end{figure}

Consider a one-dimensional domain $\ooomg = [-\horizon,1]$ decomposed as in Figure~\ref{fig:Omegab-domain}~(top) with $\horizon = 0.05$. The blending region is chosen as $\omg_b  = (0.4,0.6)$. We then have $\omg_p = (-0.05, 0)$, $\omg_{nl} = (0,0.35)$,  $\Omega_t = (0.35,0.65)$, and $\omgl=(0.65,1)$. We solve both a linear patch test with analytical solution $u(x) = x$ and a quadratic patch test with analytical solution $u(x) = x^2$ .The modeling choices and numerical implementation follow~\cite{Seleson2013CMS}; specifically, a linear bond-based peridynamic model is used in $\omg_{nl}$ with a meshfree discretization; classical linear elasticity is employed in  $\omg_{l}$ with a finite difference discretization; and the corresponding blended model given by \eqref{eq: force-based blending operator bulk} is implemented in $\Omega_t$ with a hybrid discretization, which features a meshfree discretization for integral and a finite difference discretization for derivatives. A uniform grid with grid spacing $\Delta x = \horizon / 4$ is utilized. 
A piecewise constant blending function is employed to demonstrate that no regularity is required for the blending function in this case; however, similar results are obtained for  piecewise linear and cubic blending functions. Nonlocal Dirichlet boundary conditions are imposed in the left nonlocal boundary for the peridynamic model, whereas a local Dirichlet boundary condition is imposed on the right local boundary for the local model. 
The results are presented in Figure~\ref{Fig:Patch Tests Blended Model} and demonstrate that the blended model passes the linear and quadratic patch tests. 
A computation of the first derivative of the solution for the linear patch test and the second derivative of the solution for the quadratic patch test indicates that these tests are passed to machine precision. 

%


Additional static numerical examples demonstrating the performance of the blending method appear  for one-dimensional problems in~\cite{Seleson2013CMS} and  for two-dimensional problems in~\cite{seleson2015concurrent}. The one-dimensional  studies 
include a point load case, where 
the sensitivity of the numerical results on the model parameters, including the size of the nonlocal sub-domain where the point load is applied, the size of the blending region, the choice of blending function, and the horizon size, is investigated. Moreover, additional studies are conducted to quantify the effect of the location of the point load with respect to the coupling configuration as well as the speedup attained by employing a coarse discretization in the local sub-domain. 
%
Two-dimensional simulations include patch-test and point load examples as well as the deformation of a square plate with a rectangular slit  at its center.




\subsection{Splice method}\label{subsec:splice}
The splice method was originally presented for state-based peridynamics in~\cite{silling2015variable} as a means to couple two peridynamic models with different horizons.  
Then, it was applied to LtN coupling
by 
taking the horizon in one of the two peridynamic models to zero, so that the model can be effectively replaced by a classical local model. A methodology resembling the splice approach, formulated instead at the discrete level, was proposed to couple a discretized bond-based peridynamic model with classical meshless methods in~\cite{shojaei2016coupled,shojaei2016coupling,shojaei2017coupling} and with classical finite elements in~\cite{galvanetto2016effective,zaccariotto2018coupling,zaccariotto2017enhanced}. The latter was extended to  state-based peridynamics in~\cite{ni2019coupling}. Conceptually, the splice approach is probably the simplest LtN coupling method because each material point is described with a reference fully local or fully nonlocal model without the need of special coupling techniques.

\subsubsection{Mathematical formulation}  
Let a domain $\omg$ be decomposed in two non-overlapping sub-domains $\omgn$ and $\omg_l$, such that $\omg = \omgn \cup \omgl$. These two sub-domains are connected by the interface $\Gamma_v$. 
This configuration resembles the one used in the partitioned procedure for the non-overlapping case and is illustrated in Figure \ref{fig:blended-domains} (except that the same nonlocal formulation is employed in the whole nonlocal sub-domain, which includes both $\Omega_{nl}$ and $\Omega_t$). 
%

We define the splice coupling operator as follows (see \eqref{eq: PD operator state-based limit} and \eqref{eq: PD operator state-based}):
\begin{equation}\label{eq: splice operator}
\mcL^{\rm splice}\ub(\xb)  :=  
\left\{ \begin{array}{cc}
\displaystyle
 \int_{ B_\horizon(\xb)} \left\{
\statevecT[\xb]\langle \xbp - \xb\rangle  - \statevecT[\xbp]\langle \xb - \xbp \rangle
  \right\}d\xbp 
& \quad \xb \in \omg_{nl},\\[0.15in]
 \nabla\cdot \nub^0(\xb) 
& \quad \xb \in \omg_l. 
\end{array}
\right.
\end{equation}
The basic idea is that each material point is represented by either a fully local or a fully nonlocal model. Near the interface $\Gamma_v$, a point in $\omg_{nl}$ interacts with points in $\omg_l$. While such points belong to a local sub-domain, from the perspective of a point in $\omg_{nl}$ the interaction is described by a fully nonlocal model. A similar, but reversed, situation may occur after discretization for nodes in $\omg_l$ that could interact with some nodes in $\omg_{nl}$.

\paragraph{Properties}
\begin{itemize}
\item 
Because each material point in the splice method is described by either a fully local or a fully nonlocal model, and those models are consistent up to third-order polynomials, the splice model passes a linear, quadratic, and cubic patch tests.
\item 
The asymptotic compatibility property of the splice model is inherited from the asymptotic compatibility of the reference nonlocal model with an order of convergence of $\mcO(\horizon^2)$.

\item The numerical implementation of the splice method is generally non-intrusive and only requires passing information about the deformation between the local and nonlocal sub-domains across $\Gamma_v$.

%
\end{itemize}

\paragraph{The time-dependent problem.} The splice method is directly applicable to time-dependent problems, as demonstrated in, e.g.,~\cite{silling2015variable}.


\subsubsection{Applications and results}

\begin{figure}[htbp]
\begin{center}
\subfigure[Linear Patch Test]{
\includegraphics[scale = 0.3]{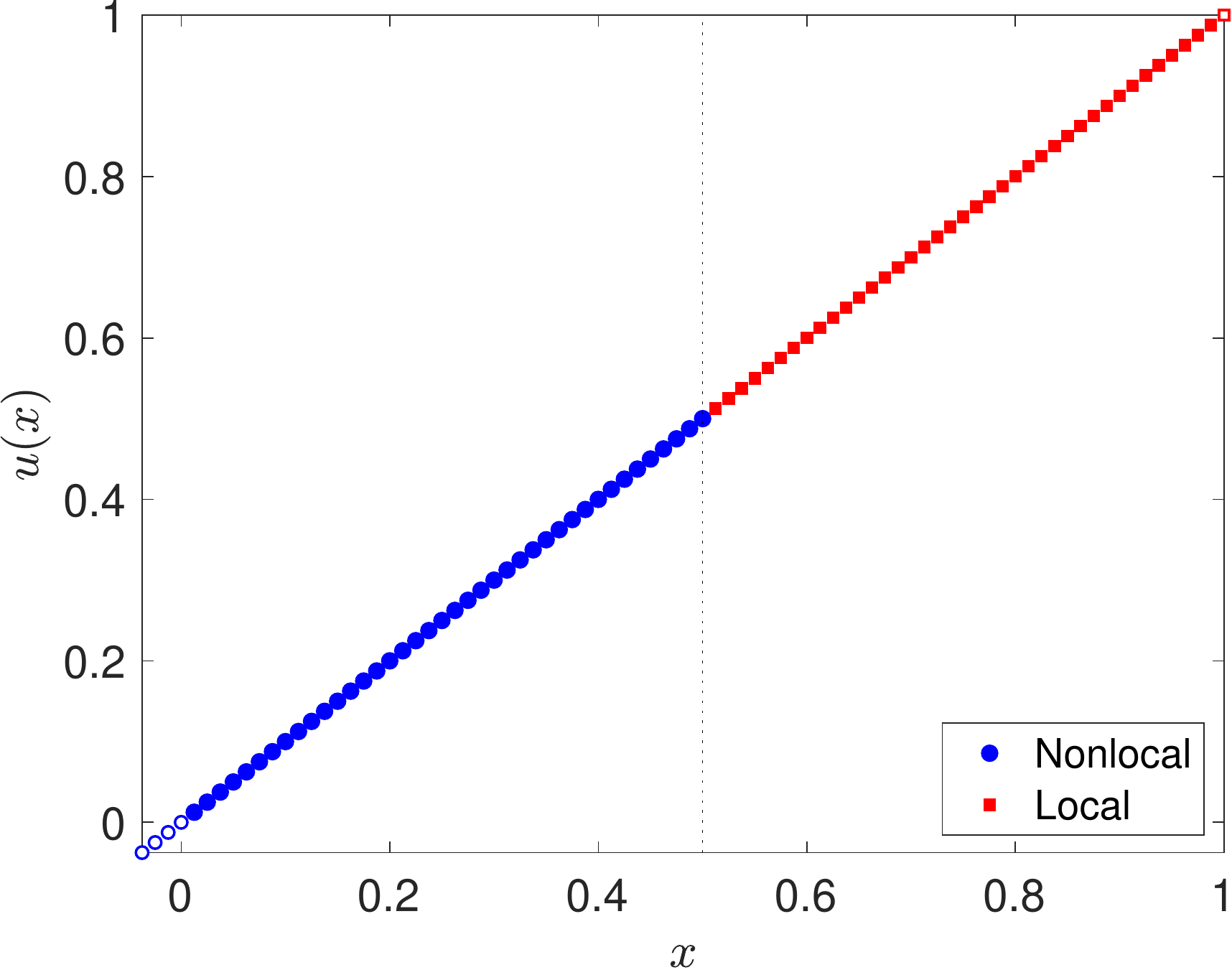}}
\qquad
\subfigure[Quadratic Patch Test]{
\includegraphics[scale = 0.3]{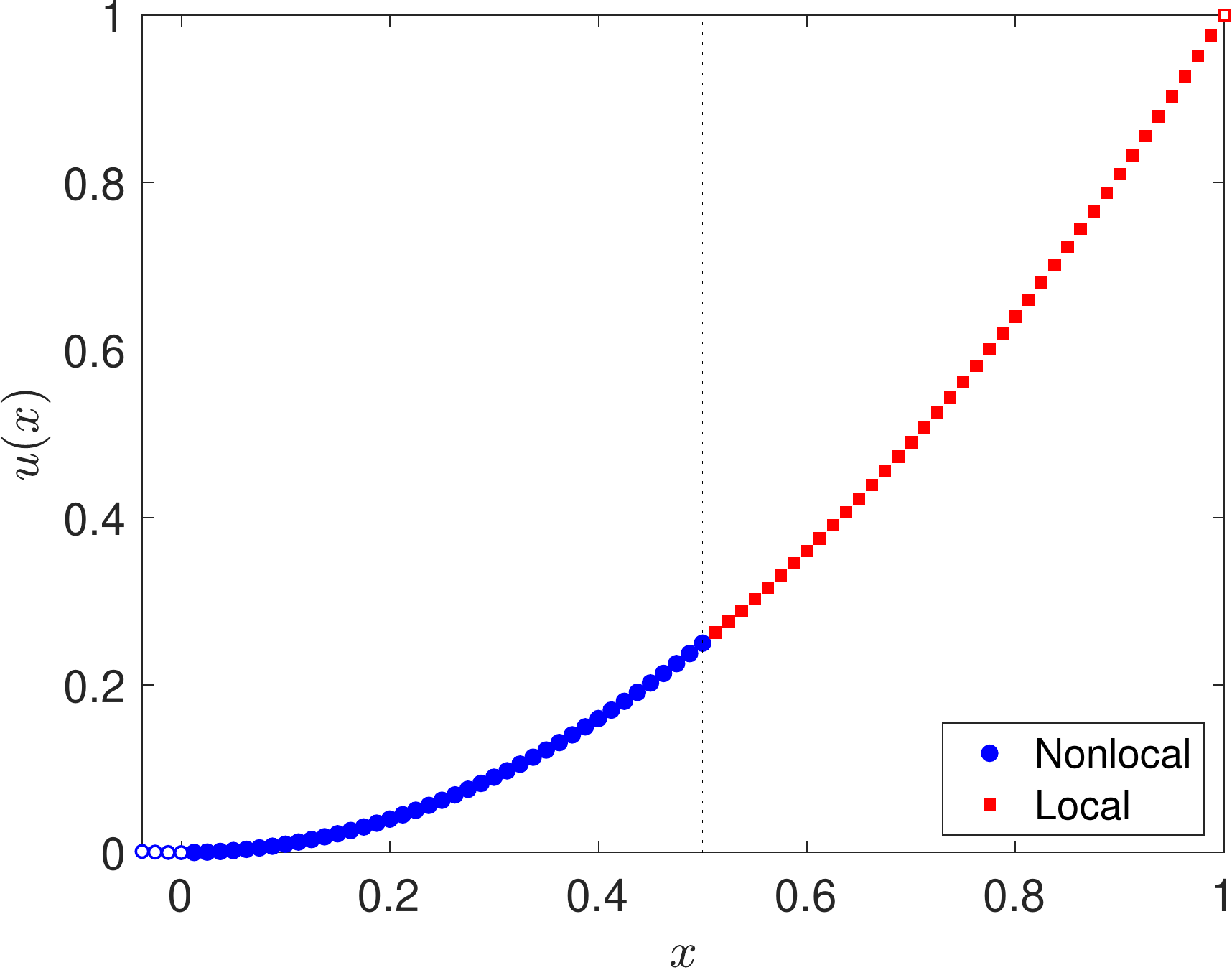}}
\vspace{-0.14 in}
\caption{Patch tests for the  splice method: linear (left) and quadratic (right). Nodes in the nonlocal sub-domain are represented by blue filled circles, whereas nodes in the left nonlocal boundary are empty blue circles. Nodes in the local sub-domain are represented by red filled squares, whereas the node in the right local boundary is an empty red square.}
\label{Fig:Patch Tests Splice Model}
\end{center}
\end{figure}

Consider a one-dimensional domain $\ooomg = [-\horizon,1]$ with $\horizon = 0.05$ decomposed into $\ooomg =\omg_p \cup \omg_{nl} \cup \omg_l $, where  $\omg_p = (-0.05, 0)$, $\omg_{nl} = (0,0.5)$, and $\omgl=(0.5,1)$. 
%
%
We solve both a linear patch test with analytical solution $u(x) = x$ and a quadratic patch test with analytical solution $u(x) = x^2$ . We employ a linear bond-based peridynamic model in $\omg_{nl}$ discretized with a meshfree method and a classical linear elasticity model in  $\omgl$ discretized with a finite difference method. 
%
%
A uniform grid with grid spacing $\Delta x = \horizon / 4$ is utilized. 
%
%
Nonlocal Dirichlet boundary conditions are imposed in the left nonlocal boundary for the peridynamic model, whereas a local Dirichlet boundary condition is imposed on the right local boundary for the local model. 
The results are presented in Figure~\ref{Fig:Patch Tests Splice Model} and demonstrate that the splice model passes the linear and quadratic patch tests. 
A computation of the first derivative of the solution for the linear patch test and the second derivative of the solution for the quadratic patch test indicates that these tests are passed to machine precision.

Several numerical studies concerning the splice method have been carried out in the above provided references. First, in~\cite{silling2015variable}, the performance of the splice method is studied for a spall initiated by the impact of two brittle elastic plates. 
In \cite{galvanetto2016effective}, static numerical tests are performed, including one-dimensional linear and quadratic patch tests as well as three simulations in two dimensions, which concern the deformation of a plate undergoing rigid body motion, horizontal stretch, and shear. 
%
%
Wave and crack propagation examples were discussed in \cite{zaccariotto2018coupling} in one, two, and three dimensions. 
Additional dynamic examples appear in  \cite{shojaei2016coupled} for a cantilever beam subjected to a periodic excitation as well as crack propagation and branching in a pre-cracked plate under traction. 
%
%
%
Advanced numerical examples including brittle failure analysis are given in \cite{ni2019coupling}, including the Wedge-splitting test and the Brokenshire torsion experiment.
Finally, in \cite{shojaei2017coupling} a switching nodal technique that can be used for adaptivity is presented. 
%


\section{Varying Horizon approaches}\label{sec:VH}

The concept of varying the horizon in a nonlocal model for LtN coupling purposes is motivated by the convergence of nonlocal models to local counterparts in the limit as $\horizon \to 0$, under suitable regularity assumptions, as discussed in Sections \ref{subsubsec:diffusion} and~\ref{subsubsec:mechanics}.  
Earlier studies with variable horizon for adaptive refinement in peridynamics are presented in~\cite{bobaru2011adaptive,bobaru2009convergence}. 
In the context of nonlocal diffusion, a spatially varying-horizon formulation for analysis of interface problems is introduced in~\cite{seleson2013interface} and specialized to a sharp transition between nonlocal and local sub-domains. 
The concept of smoothly varying the horizon is presented and analyzed in~\cite{silling2015variable} and demonstrated for the coupling of two peridynamic models with different horizons. Later on, the shrinking horizon approach in  nonlocal diffusion is analyzed for LtN coupling in  \cite{TTD19,TD17trace}. 
As demonstrated below, the shrinking horizon method is not linearly patch-test consistent. To overcome this issue, while allowing to spatially vary the horizon, the partial stress method is proposed in~\cite{silling2015variable}. 

In this section, we begin by discussing the shrinking horizon method based on \cite{silling2015variable,TTD19} in Section \ref{subsec:shrink}, followed by a description of the partial stress method presented in \cite{silling2015variable} in Section \ref{subsec:partial-stress}.

\subsection{Shrinking horizon}
\label{subsec:shrink}

In \cite{silling2015variable}, the idea of spatially varying the horizon is discussed. Specifically, the conception of material homogeneity under a change in horizon is presented, resulting in the concept of \textit{variable scale homogeneous} body. This provides an appropriate rescaling of the peridynamic force vector state, so that the corresponding strain energy density remains invariant with respect to changes in the horizon, under uniform deformations. 
In the related work in \cite{TTD19,TD17trace}, the authors discuss the validity of the nonlocal diffusion models with a shrinking 
 horizon as the material points approach an interface, where the nonlocal models were localized and connected to local models on the other side, effectively resulting in an LtN coupling. The coercivity and the trace space of the corresponding energy functionals were discussed to guarantee the well-posedness of such nonlocal models with localization on 
 {the interface}. 
%
 The coupled model could be solved as a whole system, or be solved using a classical non-overlapping domain decomposition method where the nonlocal and local models could be treated separately with suitable {transmission} conditions on the 
 {interface}. 
 The shrinking horizon method is an energy-based VH approach. This method 
 does not pass the patch tests; however, the errors can be controlled by choosing the horizon function properly. 
 Below, for brevity, we mainly discuss the formulation in the context of nonlocal diffusion, but make connections, where  appropriate, to the analogue peridynamic formulation. 

\subsubsection{Mathematical formulation}
Let $\Omega_{nl}$ and $\Omega_{l}$ be two open domains in $\R^n$ that satisfies $\overline{\Omega_{nl}}\cap\overline{\Omega_{l}}=\Gamma\subset\R^{n-1}$, as illustrated in Figure \ref{fig:variable-domain}.
The nonlocal energy space $\mcS^{\rm{ND}}(\ooomg_{nl})$ is then tied to a kernel function $\gamma$ with heterogeneous localization on the boundary $\Gamma$. More specifically, we assume that $\gamma$ is given by 
\begin{equation} \label{eq:shrink-kernel}
    \gamma(\xb, \xbp) =  
    {B_{\delta(\xb)}(\xbp-\xb)\,k(\xb, \xbp)},
\end{equation}
similarly to~\eqref{eq:loc-kernel}, except that 
$\delta(\xb)$ is now a spatially-dependent function that 
decays to zero as $\xb$ approaches $\Gamma$.
Using a properly scaled function kernel $k(\xb, \xbp)$, such as 
\[
k(\xb, \xbp)= \frac{{6}}{|\partial B_1| (\delta(\xb))^{n+2}} \quad 
\text{or } \quad
k(\xb, \xbp)= \frac{{4}}{|\partial B_1| (\delta(\xb))^{n+1}|\xb-\xbp|}\,,
\]
{where $|\partial B_1|$ stands for the area of the surface of the unit ball in $\R^n$,} we can still get $\Delta u$ as the limit operator of $\mcL^{\rm{ND}} u$ by taking $\delta(\xb) \to 0$ for every $\xb\in \Omega_{nl}$.
A simple choice of the horizon function $\delta(\xb)$ is given by 
\begin{equation}\label{deltax:linear}
    \delta(\xb) = \min(\delta, \text{dist}(\xb, \Gamma))\,,
\end{equation}
with $\delta$ being the maximum value of $\delta(\xb)$. This choice of the horizon function is depicted in Figure \ref{fig:variable-domain}. In one dimension, the horizon function $\delta(\xb)$ given by \eqref{deltax:linear} is a piece-wise linear function. However, numerical examples in Section \ref{subsubsec:shrink_num} show that better horizon function may be used in order to attain optimal order of convergence to the local limit as $\max_{\xb}\delta(\xb)\to0$.

\begin{remark}
A related scaling of the peridynamic force vector state in \cite{silling2015variable} is given by
\begin{equation}
\widehat{\underline{\bf T}}(\underline{\bf Y}[\xb],\xb)\langle \xib\rangle 
=
\frac{1}{(\horizon(\xb))^{1+n}}
\widehat{\underline{\bf T}}_1(\underline{\bf Y}_1[\xb])
\left<\frac{\xib}{\horizon(\xb)}\right>, 
\end{equation}
where $n$ is the spatial dimension, $\widehat{\underline{\bf T}}_1$ is a reference material model, $\underline{\bf Y}$ is the deformation state, and $\underline{\bf Y}_1$ is the reference deformation state. The hat notation is used to explicitly indicate dependence on the deformation. 
\end{remark}

With the localization of the nonlocal interactions at the boundary $\Gamma$, it is shown in \cite{TD17trace} that the nonlocal energy space  $\mcS^{\rm{ND}}(\ooomg_{nl})$ has $H^{1/2}(\Gamma)$ as the trace space on $\Gamma$, which is exactly the trace space of $H^1$ functions. As a result of the trace theorem, we can define the combined energy space
\[
\mcW(\ooomg)=\{ u\in \mcS(\ooomg_{nl})\cap H^1(\Omega_l)\, :\, u_-=u_+ \text{ on } \Gamma, u|_{\omgb} =0  \}\,,
\]
where $u_-(\xb)$ and $u_+(\xb)$ are defined as $\lim_{\yb\to \xb, \yb\in\Omega_{nl}}u(\yb)$ and $\lim_{\yb\to \xb, \yb\in\Omega_{l}}u(\yb)$, respectively. The total energy is a combination of the nonlocal and local parts given by
\[
\begin{split}
E(u, f)=\frac{1}{4}\int_{\ooomg_{nl}}\int_{\ooomg_{nl}} \gamma(\xb,\xbp) ((\mcD^\ast u)(\xb,\xbp))^2 d\xb d\xbp & +\frac{1}{2} \int_{\Omega_l}|\nabla u(\xb)|^2d\xb \\
&- \int_{\Omega} f(\xb)u(\xb) d\xb \,,
\end{split}
\]
for any $u\in\mcW(\ooomg)$. The well-posedness of the coupled problem is guaranteed by the extension of the nonlocal Poincar\'e inequality to nonlocal space with variable-horization function $\gamma$ given by \eqref{eq:shrink-kernel}. 

\paragraph{Properties}
\begin{itemize}
\item The shrinking horizon method is mathematically well-posed and energy stable for diffusion problems on general domains in all dimensions. 
\item There is no overlapping region 
{between nonlocal and local sub-domains}.
Moreover, since the  nonlocal and local energy functionals have the same trace space on the interface, one can use all the {classical} non-overlapping domain-decomposition methods for solving the coupled problem. 
\item 
The method does not pass the patch-tests. However, ghost forces can be controlled by using a slowly varying horizon function. 
\item 
The order of convergence of the coupled problem to the local problem as $\delta\to0$ depends on the choice of the horizon function. For a piece-wise linear horizon function, the solution convergences at a rate of $O(\delta)$. With a slowly varying horizon function, the optimal order $O(\delta^2)$ could be achieved. 
\item 
The total energy is equivalent to either fully nonlocal or fully local energy up to linear functions.  
\end{itemize}

\subsubsection{{Applications and results}} \label{subsubsec:shrink_num}
The shrinking horizon approach produces well-posed coupled models for a general horizon function $\delta(\xb)$ that gets localized at $\Gamma$, as shown in \cite{TTD19,tian_thesis}. 
However, the particular choices of the horizon function affect the convergence rate of the solutions as $\delta(=\max_{\xb} \delta(\xb)) \to 0$. In \cite{TTD19},
one-dimensional experiments are performed to illustrate that the solutions converge only at the order $O(\delta)$ when using a piece-wise linear horizon function, as the one given by \eqref{deltax:linear}. The authors provide two remedies for increasing the order of the convergence. The first approach is to use a specific auxiliary function, and the second is to use a smooth and slowly varying horizon function.
Here we only discuss the second approach; this approach was also discussed in \cite{silling2015variable} as a means to reduce ghost forces. 
\begin{figure}[H]
\begin{center}
\subfigure[Horizon functions]{
\includegraphics[width =5.5cm]{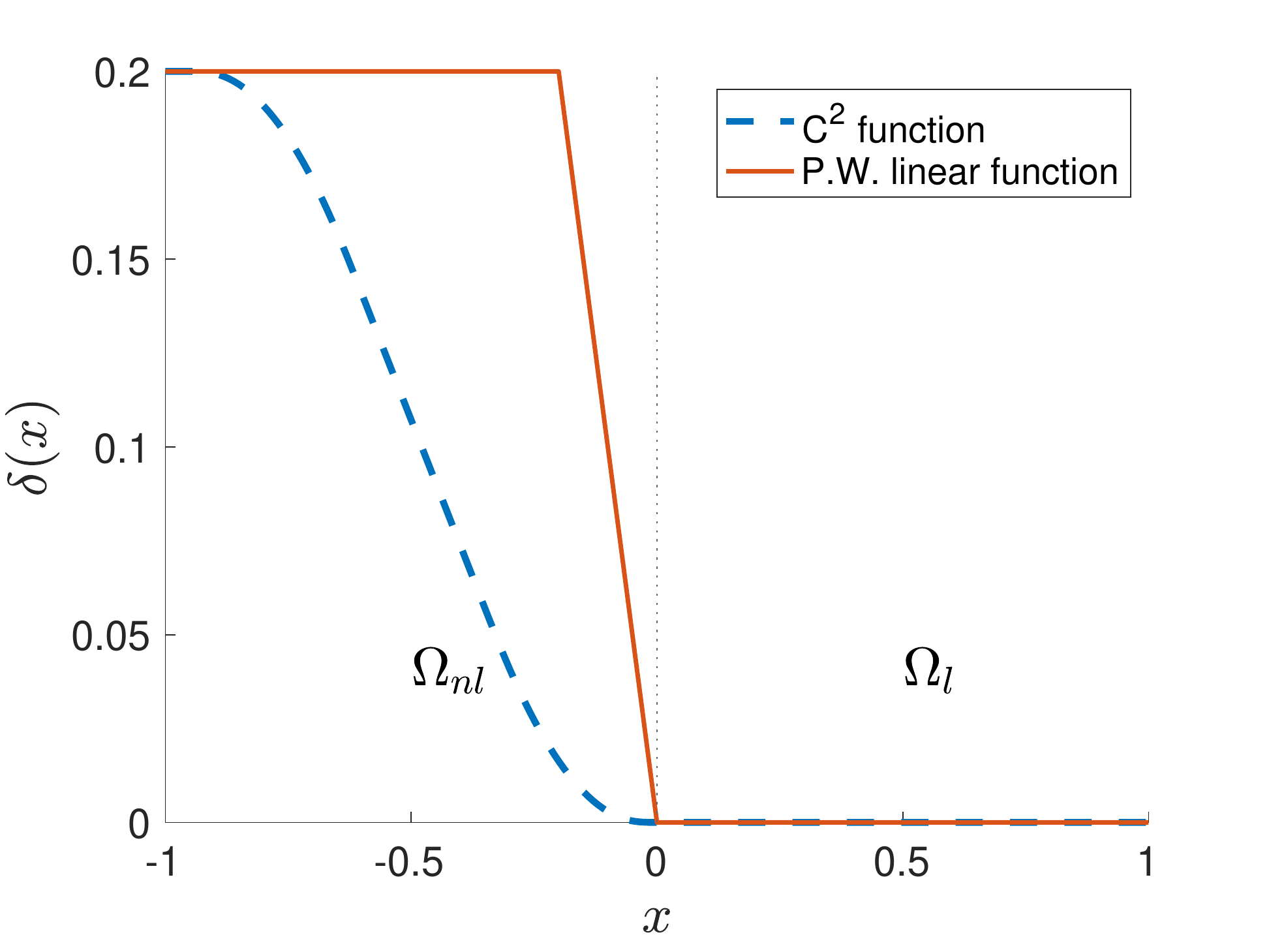}}
\qquad
\subfigure[Ghost forces]{
\includegraphics[width= 5.5cm]{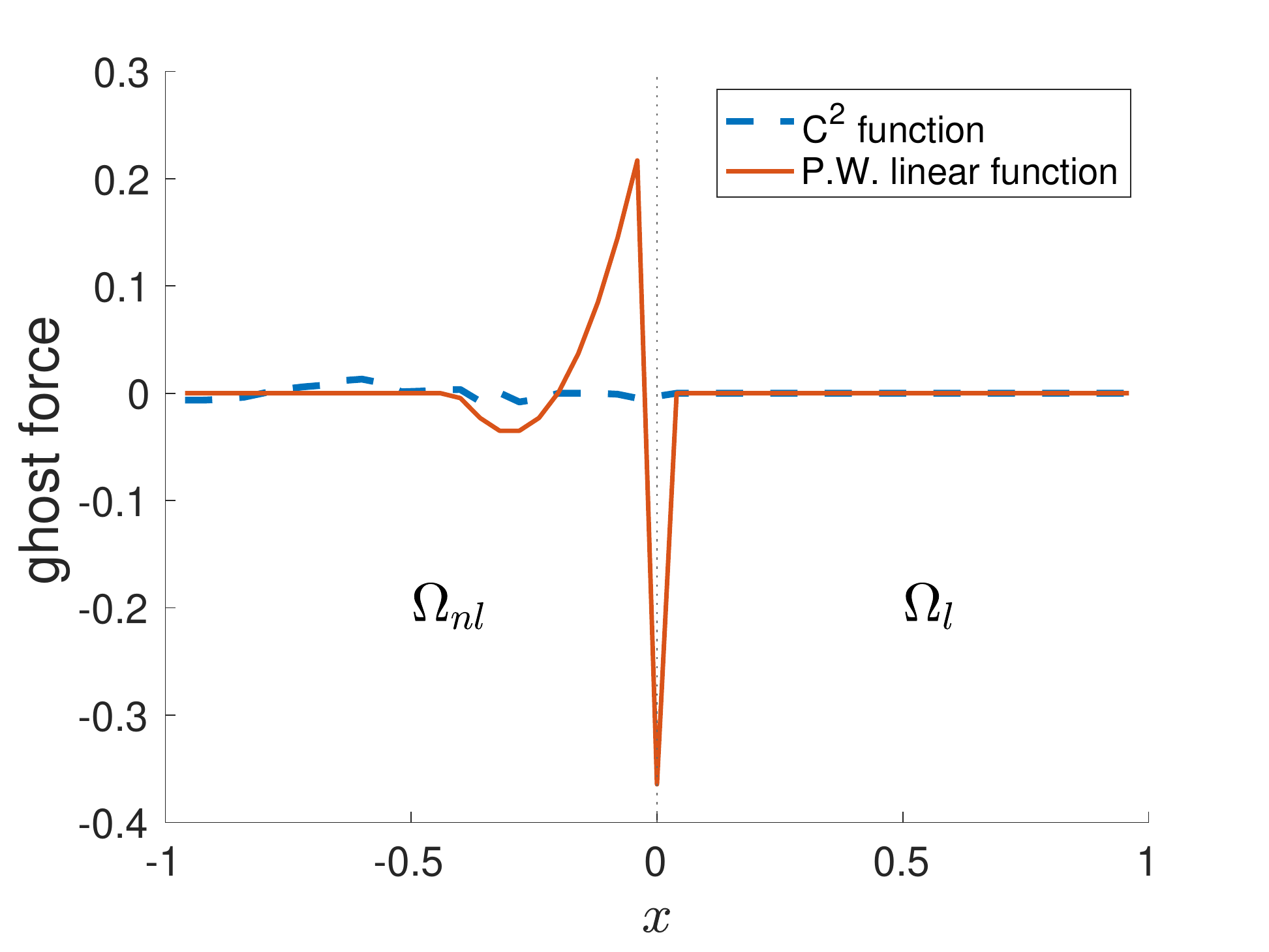}}
\vspace{-0.14 in}
\caption{Left: examples of $\delta(x)$ in one dimension. The orange solid line represents a piece-wise linear horizon function and the blue dashed line represents a $C^2$ horizon function. Right: the ghost forces under linear displacement using the piece-wise linear horizon function and the $C^2$ horizon function given by the left plot. }
\label{fig:horizon}
\end{center}
\end{figure}
\begin{figure}[H]
\begin{center}
\subfigure[Linear Patch Test]{
\includegraphics[width =5.5cm]{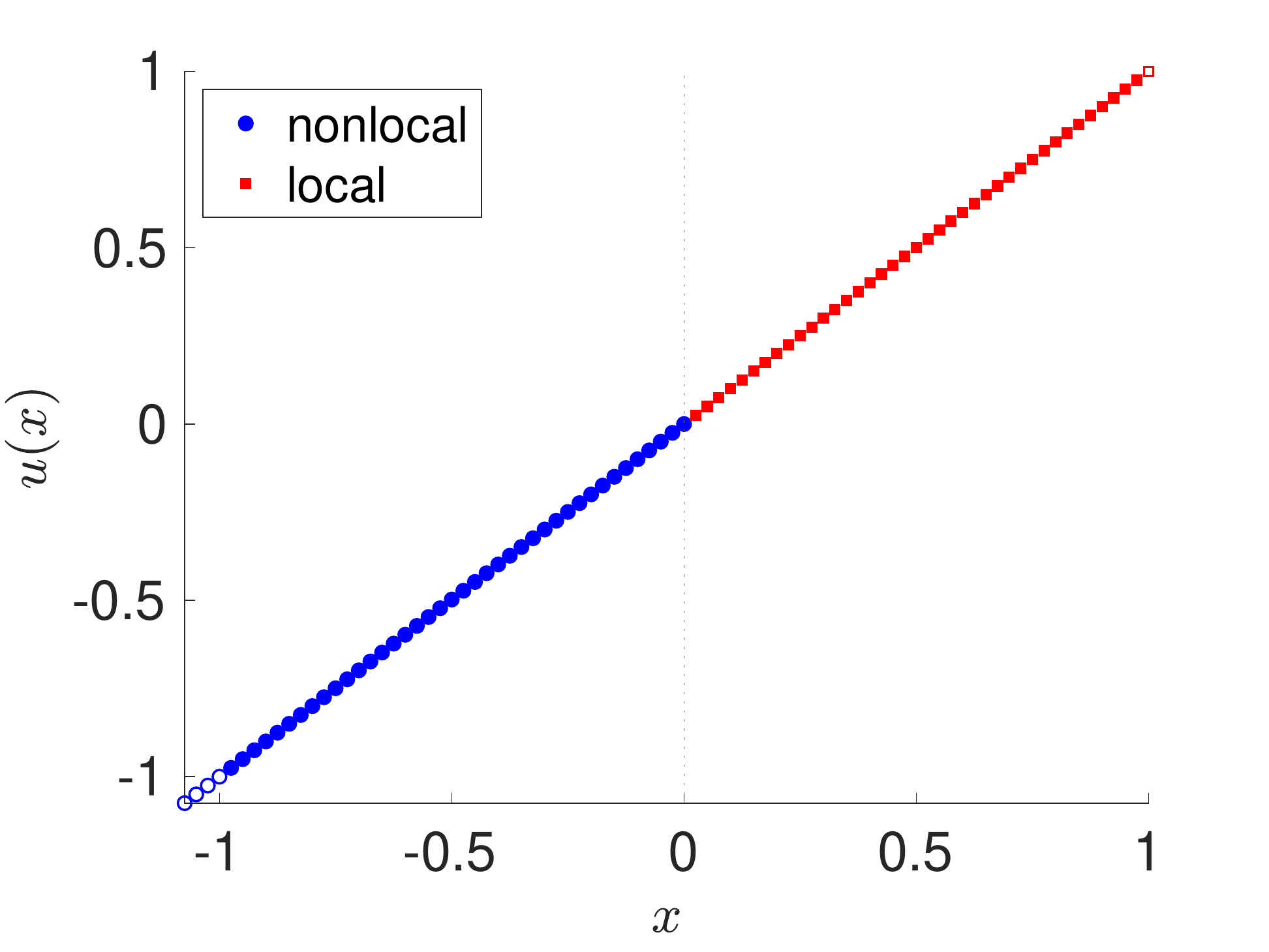}}
\qquad
\subfigure[Quadratic Patch Test]{
\includegraphics[width= 5.5cm]{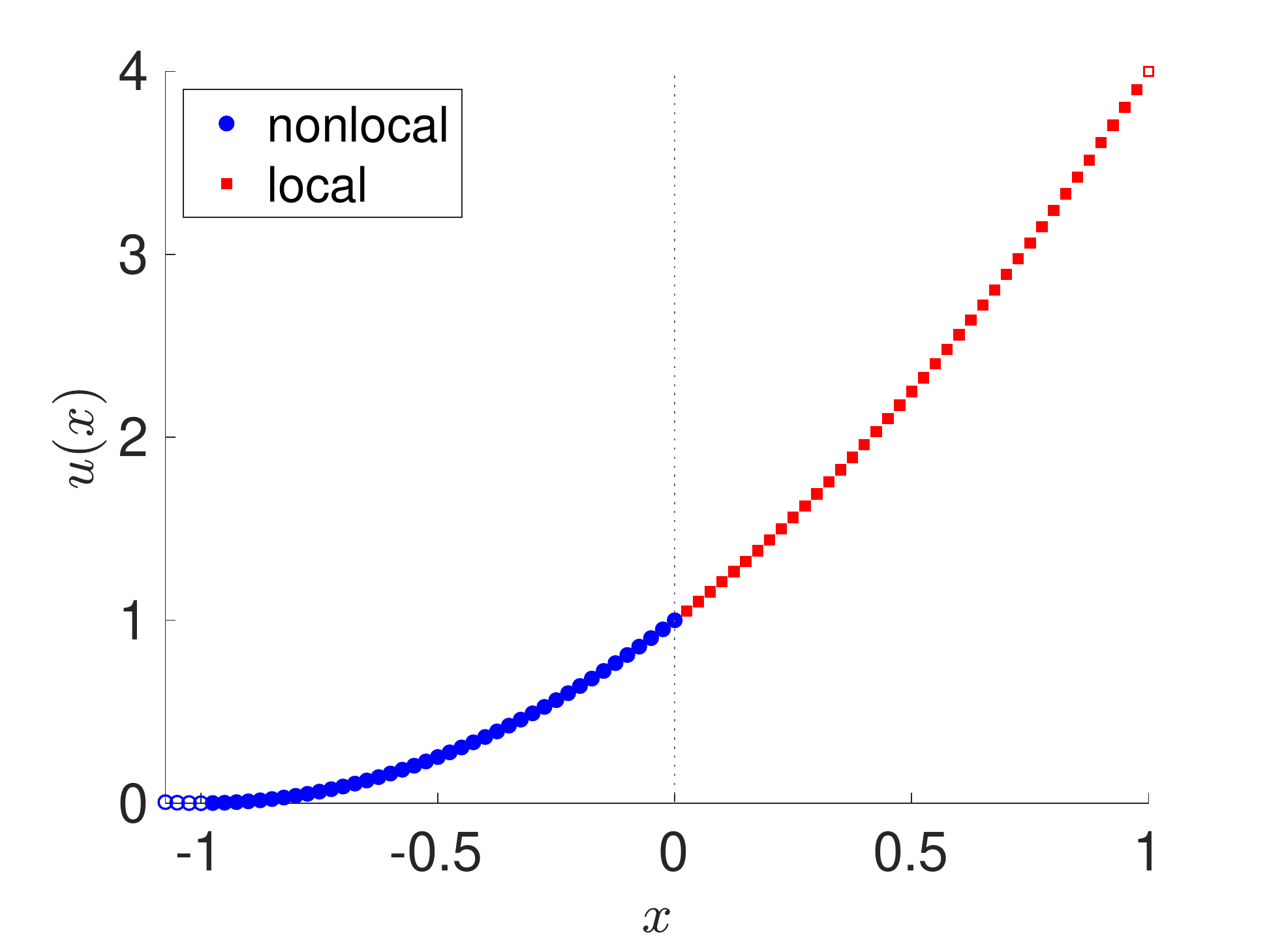}}
\vspace{-0.14 in}
\caption{Patch tests for the shrinking horizon method with $C^2$ horizon function $\delta(x)$: linear (left) and quadratic (right). 
%
Nodes in the nonlocal sub-domain are represented by blue filled circles, whereas nodes in the left nonlocal boundary are empty blue circles. Nodes in the local sub-domain are represented by red filled squares, whereas the node in the right local boundary is an empty red square.
%
}
\label{fig:patch_tests_shrink}
\end{center}
\end{figure}
The left plot of Figure \ref{fig:horizon} shows two choices of the horizon function $\delta(x)$ in one dimension, with one being a piece-wise linear function and another being a $C^2$ function. 
Since the model does not pass the linear patch test exactly, it generates ghost forces when the solution is a linear profile $u^{\rm lin}$ (see Definition \ref{def:patch-test}). The right plot of Figure \ref{fig:horizon} shows the ghost forces over the domain $\Omega=(-1,1)$ using each of the horizon functions. 
It is clear that the piece-wise linear horizon function generates a large magnitude of ghost forces around the interface between $\Omega_{nl}=(-1,0)$ and $\Omega_l=(0,1)$, while the ghost forces are under control by using the $C^2$ horizon function. Moreover, it is observed in \cite{TTD19} that the ghost forces converge to zero everywhere as $\delta\to0$ if the $C^2$ horizon function is used. 
Figure \ref{fig:patch_tests_shrink} further shows the linear and quadratic patch tests by using the $C^2$ horizon function. In these tests, the largest horizon is $\delta =0.1$  and the spatial mesh size is $\Delta x=0.025$. 
The FEM 
with piece-wise linear basis functions is used for computing the numerical solution, since it is shown to be an asymptotically compatible scheme for nonlocal variational problems \cite{Tian2014}. Although the coupled model does not pass the patch tests in theory, Figure \ref{fig:patch_tests_shrink} shows that by using the $C^2$ horizon function, the patch-test consistency could be preserved approximately. 

Furthermore, numerical examples in \cite{TTD19} show that by using the $C^2$ horizon function, the optimal order of convergence to the local limits could be achieved. Solutions of the coupled problems converge to the solutions of the local problems in the $L^2$ norm at the rate $O(\delta^2)$, and the numerical derivative of the solutions converge at the rate $O(\delta)$. In contrast, using the piece-wise linear horizon function, one can only observe first order convergence in the solutions. 

\subsection{Partial stress method}
\label{subsec:partial-stress}

As discussed in Section~\ref{subsec:shrink}, a shrinking horizon approach is not patch-test consistent, even though deviations from the patch test can be controlled by the regularity of the horizon function. In~\cite{silling2015variable}, this consideration led to the development of an alternative strategy to spatially vary the horizon, referred to as the partial stress method.

The proposition of the partial stress approach, which relates to Figure~\ref{fig:blended-domains}, for LtN coupling is to reformulate the operator~\eqref{eq: PD operator bond-based bulk} in the transition region connecting local and nonlocal sub-domains, 
in a way that spatially varying the horizon in that region does not give rise to ghost forces under uniform deformations. 

\subsubsection{Mathematical formulation}
This method introduces a new tensor-valued function called the {\it partial stress} tensor, 
\begin{equation}\label{eq: partial stress tensor}
\nub^{\rm ps}(\xb) := \int_{B_\horizon({\bf 0})} \statevecT[\xb]\langle \xib\rangle  \otimes \xib  \, d\xib,
\end{equation}
where $\statevecT$ is the force vector state from \eqref{eq: PD operator state-based}, 
and defines a corresponding {\it partial internal force density},
\begin{equation}\label{eq: PS operator bond-based bulk}
\mcL^{\rm ps} \ub (\xb) := \nabla \cdot \nub^{\rm ps}(\xb).
\end{equation}
It is important to note that the partial stress tensor coincides with the collapse stress tensor, $\nub^0$ ({\it cf.}~\eqref{eq: PD operator state-based limit}), 
in the case of a uniform deformation of a homogeneous body, which is characterized by a constant $\horizon$.

We now refer to Figure \ref{fig:blended-domains}: the domain $\omg$ is decomposed into {three} disjoint sub-domains: {$\omg=\omg_{nl}\cup\Omega_t\cup\omgl$}, i.e.,  the nonlocal sub-domain, the transition region, and the local sub-domain. The partial stress coupled problem is given by 

\begin{equation} \label{eq: PS static problem}
\left\{
\begin{aligned}
- \int_{ B_\horizon(\xb)} \left\{
\statevecT[\xb]\langle \xbp - \xb\rangle  - \statevecT[\xbp]\langle \xb - \xbp \rangle
  \right\}d\xbp
= \bb(\xb) & \quad\xb\in\omgn, \\[2mm]
-\nabla \cdot \nub^{\rm ps}(\xb) = \bb(\xb) & \quad\xb\in\omgt, \\[2mm]
-\nabla \cdot \nub^0(\xb)= \bb(\xb) & \quad\xb\in\omgl. 
\end{aligned}\right.
\end{equation}

\paragraph{Properties.} Not many properties have been discussed for the partial stress method in the literature. Here, we simply summarize two properties from~\cite{silling2015variable}: 
\begin{itemize}

\item For a uniform deformation of a homogeneous body,
\begin{equation*}
\nub^{\rm ps} \equiv\nub^0  
\qquad \mbox{and} \qquad
\nabla \cdot \nub^{\rm ps} \equiv\nabla \cdot \nub^{0}\equiv {\bf 0}.
\end{equation*}

\item The partial stress tensor \eqref{eq: partial stress tensor} and  partial internal force density \eqref{eq: PS operator bond-based bulk} converge, under suitable regularity assumptions, to the corresponding classical local counterparts, 
in the limit as $\horizon \to 0$, as $O(\horizon)$, so the method is asymptotically compatible. 
%
%
Consequently, the partial stress tensor \eqref{eq: partial stress tensor} and  partial internal force density \eqref{eq: PS operator bond-based bulk} are also compatible to the fully nonlocal counterparts with a $\horizon$-order difference.

%
\end{itemize}


\paragraph{The time-dependent problem.} The partial stress method was, in fact, only demonstrated in a dynamic setting in~\cite{silling2015variable}, where in conjunction with the  splice method was applied to the study of a spall initiated by the impact of two brittle elastic plates. 


\subsubsection{{Applications and results}} \label{subsubsec:partial_stress_num}

As described above, the partial stress method was only applied in~\cite{silling2015variable} for a dynamic problem. We therefore omit the details here and refer the reader to that work.  

\bigskip

\section{Conclusions}\label{sec:conclusion}
This paper presents a review of the state-of-the-art of LtN coupling for nonlocal diffusion and nonlocal mechanics, specifically peridynamics, and provides a classification of different LtN coupling approaches (see Figure \ref{fig:overview-chart}). Following a description of various coupling configurations and a highlight of desired properties of a general LtN coupling strategy, we report different LtN coupling methods from the literature. For each method, we briefly present its mathematical formulation and properties, and we discuss relevant applications and numerical results. 

We observe that, while many features and challenges are shared by all methods, there exist some significant differences in their formulation and implementation. 
For instance, we find that even though a LtN coupling configuration can generally be divided into a local sub-domain, a transition region, and a nonlocal sub-domain (as illustrated in Figure \ref{fig:general-domains}), each coupling method treats the transition region in its own specific way. Some methods overlap local and nonlocal descriptions, some employ a hybrid representation, some reduce the transition region to a sharp interface, and some utilize a variable horizon with or without changing the nonlocal operator. 
This variation in the treatment of the transition region has both analytical and numerical implications. 
For instance, an overlapping approach is normally non-intrusive, whereas a hybrid technique is typically intrusive; the variable horizon method may or may not be intrusive, depending on the available nonlocal implementation. 
Another important property largely emphasized in this review is the ability of a coupling method to pass the patch test. Some coupling methods do pass it exactly for up to certain polynomial order, whereas others only pass it approximately; a linear patch test is the most popular one. 
Finally, we recognize two major formulations for LtN coupling, energy-based and force-based. The former provides a natural setting to impose energy preservation; however, because such formulation requires energy minimization, it is more native to static problems and the extension to dynamics settings may not be practical. 
%
On the other hand, while force-based approaches normally equally apply to static and dynamics problems, they not always carry a well-defined energy functional. 
In Table \ref{table:coupling_summary}, we outline the methods discussed in this review, indicating relevant references and sections, and summarize some of these properties.

We conclude by stating that the goal of this review is not to provide a preferred way to perform LtN coupling, but rather to broaden the perspective of the reader that can use this review as a guide for selecting the most appropriate method based on the characteristics 
%
%
of the problem at hand, available discretization methods, and accesible data.

\begin{table}[H]
\centering
\small
\label{table:coupling_summary}
\begin{tabular}{cccccc}
\hline
{\bf Method}  & \shortstack{\bf References} & \shortstack{\bf Section}&  \shortstack{\bf Transition}& \shortstack{\bf Linear Patch Test}& \shortstack{\bf Formulation}\\
\hline
\shortstack{Optimization-based} & 
\shortstack{\cite{Bochev_14_INPROC,Delia2019,Bochev_16b_CAMWA}}
&\shortstack{ \ref{subsec:OBM}} & \shortstack{Overlap} &
\shortstack{Exact}& Force-based\\
\hline
\shortstack{Partitioned}
&\shortstack{\cite{you2019coupling,yu2018partitioned}} &\shortstack{ \ref{subsec:Partition_Robin}}  &  \shortstack{Overlap} &
\shortstack{Exact$^*$} &{Force-based}\\[-0.02in]
& &  &  \shortstack{or Sharp} & &\\
\hline
{Arlequin}&\shortstack{\cite{HanLubineau2012,ArlequinWang2019}}&\shortstack{ \ref{subsec:Arlequin}} & Hybrid
&\shortstack{Approximate} & Energy-based\\
\hline
{Morphing}&\shortstack{\cite{han2016morphing,lubineau2012morphing}} 
&\shortstack{ \ref{subsec:morphing}}
& Hybrid
& Approximate &{Force-based}\\
\hline
\shortstack{Quasi-nonlocal}&\shortstack{\cite{DuLiLuTian2018,XHLiLu2017}}&\shortstack{ \ref{subsec:quasinonlocal}} & Hybrid
& {Exact}&{Energy-based}\\
\hline
\shortstack{Blending}&\shortstack{\cite{Seleson2013CMS,seleson2015concurrent}}& \shortstack{\ref{subsec:blending}} & Hybrid
&{Exact}&{Force-based}\\
\hline
\shortstack{Splice}&\shortstack{\cite{silling2015variable,galvanetto2016effective}}& \shortstack{\ref{subsec:splice}} 
& Sharp
&{Exact}&{Force-based}\\
\hline
\shortstack{Shrinking horizon}&\shortstack{
\cite{silling2015variable,TTD19}}&\shortstack{ \ref{subsec:shrink}}
& Variable horizon
&{Approximate}&{Energy-based}\\
\hline
\shortstack{Partial stress}&\shortstack{
\cite{silling2015variable}}& \shortstack{\ref{subsec:partial-stress}} & Variable horizon
&{Exact}&{Force-based}\\[-0.02in]
& &  &  \shortstack{with partial stress} & &\\
\hline
\end{tabular}
\caption{Summary of LtN coupling methods. The $^*$ in the table indicates that the patch-test consistency for the partitioned procedure is exact up to certain conditions (see Section \ref{subsec:Partition_Robin}). 
}
\end{table}

\normalsize

\bibliographystyle{spmpsci}
\bibliography{reference}

\end{document}